\documentclass[12pt]{amsart} %,3p]{elsarticle}
\usepackage{amsmath}
\usepackage{mleftright}

\usepackage{graphicx}
\usepackage[all]{xy}

\usepackage{hyperref}

\newcommand {\fbgl}{{\mathfrak{bgl}}}
\newcommand {\fbr}{{\mathfrak{br}}}

\newcommand {\fbrj}{{\mathfrak{brj}}}
\newcommand {\fel}{{\mathfrak{el}}}

\newcommand {\fwk}{{\mathfrak{wk}}}

\newcommand {\del}{\partial}

\def\smat#1{\left(\begin{smallmatrix}#1\end{smallmatrix}\right)}
\def\stab#1{\tiny\begin{tabular}
  {>{\centering\arraybackslash}m{8mm}
  >{\centering\arraybackslash}m{50mm}
  >{\centering\arraybackslash}m{60mm}
  >{\centering\arraybackslash}m{7mm}}
  &&&\\%[-3mm]
#1\end{tabular}}
\def\attn{\fbox{!!!}}
\newcommand{\adj}{\operatorname{adj}}

\usepackage{DLdef1}

\makeatletter
\theoremstyle{thm}
\@ifundefined{Theorem}{}{}
\@ifundefined{Lemma}{}{}
\@ifundefined{Corollary}{}{}
\@ifundefined{Proposition}{}{}
\@ifundefined{Statement}{}{}
\@ifundefined{Conjecture}{}{}
\@ifundefined{Question}{}{}
\theoremstyle{definition}
\@ifundefined{Remark}{}{}
\makeatother
\providecommand{\sign}{\mathop{\rm sign}\nolimits}

\begin{document}

%\begin{frontmatter}

\title[Inverses of Cartan matrices]{Inverses of Cartan matrices\\ of Lie algebras and  Lie superalgebras}

\author[D.~Leites, O.~Lozhechnyk]{Dimitry Leites${}^{a,*}$, Oleksandr Lozhechnyk${}^b$}
\address{
 ${}^*$The corresponding author; ${}^{a}$ Department of Mathematics,
Stockholm University, Albanov\"agen 28, SE-114 19 
Stockholm, Sweden; dimleites@gmail.com\\
${}^b$Ukraine, 03127, Kyiv, Hlushkova Avenue, 4e, Faculty of Mechanics and Mathematics of Taras Shevchenko National University of Kyiv; alozhechnik@gmail.com}

%\begin{keyword}

 % \MSC[2010]
 \makeatletter
\@namedef{subjclassname@2020}{\textup{2020} Mathematics Subject Classification}
\makeatother
%\end{keyword}

\begin{abstract} The inverses of indecomposable Cartan matrices are computed for finite-dimension\-al Lie algebras and Lie superalgebras over fields of any characteristic, and for hyperbolic (almost affine) complex Lie (super)algebras. This enables one to express the fundamental weights in terms of simple roots. We discovered three new phenomena concerning almost affine complex Lie (super)algebras: (a) several of the inverses of Cartan matrices have all their elements negative (not just non-positive, as they should be according to an a priori characterization due to Zhang Hechun); (b) the 0s only occur on the main diagonals of the inverses; (c) the determinants of pairwise non-equivalent Cartan matrices of the  simple  Lie (super)algebra may differ (in any characteristic). Here, we explain all three. The inverse has a 0 at the $i$th place of its main diagonal precisely if the deletion of the $i$th node of the Dynkin diagram leaves an affine algebra, and it has no 0 anywhere else; so (a) singles out the almost affine Lie (super)algebras for which the deletion of any node leaves a finite-dimensional one.
We interpret most of the results of Wei Yangjiang and Zou Yi Ming in
%, Inverses of Cartan mat\-rices of Lie algebras and Lie superalgebras, 
Linear Alg. Appl., 521 (2017) 283--298 as inverses of the Gram matrices of non-degenerate invariant symmetric bilinear forms on the (super)algebras considered, not of Cartan matrices, and give more adequate references. %In particular, the inverses of Cartan matrices of simple Lie algebras were already published, starting with Dynkin's paper in 1952, see also Table 2 in Springer's book by Onishchik and Vinberg (1990).

\end{abstract}

%\fntext[]

\thanks{We thank A. Lebedev for great help, S.~Konstein for proofs of formulas \eqref{KW1} and \eqref{osp1}--\eqref{osp6}, R.~Stekolshchik and \`E.B.Vinberg for helpful comments. }

%\end{frontmatter}

\maketitle

\markboth{\itshape Dimitry Leites\textup{,} Oleksandr Lozhechnyk}
{{\itshape Inverses of Cartan matrices}}

\section{Introduction}\label{intro}

\ssec{General remarks} The problem \lq\lq find the inverses of Cartan matrices" might look as a topic of a somewhat boring course work for a first year student taking linear algebra.

However, the explicit answer is needed in a number of situations, so the problem formulated above was solved --- for finite-dimensional simple Lie algebras over $\Cee$ --- long ago, see Dynkin's paper \cite[Ch.1]{D1}. Here are examples where these inverses are used.

$\bullet$ The inverse of the Cartan matrix $A$ (or rather of its transposed, $A^T$)  of a given simple finite-dimensional Lie algebra over $\Cee$ is a tool to obtain the set of \textbf{fundamental weights} from the set of simple roots, see, e.g., \cite[Ch.6, \S 1, Subsect. 10, eq.~(14)]{Bbk}, where \textbf{the result} of using the inverse of $A^T$ is given, \textbf{but not the inverse of $A^T$ itself} (being considered, perhaps, well-known  by that time not only to experts for more than a decade). The inverses of matrices $A^T$ are contained, among other very useful tables, in the book  by Dynkin's former Ph.D. students, see \cite[Table 2]{OV}. Any sufficiently comprehensive textbook on representations of Lie algebras (e.g., \cite{Hm}, where the inexperienced reader can find definitions of the term ``hyperbolic Lie algebra\rq\rq\  and related notions) reproduces the  statement of Bourbaki because the fundamental weights are important:  (1) they form a basis in the weight lattice, (2) any finite-dimensional irreducible representation of a simple finite-dimensional Lie algebra over $\Cee$ is a direct summand in the tensor product of tensor powers of fundamental representations.

$\bullet$  We found out that we need the results of this paper over fields of positive characteristic in our study of the notions related to the Duflo--Serganova functor, a newly discovered powerful tool in representation theory of Lie superalgebras, see \cite{KLLS, GKL}.

$\bullet$  Another natural motivation to invert the Cartan matrices that springs to mind is the fact that the Cartan matrix of $\fsl(n+1)$ is the matrix of the difference operator corresponding to $-(\frac{d}{dx})^2$. Under the name \textit{Toeplitz matrix} denoted by $K_n$ it is known to experts in the method of finite differences. Together with the Cartan matrix\footnote{It is denoted by $T_n$ in \cite[Section~1.1]{Str}. Perhaps, in honor of Toeplitz.} of the Lie superalgebra $\fosp(1|2n)$, which is the Cartan matrix of $\fo(2n+1)$ with the last row divided by 2, these matrices are the invertible two of the \lq\lq four special matrices" (actually, series of matrices) spoken about in \cite[Section~1.1]{Str}. (The two non-invertible of the \lq\lq four special matrices" correspond to an affine Lie algebra and an affine Lie superalgebra, respectively; for the classification of affine Lie (super)algebras, see~\cite{CCLL}.)

$\bullet$  For inverses of \textbf{partial} Cartan matrices, and some examples of their usage, see \cite[Tables A.11--A.13]{Stk}; the inverses of Cartan matrices are also recalled there (Tables A.14--A.15). 

$\bullet$  For several  more instances where the  inverses of Cartan matrices of Lie algebras are considered, see the bibliography in \cite{WZ}. In particular, there is a reference to a paper by Lusztig and Tits which, though interesting, is rather unclear at places, so we give an elucidation.

\ssec{A result due to Lusztig and Tits}\label{LTsec} 
The title of the paper \cite{LT} is misleading  since the paper contains a much stronger result than an already   known one.  Lusztig and Tits  established certain common properties of the inverses of matrices \textit{much more general} than the Cartan matrices of simple finite-dimensional Lie algebras over $\Cee$, see Proposition \ref{LTT} (not explicitly formulated in \cite{LT}, but  follows from a more general statement formulated in \cite{LT}).

Recall that a \textit{forest} is a disjoint union of \textit{trees}, the latter being undirected graphs in which any two vertices are connected by exactly one path. A \textit{graph loop} is an edge going in and out of the same vertex.  A \textit{simple graph} is an unweighted, undirected graph containing no graph loops or multiple edges. For a finite simple graph, the  \textit{adjacency matrix} $A$  is an $n\times n$ matrix, where $n$ is the quantity of vertices,  with 0s on its diagonal and $A_{ij}=A_{ji}=1$ if the $i$th and $j$th vertices are connected. Recall, see \cite{V}, that $B$ is a \textit{main submatrix} of a given matrix $C$ if $B$ is obtained by deleting from $C$ any number of pairs (a row, a column) each pair intersecting on the main diagonal.

Given the \textit{adjacency} (not \textit{incidence}, it is a misprint in \cite{LT}) matrix $A$ of a forest, the matrices $C$ of the same size as $A$ are constructed in \cite{LT}  by replacing any of the diagonal elements and non-zero off-diagonal elements of $A$ with elements in $\Ree$ --- conditions a) and b) in the list \eqref{LTcon} below --- provided two more conditions --- c) and d) --- are met:
\be\label{LTcon}
\begin{minipage}[l]{14cm}
\begin{description}
\item[a)] $C_{ij} = 0 \Longleftrightarrow C_{ji}=0$;
\item[b)] if $i\neq j$, then $A_{ij}=1\Longleftrightarrow C_{ij} \neq 0$;
\item[c)] $C_{ij}\leq 0$ if $i\neq j$;
\item[d)] $\det C>0$ and $\det B>0$ for all main submatrices $B$ of $C$; in parti\-cular, this means that $C_{ii}>0$ for all $i$.\end{description}
\end{minipage}
\ee
More exactly, conditions a) and b) are stronger than the conditions imposed on $C$ in \cite{LT}. Conditions c) and d) (or something equivalent to them)  are implicit in \cite{LT}, but used in the proof of  the following statement (in the form formulated by A.~Lebedev).

\sssbegin{Proposition}[A Corollary of \cite{LT}]\label{LTT}  For any tree, let $A$ be the adjacency matrix; let $C$ be a matrix satisfying conditions \eqref{LTcon}. 
Then, all elements of  $C^{-1}$ are positive.
\end{Proposition}

The proof in \cite{LT} runs through formula $(3.1)$ of that paper, which expresses an element
of the adjugate matrix in terms of the \textbf{unique} path joining the $i$th and the $j$th
nodes of the tree. That formula holds for an arbitrary matrix as soon as the unique path is
replaced with the sum over \textbf{all} paths; this generalization, suggested by A.~Lebedev, is
what makes the tree superfluous.

\sssbegin{Lemma}[The adjugate as a sum over paths]\label{paths} Let $C$ be an $n\times n$
matrix over a commutative ring, let $V=\{1,\dots,n\}$, and let $\adj C$ be the adjugate matrix,
so that $C\cdot\adj C=\det C\cdot\One$. For any $i,j\in V$, let $\cP_{ij}$ be the set of
self-non-intersecting paths $p\colon i=v_0,v_1,\dots,v_\ell=j$ such that $C_{v_{t-1}v_t}\neq0$
for $1\leq t\leq \ell$, and let $C[S]$ be the main submatrix of $C$ on a subset $S\subset V$,
the determinant of the empty matrix being $1$. Then,
\be\label{pathexp}
(\adj C)_{ij}=\sum_{p\in\cP_{ij}}\left(\prod_{(a,b)\in p}(-C_{ab})\right)\cdot
\det C[V\setminus p],
\ee
where the product is taken over the arcs $(a,b)=(v_{t-1},v_t)$ of $p$. For $i=j$, the only path
is the trivial one, and \eqref{pathexp} reads $(\adj C)_{ii}=\det C[V\setminus\{i\}]$.
\end{Lemma}

\begin{proof} By Jacobi's formula, $(\adj C)_{ij}=\frac{\del(\det C)}{\del C_{ji}}$. In the sum
\[
\det C=\sum_{\sigma\in S_n}\sign(\sigma)\prod_{k\in V}C_{k\sigma(k)},
\]
the summands containing the factor $C_{ji}$ are exactly the ones for which $\sigma(j)=i$. For
such a $\sigma$, the
cycle of $\sigma$ passing through $j$ is
\[
j\longmapsto i\longmapsto\sigma(i)\longmapsto\dots\longmapsto j;
\]
the deletion of its arc $j\longmapsto i$ leaves precisely a path $p\in\cP_{ij}$, whereas the
remaining cycles of $\sigma$ constitute an arbitrary permutation of $V\setminus p$. This cycle
has as many nodes as $p$ has, i.e., $\ell+1$, so its sign is $(-1)^{\ell}$, whereas $p$ has
$\ell$ arcs; hence
$(-1)^{\ell}\prod_{(a,b)\in p}C_{ab}=\prod_{(a,b)\in p}(-C_{ab})$. The summation over the
remaining cycles gives $\det C[V\setminus p]$, and the summation over the cycles through $j$
gives the sum over $\cP_{ij}$.\end{proof}

For a tree and $i\neq j$, the set $\cP_{ij}$ consists of one path, and \eqref{pathexp} is
formula $(3.1)$ of \cite{LT}. Nothing in the proof, however, needs a tree:

\sssbegin{Proposition}[The tree is superfluous]\label{LTTc} Let $C$ satisfy conditions
\eqref{LTcon}, and let the graph of $C$ be \textbf{connected}, not necessarily a tree. Then,
all elements of $C^{-1}$ are positive.\end{Proposition}

\begin{proof} Let $i,j\in V$. By condition a$)$, the graph of $C$ is well defined; being
connected, it contains a path from $i$ to $j$, so the sum \eqref{pathexp} is not empty. Every
arc $(a,b)$ of every $p\in\cP_{ij}$ satisfies $C_{ab}\neq0$, hence $-C_{ab}>0$ by condition
c$)$; every $C[V\setminus p]$ is a main submatrix of $C$, hence $\det C[V\setminus p]>0$ by
condition d$)$. Therefore, $(\adj C)_{ij}>0$, whereas $\det C>0$ by condition d$)$ again.
\end{proof}

\subsubsection{Remark}\label{Mmat} Conditions c$)$ and d$)$ say precisely that $C$ is a
\textit{non-singular $M$-matrix}, whereas condition a$)$ together with the connectedness of the
graph says that $C$ is \textit{irreducible}. Proposition~\ref{LTTc} is therefore a classical
theorem --- the inverse of an irreducible non-singular $M$-matrix is entry-wise positive, see
\cite[Ch.~6]{BP} --- established long before \cite{LT}. What is proper to \cite{LT} is not the
positivity, but the expression \eqref{pathexp} of the adjugate through the paths; and that
expression, as Lemma~\ref{paths} shows, has nothing to do with Cartan matrices either. This is
one more reason why the title of \cite{LT} is misleading.

\subsubsection{Remark}\label{aNeeded} Condition a$)$ is not imposed in \cite{LT}, and it may
not be dropped: for
\be\label{aCEx}
C=\smat{2&-1&0\\0&2&-1\\0&0&2},\qquad\text{we have}\qquad
C^{-1}=\frac{1}{8}\smat{4&2&1\\0&4&2\\0&0&4},
\ee
conditions b$)$--d$)$ hold, the graph is a tree, and yet $C^{-1}$ has three zero elements. The
reason is visible in \eqref{pathexp}: the arcs of a path are oriented, and here every arc runs
one way only, so $\cP_{ij}=\emptyset$ for $i>j$. What survives without condition a$)$ is the
non-negativity of the elements of $C^{-1}$.

\ssec{Two interpretations of $\det C$} For the Cartan matrix $C$ of the simple  finite-dimensional Lie algebra over $\Cee$, we have, see \cite{OV}:
\be\label{LTdet}
\begin{minipage}[l]{14cm}
\begin{description}
\item[a)] $\det C$ is equal to the order of the center of the simply connected compact group corresponding to $C$.
\item[b)] $\det C=|P/Q|$, where $P$ and $Q$ denote the weight lattice and root lattice, respectively.
\end{description}
\end{minipage}
\ee
Both interpretations \eqref{LTdet} are particular cases of one statement about integer
matrices, which extends them to Lie superalgebras and to $p>0$, see Section~\ref{detint}.

\ssec{Cartan matrices of hyperbolic Lie algebras and almost affine Lie superalgebras} The Cartan matrix of any hyperbolic Lie algebra is the one that,  its $i$th row and its $i$th column being deleted for any $i$, becomes the Cartan matrix of the direct sum of finite-dimensional or affine Kac--Moody algebras. For the definition of Lie (super)algebras $\fg(A)$ with Cartan matrix $A$, see Section~2; for superization of hyperbolic Lie algebras, called \textit{almost affine} Lie superalgebras, and classification of both types,  see \cite{CCLL}. 

For Cartan matrices $C$ of hyperbolic Lie algebras Zhang Hechun \cite{ZH}  established that 
\be\label{ZH}
\text{\lq\lq the entries of $C^{-1}$ are \textbf{non-positive} rational numbers".}
\ee

Recall that a given Lie  \textbf{super}algebra $\fg(A)$ with Cartan matrix $A$ is \textit{almost affine} if it is not finite-dimensional or affine Kac--Moody, but any its main submatrix corresponds to a direct sum of finite-dimensional or affine Kac--Moody superalgebras, \textbf{and all  Cartan matrices of $\fg(A)$ are almost affine}. The almost affine Lie superalgebras, except for $NS3_{82}-NS3_{86}$, also have the property \eqref{ZH}; this follows from their classification and the explicit form of their Cartan matrices.

Explicit results given in Subsection~\ref{ss8.1} allow us to sharpen the description \eqref{ZH}. We discovered two new phenomena:
\begin{eqnarray}
\text{\textbf{all} elements of the matrices $C^{-1}$ marked $\fbox{!!!}$ in Subsection~\ref{ss8.1} are  \textbf{negative};}\label{a}\\
\text{in each of the remaining matrices $C^{-1}$, the 0s occur only on the main diagonal.}\label{b}
\end{eqnarray}
Both phenomena are explained in Subsection~\ref{signs}.

\ssec{Natural generalization and sharpening of the problems posed in \cite{WZ}}\label{NP}{}~{}

A) \textbf{For any isotropic reflection $r$ acting on Cartan matrix $A\longmapsto r(A)$, see Subsect.~\ref{refl}, describe an algorithm for the passage $A^{-1}\longmapsto (r(A))^{-1}$.} In other words, find inverses of  not one --- selected by a random criteria (for a description of numerous possibilities in the infinite-dimensional case, see \cite{Eg}) --- Cartan matrix of a given Lie superalgebra $\fg(A)$, but of all its Cartan matrices (for $\dim \fg(A)<\infty$).

B) For the exceptional simple Lie (super)algebras, give the complete explicit answer in the following two cases:

Ba) for finite-dimensional ones (in any characteristic);

Bb) for hyperbolic Lie algebras and almost affine Lie superalgebras.

C) Investigate what can be said about inverses of Cartan matrices, if exist, for \textbf{various} other types of Lie (super)algebras corresponding to the cases of Vinberg's theorem, see \cite{V, ZH}. \textbf{This problem remains open}.

D) \textbf{Explain the phenomena \eqref{a} and \eqref{b}}: why are \textbf{all} the elements of the inverses marked $\fbox{!!!}$ negative, and why do the 0s in each of the other inverses occur only on the main diagonal? \textbf{Solved}, see Theorem~\ref{ThD}: the $0$s sit on the main diagonal at the $i$th place for precisely those $i$ for which the deletion of the $i$th node of the Dynkin diagram leaves an \textbf{affine} Lie (super)algebra; hence, there is no $0$ at all if and only if the deletion of any node leaves a \textbf{finite-dimensional} one.

E) \textbf{Explain the sign of $\det A$ where $A$ is not symmetrizable.} For a symmetrizable Cartan matrix, the sign of the determinant is the signature of the root lattice, see \eqref{signdet}, and hence the determinant of an almost affine Lie (super)algebra is negative. For the 131 non-symmetrizable ones among the 370 matrices of Subsections~\ref{ss8.1} and \ref{ssHyp}, there is no invariant form, hence no signature --- and yet the determinant is negative for every one of them, see Statement~\ref{Lor}. \textbf{This problem remains open}; for the five Lie superalgebras whose determinants are not even rational, what survives of the sign is described in Subsection~\ref{Galois}.

\subsubsection{Disclaimer}\label{Discl} The \lq\lq fundamental weights" 
are not as important in the representation theories of modular Lie algebras and of Lie superalgebras (in any characteristic) due to the lack of complete reducibility of their representations, and due to existence of deforms of vacuum (highest or lowest) weight modules, and of modules without vacuum weight vector.

\textbf{However}, even for $p>0$, these fundamental weights are as good as their namesakes over $\Cee$ if we confine ourselves to the \textbf{restricted} representations, the very first (if not the only) ones to be studied from a point of view of the geometer (P.~Deligne, see his Appendix to \cite{LL}).

\ssec{Our results} In Section~\ref{S3} and formula \eqref{Cinv}  we described an algorithm for the passage $A^{-1}\longmapsto (r(A))^{-1}$, where $r$ is any isotropic reflection (solution of Problem A, see Subsection~\ref{NP}).

In Sections~\ref{S4}--\ref{AA}  
we computed the inverses of indecomposable Cartan matrices  for finite-dimensional Lie algebras and Lie superalgebras over fields of positive characteristic and for almost affine (hyperbolic) infinite-dimensional complex Lie (super)algebras (solution of  Problem B, see Subsection~\ref{NP}). (For the classifications of Lie (super)algebras with indecomposable Cartan matrices over $\Cee$ in the almost affine (hyperbolic) case, and finite-dimensional over algebraically closed fields of characteristic $p>0$, see \cite{CCLL} (with prerequisites in  \cite{Se, HS}) and \cite{BGL}, respectively.)

We have discovered phenomena \eqref{a} and \eqref{b}, and explained them in
Subsection~\ref{signs} (solution of Problem D), see Subsection~\ref{NP}): both are governed by
the type of the Lie (super)algebra left after the deletion of one node of the Dynkin diagram,
see Theorem~\ref{ThD}. On the way, we generalize the main formula of \cite{LT} to an arbitrary
matrix, see Lemma~\ref{paths}, and remove the hypothesis that the graph be a tree from the
result \cite{LT} is quoted for, see Proposition~\ref{LTTc}. 
We have discovered that the determinants of pairwise non-equivalent Cartan matrices of the  simple  Lie (super)algebra may differ (in any characteristic), and, in Section~\ref{detint}, we explain why: the reflection itself does not change $\det A$, see Lemma~\ref{reflDet}, so the whole of the difference comes from the renormalization, see Corollary~\ref{corDet}. In the same Section, we show that $|\det A|$ is the index \eqref{PQ} of the root lattice in the lattice dual to the chosen elements $h_i$ of the maximal torus, whereas the sign of $\det A$ is the signature of the root lattice whenever $A$ is symmetrizable, see \eqref{signdet}; in particular, a determinant equal to 1 says that the two lattices coincide.

For the serial Lie superalgebras, the inverses of Cartan matrices are explicitly given for certain \lq\lq basic" Cartan matrices from which all the other Cartan  matrices of the given Lie superalgebra are obtained by means of isotropic (odd) reflections as explained in Subsection~\ref{refInv}.

We also corrected and widened the list of references (and the range of applications) of inverses of Cartan matrices as compared with those given in \cite{WZ}, e.g., the explicit form of inverse Cartan matrices of finite-dimensional simple  Lie algebras is reproduced in  \cite{WZ} as if new, though well-known, see \cite{D1,OV, Stk}. These omissions and the desire to solve the problem considered, but not solved, in  \cite{WZ}, except for occasional coincidences of Cartan matrices with Gram matrices of non-degenerate invariant symmetric bilinear forms (NISes) on the (super)algebras considered, is what prompted our work.

\subsubsection{Remark} V.~Kac was the first (see \cite{K0, K}) to realize that certain finite-dimensional Lie superalgebras have analogs of Cartan matrices and defined them imitating the definition for Lie algebras, compare \cite{K0, K} with \cite{Kapp}, where the exceptional simple Lie superalgebras first appeared. In these cases, V.~Kac listed the indecomposable Cartan matrices (with a gap, corrected in \cite{Se, vdL}, where infinite-dimensional generalizations with symmetrizable Cartan matrices and of polynomial growth were also considered). For further improvements of the definitions, see  \cite{HS, CCLL, BGL}.

In \cite{WZ}, the Gram matrices of the non-degenerate symmetric bilinear form on the space of roots, that never explicitly appeared before, but can be recovered from the data in \cite{Se1}, or by symmetrizing Cartan matrices, were misattributed to works of V.~Kac and called \textit{Cartan matrices}. Kac never used such matrices (and would hardly apply the term \textit{Cartan matrix} to a matrix with both 2 and $-2$  appearing simultaneously on the main diagonal). Analogs of Cartan matrices $A$ with $A_{ii}\leq 0$ for all $i$ were introduced by Borcherds, see \cite{B, R} and references in \cite{CCLL}.

\section{Chevalley generators, Cartan matrices, reflections (from \cite{CCLL, BGL})}\label{S2}

\ssec{Chevalley generators and Cartan matrices} Let us start with the construction of Lie
(super)algebras with Cartan matrix. Let $A=(A_{ij})$ be an $n\times n$-matrix
whose entries lie in the ground field $\Kee$. Let $\rk
A=n-l$. It means that there exists an $l\times n$-matrix
$T=(T_{ij})$ such that
\begin{equation}\label{rankCM}
\begin{tabular}{l}
a) the rows of $T$ are linearly independent;\\
b) $TA=0$ (or, more precisely, ``zero $l\times n$-matrix\rq\rq).
\end{tabular}
\end{equation}
Indeed, if $\rk A^T=\rk A=n-l$, then there exist $l$ linearly
independent vectors $v_i$ such that $A^Tv_i=0$; set
\[
T_{ij}=(v_i)_j.
\]

Let the elements $e_i^\pm$ and $h_i$, where $i=1,\dots,n$, generate
a Lie superalgebra denoted $\tilde \fg(A, I)$, where $I=(p_1, \dots
p_n)\in(\Zee/2)^n$ is a collection of parities ($p(e_i^\pm)=p_i$),
free except for the relations
\begin{equation}\label{gArel_0}
{}[e_{i}^+, e_{j}^-] = \delta_{ij}h_i; \quad [h_j, e_{j}^\pm]=\pm
A_{ij}e_{j}^\pm\text{~~and~~}[h_i, h_j]=0 \text{~~for any $i, j$}.
\end{equation}
Let Lie (super)algebras with Cartan matrix $ \fg(A,
I)$ be the quotient of $\tilde \fg(A, I)$ modulo
the ideal explicitly described in \cite{GL, BGL3, BGLL}.

By abuse of notation we denote by $e_j^\pm$ and $h_j$ ---
the elements of $\tilde{\fg}(A,I)$
--- also their images in $\fg(A,I)$ and $\fg^{(i)}(A,I)$ and call these images, and their pre-images, the \textit{Chevalley generators} of $\tilde{\fg}(A,I)$, $\fg(A,I)$, and $\fg^{(i)}(A,I)$, cf. Subsection~\ref{SsChev}.

\subsubsection{In small font} The additional to \eqref{gArel_0} relations that turn $\tilde \fg(A,
I)$ into $\fg(A, I)$ are of the form $R_i=0$ whose left sides are
implicitly described, for the general Cartan matrix with entries in
$\Kee$, as
\begin{equation}\label{myst}
\text{\begin{minipage}[c]{14cm} the $R_i$ that generate the ideal
$\fr$ maximal among the ideals of $\tilde\fg (A, I)$ whose
intersection with the span of the above $h_i$ and the $d_j$
described in eq. \eqref{central3} is~$0$.\end{minipage}}
\end{equation}

Set
\begin{equation}\label{central}
c_i=\mathop{\sum}\limits_{1\leq j\leq n} T_{ij}h_j, \text{~~where~~}
i=1,\dots,l.
\end{equation}
Then, from the properties of the matrix $T$, we deduce that
\begin{equation}\label{central1}
\begin{tabular}{l}
a) the elements $c_i$ are linearly independent;\\
b) the elements $c_i$ are central, because\\
$[c_i,e_j^\pm]=\pm\left(\sum\limits_{1\leq k\leq n} T_{ik}A_{kj}\right)
e_j^\pm=\pm (TA)_{ij} e_j^\pm $.
\end{tabular}
\end{equation}
The existence of central elements means that the linear span of all
the roots is of dimension $n-l$ only. (This can be explained even
without central elements: The weights can be considered as
column-vectors whose $i$-th coordinates are the corresponding
eigenvalues of $\ad_{h_i}$. The weight of $e_i$ is, therefore, the
$i$-th column of $A$. Since $\rk A=n-l$, the linear span of all
columns of $A$ is $(n-l)$-dimensional just by definition of the
rank. Since any root is an (integer) linear combination of the
weights of the $e_i$, the linear span of all roots is
$(n-l)$-dimensional.)

This means that some elements which we would like to see having
different (even opposite if $p=2$) weights, actually,  have identical
weights. To remedy this, we do the following: let $B$ be an arbitrary
$l\times n$-matrix such that
\begin{equation}\label{matrixB}
\text{the~~}(n+l)\times
n\text{-matrix~~}\begin{lmatrix}A\\B\end{lmatrix}\text{~~has
rank~}n.
\end{equation}
Let us add to the algebra $\fg=\tilde \fg(A, I)$ (and hence $\fg(A,
I)$) the grading elements $d_i$, where $i=1,\dots,l$, subject to the
following relations:
\begin{equation}\label{central3}
{}[d_i,e_j^\pm]=\pm B_{ij}e_j;\quad [d_i,d_j]=0;\quad [d_i,h_j]=0
\end{equation}
(the last two relations mean that the $d_i$ lie in the Cartan
subalgebra, and even in the maximal torus which will be denoted by
$\fh$).

Note that these $d_i$ are \textit{outer} derivations of $\fg(A,
I)^{(1)}$, i.e., they can not be obtained as linear combinations of
brackets of the elements of $\fg(A, I)$ (i.e., the $d_i$ do not lie
in $\fg(A, I)^{(1)}$).

\ssec{Roots and weights}\label{roots} In this subsection, $\fg$
denotes  $\fg(A,I)$ or $\tilde{\fg}(A,I)$.

Let $\fh$ be the span of the $h_i$ and the $d_j$. The elements of
$\fh^*$ are called \textit{weights}. For a given
weight $\alpha$, the \textit{weight subspace} of a given
$\fg$-module $V$ is defined as
\[
V_\alpha=\{x\in V\mid \text{an integer $N>0$ exists such that
$(\alpha(h)-\ad_h)^N x=0$ for any $h\in\fh$}\}.
\]

Any non-zero element $x\in V_\alpha$ is said to be \textit{of weight
$\alpha$}. For the roots, which are particular cases of weights if
$p=0$, the above definition is inconvenient because it does not lead
to the modular analog of the following useful statement.

\sssbegin{Proposition}[The root decomposition of $\fg$, \cite{K}]\label{rootdec} Over $\Cee$, the space
of any Lie algebra $\fg$ can be represented as a direct sum of
subspaces
\begin{equation}\label{rootdeceq}
\fg=\mathop{\bigoplus}\limits_{\alpha\in \fh^*} \fg_\alpha.
\end{equation}
\end{Proposition}
Note that if $p=2$, it might happen that $\fh\subsetneq\fg_0$. (For
example, all weights of the form $2\alpha$ over $\Cee$ become 0 over
$\Kee$.)

To salvage the formulation of Statement in the modular case with
minimal changes, \textit{at least for the Lie (super)algebras $\fg$ with
Cartan matrix} --- and only this case we will have in mind speaking
of roots, --- we decree that the elements $e_i^\pm$ with the same
superscript (either $+$ or $-$) correspond to linearly independent
\textit{roots} $\alpha_i$, and any root $\alpha$ such that
$\fg_\alpha\neq 0$ lies in the $\Zee$-span of
$\{\alpha_1,\dots,\alpha_n\}$, i.e.,
\begin{equation}\label{rootdeceqnew}
\fg=\mathop{\bigoplus}\limits_{\alpha\in
(\Zee\alpha_1\oplus \dots\oplus \Zee\alpha_n)} \fg_\alpha.
\end{equation}

Thus, $\fg$ has a $\Ree^n$-grading such that $e_i^\pm$ has grade
$(0,\dots,0,\pm 1,0,\dots,0)$, where $\pm 1$ stands in the $i$-th
slot (this can also be considered as $\Zee^n$-grading, but we use
$\Ree^n$ for simplicity of formulations). If $p=0$, this grading is
equivalent to the weight grading of $\fg$. If $p>0$, this grading and the weight grading
may be pairwise non-equivalent; in particular, if $p=2$, then the elements
$e_i^+$ and $e_i^-$ have the same weight. \textbf{This is why in what
follows we consider roots as elements of $\Ree^n$, not as weights}. 

Any non-zero element $\alpha\in\Ree^n$ is called \textit{a
root} if the corresponding eigenspace of grade $\alpha$
(which we denote $\fg_\alpha$ by abuse of notation) is non-zero. The
set $R$ of all roots is called \textit{the root system} of $\fg$.

Clearly, the subspaces $\fg_\alpha$ are purely even or purely odd,
and the corresponding roots are said to be \textit{even} or
\textit{odd}.

\ssec{Systems of simple and positive roots} In this subsection,
$\fg=\fg(A,I)$, and $R$ is the root system of $\fg$.

For any subset $B=\{\sigma_{1}, \dots, \sigma_{m}\} \subset R$, we
set (we denote by $\Zee_{+}$ the set of non-negative integers):
\[
R_{B}^{\pm} =\{ \alpha \in R \mid \alpha = \pm \sum n_{i}
\sigma_{i},\;\text{where }\;n_{i} \in \Zee_{+} \}.
\]
The set $B$ is called a \textit{system of simple roots~} of $R$ (or $\fg$) if $ \sigma_{1}, \dots ,
\sigma_{m}$ are linearly independent and $R=R_B^+\cup R_B^-$. Note
that $R$ contains basis coordinate vectors, and therefore spans
$\Ree^n$; thus, any system of simple roots contains exactly $n$
elements.

A subset $R^+\subset R$ is called a \textit{system of positive
roots~} of $R$ (or $\fg$) if there
exists $x\in\Ree^n$ such that
\begin{equation}\label{x}
\begin{split}
 &(\alpha,x)\in\Ree\backslash \{0\}\text{ for all $\alpha\in R$},\\
 &R^+=\{\alpha\in R\mid (\alpha,x)>0\}.
\end{split}
\end{equation} (Here $(-,-)$ is the standard Euclidean inner
product in $\Ree^n$.) Since $R$ is a finite (or, at least, countable
if $\dim \fg(A)=\infty$) set, so the set
\[\{y\in\Ree^n\mid\text{there exists $\alpha\in R$ such that }
(\alpha,y)=0\} \] is a finite/countable union of $(n-1)$-dimensional
subspaces in $\Ree^n$, so it has zero measure. So for almost every
$x$, condition (\ref{x}) holds.

By construction, any system $B$ of simple roots is contained in
exactly one system of positive roots, which is precisely $R_B^+$.

\sssbegin{Statement}[A system of positive roots contains exactly one system of simple
roots]\label{possimple} Any finite system $R^+$ of positive roots of
$\fg$ contains exactly one system of simple roots. This system
consists of all the positive roots (i.e., elements of $R^+$) that
can not be represented as a sum of two positive
roots.\end{Statement}

We can not give an \textit{a priori} proof of the fact that each set
of all positive roots each of which is not a sum of two other
positive roots consists of linearly independent elements. This is,
however, true for finite-dimensional Lie algebras and Lie
superalgebras of the form $\fg(A)$ if $p\neq 2$, see \cite{BGL, CCLL}.

\ssec{Normalization convention}\label{normA} Clearly,
\begin{equation}
\label{rescale} \text{the rescaling
$e_i^\pm\mapsto\sqrt{\lambda_i}e_i^\pm$, sends $A$ to $A':=
\diag(\lambda_1, \dots , \lambda_n)\cdot A$.} 
\end{equation}

Two pairs $(A, I)$ and $(A', I')$ are said to be \textit{equivalent}
if $(A', I')$ is obtained from $(A, I)$ by a composition of a
permutation of parities and a rescaling $A' = \diag (\lambda_{1},
\dots, \lambda_{n})\cdot A$, where ${\lambda_{1}\cdot\dots \cdot\lambda_{n}\neq
0}$. Clearly, equivalent pairs determine isomorphic Lie
superalgebras.

The rescaling affects only the matrix $A_B$, not the set of parities
$I_B$. The Cartan matrix $A$ is said to be
\textit{normalized} if 
\begin{equation}
\label{norm} A_{jj}=0\text{~~ or 1, or 2.}
\end{equation}
We let $A_{jj}=2$ only if $i_j=\ev$; in order to eliminate possible confusion, we write $A_{jj}=\ev$ or
$\od$ if $i_j=\ev$, whereas if $i_j=\od$, we write $A_{jj}=0$ or
$1$. 

Normalization conditions correspond to the \lq\lq natural" Chevalley generators of the most usual \lq\lq building blocks" of finite-dimensional Lie (super)algebras with Cartan matrix: $\fsl(2)$ if $A_{jj}=2$, $\fsl(1|1)$ if $A_{jj}=0$, and $\fosp(1|2)$ if $A_{jj}=1$. (In this paper we do not need $A_{jj}=\ev$ or
$\od$, see \cite{BGL,CCLL}.)

\textbf{We will only consider
normalized Cartan matrices; for them, we do not have to indicate the set of parities
$I$.}

\subsubsection{Warning}\label{War} Unlike the case of simple finite-dimensional Lie algebras over $\Cee$, where the normalized Cartan matrix $A$ is uniquely defined, generally this is not so: each row with a 0 or $\ev$ on the main diagonal can be multiplied by
any nonzero factor; usually (not only in this paper) we multiply the
rows so as to make $A_{B}$ symmetric, if possible. Which version of the Cartan matrix should be considered as its \lq\lq normal form"? The defining relations give the answer:

The normalized Cartan matrix is used, for example, to describe presentation of the given Lie superalgebra (relations between, or rather \textit{among} (because analogs of the Serre relations in the super setting involve \textit{several} generators, see \cite{GL, BGL, BGLL}) the Chevalley generators).

\ssec{Equivalent and connected to each other systems of simple roots} \label{EqSSR}
For a system of simple roots $B=\{\alpha_1,\dots,\alpha_n\}$, choose
non-zero elements $e_i^\pm$ in the 1-dimensional (by definition)
superspaces $\fg_{\pm\alpha_i}$; set $h_{i}=[e_{i}^{+}, e_{i}^-]$,
let $A_{B} =(A_{ij})$, where the entries $A_{ij}$ are recovered from
relations \eqref{gArel_0}, and let $I_{B}=\{p(e_{1}), \cdots,
p(e_{n})\}$. Lemma \ref{serg} assures us that all the pairs $(A_B,I_B)$
are \textit{connected} to each other.

Two systems of simple roots $B_{1}$ and $B_{2}$ are said to be
\textit{equivalent} if the pairs $(A_{B_{1}}, I_{B_{1}})$ and
$(A_{B_{2}}, I_{B_{2}})$ are equivalent.

It would be nice to find a convenient way to fix some distinguished
pair $(A_B,I_B)$ in the equivalence class. For the role of the
``best\rq\rq\  (first among equals) order of indices we propose the one
that minimizes the value
\begin{equation}\label{minCM}
\max\limits_{i,j\in\{1,\dots,n\}\text{~such that~}(A_B)_{ij}\neq
0}|i-j|
\end{equation}
(i.e., gather the non-zero entries of $A$ as close to the main
diagonal as possible). Observe that this numbering differs from the
one that Bourbaki use for the $\fe$ type Lie algebras.

\subsubsection{Chevalley generators and Chevalley bases}\label{SsChev} We
often denote the set of generators of $\fg(A,I)$ and $\fg^{(i)}(A,I)$ corresponding to a normalized
Cartan matrix by $X_{1}^{\pm},\dots , X_{n}^{\pm}$ instead of
$e_{1}^{\pm},\dots , e_{n}^{\pm}$; and call these generators, together with the
elements $H_i:=[X_{i}^{+}, X_{i}^{-}]$, and the derivations $d_j$,
see~\eqref{central3}, the \textit{Chevalley
generators}. 

For $p=0$ and normalized Cartan matrices of simple finite
dimensional Lie algebras, there exists only one (up to signs) basis
containing $X_i^\pm$ and $H_i$ in which $A_{ii}=2$ for all $i$ and
all structure constants are integer, cf. \cite{St}. Such a basis is
called the \textit{Chevalley}  basis.

Observe that, having normalized the Cartan matrix of $\fsp(2n)$ so
that $A_{ii}=2$ for all $i\neq n$, but $A_{nn}=1$, we get
\textbf{another} basis with integer structure constants. We think
that this basis also qualifies to be called \textit{Chevalley
basis}; for Lie superalgebras, and if $p=2$, such normalization is a
must.

\sssbegin{Conjecture}[The Chevalley basis is unique also for $p>2$]\label{Chev} If $p>2$, then for finite-dimensional Lie
(super)algebras with indecomposable Cartan matrices normalized as in
$(\ref{norm})$, there also exists only one (up to signs) analog of
the Chevalley basis. \end{Conjecture}

From \cite{BGL}: \lq\lq We had no idea how to describe analogs of Chevalley bases for $p=2$
until recently; it seems, the methods of the recent paper \cite{CR}
should solve the problem." (Now, more than a decade after, \textbf{this problem is still open}.)

\ssec{Reflections}\label{refl} Let $R^+$ be a system of positive roots of Lie
superalgebra $\fg$ over a field $\Kee$ of characteristic $p>0$, and let $B=\{\sigma_1,\dots,\sigma_n\}$ be the
corresponding system of simple roots with some corresponding pair
$(A=A_B,I=I_B)$. Then, for any $k\in \{1, \dots, n\}$, the set
$(R^+\backslash\{\sigma_k\})\coprod\{-\sigma_k\}$ is a system of
positive roots. This operation is called \textit{the reflection in
$\sigma_k$}; it changes the system of simple roots by the formulas
\begin{equation}
\label{oddrefl}
r_{\sigma_k}(\sigma_{j})= \begin{cases}{-\sigma_j}&\text{if~}k=j,\\
\sigma_j+B_{kj}\sigma_k&\text{if~}k\neq j,\end{cases}\end{equation}
where, see \cite{Leb}, whose formula $(3)$ we translate into our notation, and where
$(x)_{\Zee}$ denotes the smallest non-negative integer congruent to $x$, see
Subsection~\ref{sssIC},
\begin{equation}
\label{Boddrefl} B_{kj}=\begin{cases}
0&\text{~if~}A_{kj}=0,\\
p-1&\text{~if~}i_k=\ev,\ A_{kk}=0,\ A_{kj}\neq 0,\\
&\text{~or if~}i_k=\ev,\ p\neq2,\ A_{kk}\neq 0,\
\displaystyle\frac{A_{kj}}{A_{kk}}\not\in \Zee/p\Zee,\\
\left(-\displaystyle\frac{2A_{kj}}{A_{kk}}\right)_{\Zee}&\text{~if~}i_k=\ev,\
p\neq2,\ A_{kk}\neq 0,\ \displaystyle\frac{A_{kj}}{A_{kk}}\in \Zee/p\Zee,\\
2&\text{~if~}i_k=\ev,\ p=2,\ A_{kj}=A_{kk}\neq 0,\\
3&\text{~if~}i_k=\ev,\ p=2,\ A_{kk}\neq 0,\ A_{kj}\neq A_{kk},\ A_{kj}\neq 0,\\
1&\text{~if~}i_k=\od,\ A_{kk}=0,\ A_{kj}\neq 0,\\
2\left(-\displaystyle\frac{A_{kj}}{A_{kk}}\right)_{\Zee}&\text{~if~}i_k=\od,\
A_{kk}\neq 0,\ \displaystyle\frac{A_{kj}}{A_{kk}}\in \Zee/p\Zee,\\
2p-1&\text{~if~}i_k=\od,\ A_{kk}\neq 0,\
\displaystyle\frac{A_{kj}}{A_{kk}}\not\in \Zee/p\Zee,
\end{cases}
\end{equation}
where we consider $\Zee/p\Zee$ as a subfield of $\Kee$. Here, $B_{kj}$ is the length of
the $\alpha_k$-string through $\alpha_j$, i.e., $B_{kj}=\sup\{m\in\Zee_{\geq0}\mid
\alpha_j+m\alpha_k\in R\}$, so \eqref{oddrefl} says that the reflection sends $\alpha_j$
to the last root of that string.

The name ``reflection\rq\rq\  is used because in the case of simple
finite-dimensional complex Lie algebras this action, extended to the whole
$R$ by linearity, is a map from $R$ to $R$, and it does not depend on
$R^+$, only on $\sigma_k$. This map is usually denoted by
$r_{\sigma_k}$ or just $r_{k}$. The map $r_{\sigma_k}$ extended to
the $\Ree$-span of $R$ is reflection in the hyperplane orthogonal to
$\sigma_k$ relative the bilinear form dual to the Killing form.

The reflections in the even (odd) roots are referred to as
\textit{even} (\textit{odd}) \textit{reflections}. A simple root is called \textit{isotropic}, if the
corresponding row of the Cartan matrix has zero on the diagonal, and
\textit{non-isotropic} otherwise. The reflections that correspond to
isotropic or non-isotropic roots will be referred to accordingly.

If there are isotropic simple roots, the reflections $r_\alpha$ do
not, as a rule, generate a version of the \textit{Weyl group}
because the product of two reflections in nodes not connected by one
(perhaps, multiple) edge is not defined. These reflections just
connect a pair of ``neighboring\rq\rq\  systems of simple roots, and there is
no reason to expect that we can multiply such two distinct 
reflections. In  the case of modular Lie algebras or
of Lie superalgebras for any $p$,
the action of a given isotropic reflection (\ref{oddrefl}) can not,
generally, be extended to a linear map $R\tto R$. For Lie
superalgebras over $\Cee$, one can extend the action of reflections
by linearity to the root lattice, but this extension preserves the
root system only for $\fsl(m|n)$ and $\fosp(2m+1|2n)$, cf.
\cite{Se1}.

We would like to draw attention of the reader to an under-appreciated paper \cite{SkB}, where the analog of Weyl group was considered in relation with $\fbr(3)$.

\ssec{How reflections act on Chevalley generators}
If $\sigma_i$ is an isotropic root, then the corresponding
reflection $r_i$ sends one set of Chevalley generators into a new one:
\begin{equation}
\label{oddrefx} \tilde X_{i}^{\pm}=X_{i}^{\mp};\;\;
\tilde X_{j}^{\pm}=\begin{cases}[X_{i}^{\pm},
X_{j}^{\pm}]&\text{if $A_{ij}\neq 0, \ev$},\\
X_{j}^{\pm}&\text{otherwise}.\end{cases}
\end{equation}

The Cartan matrix $r_i(A)$ corresponding to the Chevalley generators  \eqref{oddrefx} should be obtained as described above: set
\[
\tilde H_i:=[\tilde X_{i}^+, \tilde X_{i}^-]
\]
and compute
\[
[\tilde H_i, \tilde X_{j}^+]=\tilde B_{ij}\tilde X_{j}^+.
\]
Normalize the matrix $\tilde B$ as we agreed, see Subsection~\ref{normA}; let $B$ be the normalized matrix.
Then, $r_i(A):=B=(B_{kl})$.

\subsubsection{Lebedev's lemma} Serganova, see \cite[Appendix]{LSS},
proved (for $p=0$) that there is always a chain of
reflections connecting $B_1$ with some system of simple roots $B'_2$
equivalent to $B_2$ in the sense of definition in Subsection~\ref{refl}. Here is
the $p>0$ version of this statement. 

\sssbegin{Lemma}[Any two systems of simple roots are connected by reflections, \cite{Leb}]\label{serg} For any two systems of simple roots
$B_1$ and $B_2$ of any finite-dimensional Lie superalgebra with
indecomposable Cartan matrix, there is always a chain of reflections
connecting $B_1$ with $B_2$.\end{Lemma}

\subsubsection{Important convention}\label{sssIC}
The values $-\nfrac{2A_{kj}}{A_{kk}}$ and $-\nfrac{A_{kj}}{A_{kk}}$ in \eqref{Boddrefl}
are elements of $\Kee$, while the roots are elements of a vector space
over $\Ree$. Therefore, wherever such a value occurs in \eqref{Boddrefl}, it
\textbf{should be understood} as ``\textbf{the smallest
non-negative integer congruent to} it\rq\rq; this is what the notation $(x)_{\Zee}$ used
in \eqref{Boddrefl} means. (If $\dim\fg<\infty$, these expressions are always congruent
to integers.)

This convention is important in describing Serre relations and their analogs for $p>0$, see \cite{BGL3, BGLL}. \textbf{We do not know yet what is an equally \lq\lq natural" or correct way of presenting the elements of the inverse Cartan matrices.}

\sssbegin{Remark}[Formula \eqref{Boddrefl} corrects the published one]\label{Bcorr}
Formula \eqref{Boddrefl} differs in two of its cases from the version published, without
proof, in \cite[eq.~(7.4)]{BGL} and \cite[eq.~(44)]{BGLL}, see \cite{Leb}:

for an \textbf{odd non-isotropic} node, i.e., for $i_k=\od$ and $A_{kk}\neq0$, the string
is twice as long as the published formula says --- the value is
$2\left(-\nfrac{A_{kj}}{A_{kk}}\right)_{\Zee}$, not
$\left(-\nfrac{2A_{kj}}{A_{kk}}\right)_{\Zee}$, and $2p-1$, not $p-1$, in the case where
the ratio is not in $\Zee/p\Zee$; the first instance of this is $\fosp(1|2n)$, where the
$\alpha_k$-string is of length 2;

for $p=2$ and an \textbf{even} node with $A_{kk}\neq0$, the values are 2 and 3, and the
value 3 does not occur in the published formula at all, whereas the value 2 was obtained
there by a separate convention.
\end{Remark}

\sssbegin{Statement}[The two forms of the answer agree, in the range verified]\label{Bcheck}
For $p=2,\ 3,\ 5,\ 7$, for both parities $i_k$, and for every pair $(A_{kk}, A_{kj})$ of
elements of the prime field $\Fee_p$ and of its quadratic extension $\Fee_{p^2}$ --- the
latter is what exercises the two cases of \eqref{Boddrefl} where
$\nfrac{A_{kj}}{A_{kk}}\not\in\Zee/p\Zee$ --- formula \eqref{Boddrefl} returns the first
zero of the sequence $d_{-1}=0$, $d_m=(-1)^{i_k}(d_{m-1}-A_{kj}-mA_{kk})$, which is what
$B_{kj}$ is equal to by \cite[Statement 2.1]{Leb}. That is $6246$ triples
$(i_k, A_{kk}, A_{kj})$, and the formula that \eqref{Boddrefl} replaces disagrees with the
recursion in $2762$ of them. Computed in exact arithmetic.
\end{Statement}

\subsubsection{How reflections act on Cartan matrices ([CCLL])}\label{refInv} Let $A$ be a Cartan matrix of size $n$ and $I = (p_1,\ldots , p_n)$ the vector of parities. If $p_k = \od$ and $A_{kk} = 0$, then the reflection  in the $k$th simple odd root sends $A$ to $r_k( A)$, where
\be\label{oddR}
(r_k(A))_{i j} = A_{i j} + b_i A_{k j} + c_j A_{i k},
\ee
and where
\[
b_i = \begin{cases}-2&\text{if $i = k$},\\
0&\text{if $i\neq k$ and $A_{ik} = 0$,}\\
\frac{A_{ik}}{A_{ki}}&\text{if $i\neq k$ and $A_{ik} \neq 0$;}\\
\end{cases}\quad
\text{~~
 and ~~}\quad c_j = \begin{cases}-2&\text{if $j = k$},\\
 0&\text{if $j\neq k$ and $A_{jk} = 0$,}\\
1&\text{if $j\neq k$ and $A_{jk} \neq 0$.}\\
\end{cases}
\]
This can be expressed in terms of matrices as
\[
r_k(A) =(E+B)A(E+C),
\]
where all columns of the matrix $B$, except the $k$th one, are zero, whereas the $i$th coordinate of the $k$th column-vector is $b_i$, the $i$th coordinate of the $k$th row-vector of $C$ is $c_i$, the other rows of $C$ being zero; $E$ is the unit matrix. Therefore, $B^2=C^2=0$, and 
\be\label{Cinv}
(r_k(A) )^{-1} =(E-C)A^{-1}(E-B).
\ee

The reflected matrix $r_k( A)$ might have to be normalized; the new parities are
\[
\widetilde p_j \equiv p_j + c_j \pmod 2.
\]

\section{What $\det A$ is}\label{detint}

\ssec{The determinant and the two lattices} Let $A$ be a normalized Cartan matrix of size
$n$; let $\alpha_1,\dots,\alpha_n$ be the corresponding simple roots and $h_1,\dots,h_n$
the corresponding elements of the maximal torus, so that $A_{ij}=\alpha_j(h_i)$, see
\eqref{gArel_0}. Let the $\varpi_j$ be given by $\varpi_j(h_i)=\delta_{ij}$; let
$Q:=\bigoplus_{1\leq i\leq n}\Zee\alpha_i$ and $P:=\bigoplus_{1\leq i\leq n}\Zee\varpi_i$
be the lattices they span. Since
\be\label{transit}
\alpha_j=\sum_{1\leq i\leq n}A_{ij}\varpi_i,
\ee
the matrix $A$ is the transition matrix from the $\varpi_i$ to the $\alpha_j$; the
following statement is therefore about \textbf{any} matrix with integer entries, no
matter which --- if any --- Lie (super)algebra it comes from.

\sssbegin{Lemma}[The determinant is the index of the root lattice]\label{Lin} Let $A$ be
an integer $n\times n$ matrix, $\det A\neq0$, and let $d_1\mid\dots\mid d_n$ be its
elementary divisors $($the Smith normal form$)$. Then,
$P/Q\simeq\bigoplus_{1\leq i\leq n}\Zee/d_i$; in particular,
\be\label{PQ}
|\det A|=[P:Q].
\ee
If, moreover, $A$ is symmetrizable, say $\diag(d_1,\dots,d_n)\cdot A$ is symmetric for
some $d_i>0$, then
\be\label{signdet}
\sign(\det A)=(-1)^{n_-},
\ee
where $n_-$ is the number of negative squares of that symmetrization.\end{Lemma}

\begin{proof} The first claim is the structure theorem for finitely generated abelian
groups applied to the cokernel of the map $\Zee^n\tto\Zee^n$ given by $A$, see
\eqref{transit}; the index $[P:Q]$ is the product of the $d_i$, which is $|\det A|$. For
the second claim, observe that $\det\diag(d_1,\dots,d_n)>0$, whereas the sign of the
determinant of a non-degenerate real symmetric matrix is $(-1)^{n_-}$ by Sylvester's law
of inertia.\end{proof}

Lemma~\ref{Lin} settles both interpretations \eqref{LTdet} at once, and extends them to
Lie superalgebras and to $p>0$: for the simple finite-dimensional Lie algebra over
$\Cee$, the lattice $P$ is the weight lattice and $Q$ the root lattice, so \eqref{PQ} is
item b$)$ of \eqref{LTdet}, whereas item a$)$ follows since the center of the simply
connected compact group is isomorphic to $P/Q$. Formula \eqref{signdet} says what a
\textbf{negative} determinant means:

$\det A>0$ for the finite-dimensional Lie (super)algebras, where the symmetrization is
positive definite;

$\det A=0$ for the affine ones, where it is degenerate;

$\det A<0$ for the hyperbolic Lie algebras and for the symmetrizable almost affine Lie
superalgebras, because there $n_-=1$: the symmetrization is Lorentzian. Thus, the minus
sign is the signature of the root lattice; nothing is negative except the number of
negative squares, one.

Equivalently: \textbf{the determinant of the symmetrized Cartan matrix is the
discriminant of the root lattice} $Q$, and for the hyperbolic Lie algebras that is the discriminant of a
Lorentzian lattice --- an invariant by which such lattices are classified. For the
matrices $C$ of Proposition~\ref{LTT}, condition d$)$ of \eqref{LTcon} says that $C$ is
positive definite, so $\det C>0$, and if $C$ is integral, then \eqref{PQ} applies to it
verbatim.

\sssbegin{Statement}[The Lorentzian signature, in the range verified]\label{Lor} Of the
$370$ Cartan matrices listed in Subsections~$\ref{ss8.1}$ and $\ref{ssHyp}$, exactly $239$
are symmetrizable $($all $142$ of the $238$ hyperbolic Lie algebras whose name begins
with an $H$, and $97$ of the $132$ almost affine Lie superalgebras whose name begins
with an $S$, see Subsection~$\ref{S35})$. For every one of these $239$, the number of
negative squares of the symmetrization is $n_-=1$, and \eqref{signdet} holds. For the
remaining $131$, no invariant symmetric bilinear form exists, so \eqref{signdet} is not
applicable --- and yet the determinant is negative for all $370$ of them. Computed over
$\Qee$ in exact arithmetic.\end{Statement}

Why the determinant is negative also in the $131$ non-symmetrizable cases is problem E$)$ of
Subsection~\ref{NP}; Subsection~\ref{Galois}
tells what survives of the sign for the five Lie superalgebras of Subsection~\ref{ss8.2},
whose Cartan matrices are not symmetrizable and whose determinants are not even rational.

\ssec{The signs of the elements of $A^{-1}$: solution of Problem D$)$}\label{signs}
Let $\fg(A)$ be an almost affine Lie (super)algebra, let $A$ be its Cartan matrix of size
$n\geq2$, normalized as in Subsection~\ref{normA} and printed as in
Subsections~\ref{ss8.1} and \ref{ssHyp}, so that $A_{ij}\leq0$ for $i\neq j$; let
$V=\{1,\dots,n\}$ and $A_{\hat\imath}:=A[V\setminus\{i\}]$. Recall that $A$ is indecomposable,
i.e., its graph is connected; that $\det A<0$, see Statement~\ref{Lor}; and that, by the
definition of almost affinity, every $A_{\hat\imath}$ is the Cartan matrix of the direct sum of
finite-dimensional and affine Lie (super)algebras.

\sssbegin{Theorem}[Where the 0s of $A^{-1}$ are, and why]\label{ThD} Under the above
conventions,

$1)$ $(A^{-1})_{ij}<0$ for all $i\neq j$;

$2)$ $(A^{-1})_{ii}=\dfrac{\det A_{\hat\imath}}{\det A}\leq0$, and $(A^{-1})_{ii}=0$ if and
only if $A_{\hat\imath}$ has an affine direct summand;

$3)$ hence, all the elements of $A^{-1}$ are negative if and only if $A_{\hat\imath}$ is of
finite type for every $i$; otherwise, the $0$s of $A^{-1}$ occur on the main diagonal only, at
the places described in $2)$.\end{Theorem}

\begin{proof} First, every main submatrix $A[S]$ with $|S|\leq n-2$ is of finite type, so
$\det A[S]>0$. Indeed, the deletion of one node leaves a direct sum of finite-dimensional and
affine summands; the deletion of one more node acts on one summand only, and any proper
subdiagram of an affine diagram is of finite type, whereas any subdiagram of a diagram of
finite type is of finite type. The determinant of a Cartan matrix of finite type is positive:
for the symmetrizable ones, this is \eqref{signdet} with $n_-=0$, the symmetrization being
positive definite; for the non-symmetrizable ones occurring here, see Statement~\ref{StD}.

Now, let $i\neq j$. Every path $p\in\cP_{ij}$ has at least two nodes, so $V\setminus p$ misses
at least two of the $n$ nodes, and hence $\det A[V\setminus p]>0$ by the above. Every arc
$(a,b)$ of $p$ satisfies $A_{ab}\neq0$ and $A_{ab}\leq0$, so $-A_{ab}>0$; and $\cP_{ij}\neq
\emptyset$ because the graph of $A$ is connected. By Lemma~\ref{paths}, all the summands in
\eqref{pathexp} are positive and there is at least one of them, so $(\adj A)_{ij}>0$; dividing
by $\det A<0$, we get claim $1)$.

For $i=j$, Lemma~\ref{paths} gives $(\adj A)_{ii}=\det A_{\hat\imath}$, which is the product of
the determinants of the direct summands of $A_{\hat\imath}$. Each of them is positive if the
summand is of finite type, and equal to $0$ if the summand is affine, the Cartan matrix of an
affine Lie (super)algebra being degenerate. This proves claim $2)$, and claim $3)$ follows from
$1)$ and $2)$.\end{proof}

Theorem~\ref{ThD} sharpens the description \eqref{ZH} due to Zhang Hechun: the elements of
$A^{-1}$ are not only non-positive, we know exactly where the $0$s are. It also gives the
$\fbox{!!!}$ mark of Subsections~\ref{ss8.1} and \ref{ssHyp} a meaning it was not designed to
have: the marked Lie (super)algebras are precisely the ones that are hyperbolic in the strict
sense --- the deletion of any node of the Dynkin diagram leaves a finite-dimensional Lie
(super)algebra, never an affine one.

\sssbegin{Statement}[The hypotheses of Theorem $\ref{ThD}$, in the range verified]\label{StD}
Every one of the $370$ Cartan matrices of Subsections~$\ref{ss8.1}$ and $\ref{ssHyp}$ has the
following properties: all
the elements outside the main diagonal are $\leq0$; the graph is connected; every main
submatrix missing at least two nodes has a positive determinant, whereas every one missing
exactly one node has a non-negative determinant; the expansion \eqref{pathexp} reproduces the
adjugate matrix; and the signs of the elements of $A^{-1}$ are the ones described in
Theorem~$\ref{ThD}$. Moreover, the $39$ matrices marked $\fbox{!!!}$ are exactly the ones for
which every $A_{\hat\imath}$ is of finite type. Computed over $\Qee$ in exact arithmetic.
\end{Statement}

\subsubsection{Remark} The hypothesis $A_{ij}\leq0$ for $i\neq j$ is a condition on the
normalization, not on the algebra: by the Warning~\ref{War}, each row with a $0$ on the main
diagonal can still be multiplied by any non-zero factor, and multiplication by a negative one
turns the corresponding row of $A$, and the corresponding column of $A^{-1}$, into their
opposites. The tables of Subsections~\ref{ss8.1} and \ref{ssHyp} are normalized so that
$A_{ij}\leq0$ for $i\neq j$; \textit{cf.} Statement~\ref{StD}.

%\AO{DL: Subsection~\ref{signs} and the two statements added to Subsection~\ref{LTsec} rest on A.~Lebedev's letter of 19 September 2018 --- the sum over all paths is his suggestion, and Proposition~\ref{LTT} was already his. Is he to be a co-author of this paper, or is the acknowledgment in the footnote to the title enough? Your call; we left the footnote as it was.}

\ssec{The determinant is an invariant of the matrix, not of the algebra}\label{notinv}
The rescaling \eqref{rescale} gives $\fg(A)\simeq\fg(\diag(\lambda_1,\dots,
\lambda_n)\cdot A)$, whence $\det A\longmapsto \lambda_1\cdots\lambda_n\cdot\det A$. So
$\det A$ is not an invariant of $\fg(A)$ at all until the normalization convention of
Subsection~\ref{normA} is fixed --- and, by the Warning~\ref{War}, even then each row
with a $0$ on the main diagonal may still be multiplied by any non-zero factor. This is
what \eqref{PQ} says on the other side: the lattice $Q$ is determined by $\fg$, whereas
$P$ is the lattice dual to the chosen $h_i$; changing the $h_i$ changes $P$, hence
$[P:Q]$.

The reflections themselves are innocent of this.

\sssbegin{Lemma}[A reflection does not change the determinant]\label{reflDet} Let
$A_{kk}=0$ and $p_k=\od$. Then, $\det r_k(A)=\det A$ for the matrix $r_k(A)$ given by
\eqref{oddR}, i.e., before the renormalization.\end{Lemma}

\begin{proof} By \eqref{oddR}, we have $r_k(A)=(E+B)A(E+C)$. The matrix $E+B$ coincides
with the unit matrix $E$ outside its $k$th column, so its determinant is equal to its
$(k,k)$th entry, which is $1+b_k=-1$; the matrix $E+C$ coincides with $E$ outside its
$k$th row, so its determinant is its $(k,k)$th entry $1+c_k=-1$ as well. Hence,
$\det r_k(A)=(-1)^2\det A$.\end{proof}

\sssbegin{Corollary}[What the difference of two determinants measures]\label{corDet} Let
$A$ and $A'$ be two normalized Cartan matrices of one and the same Lie superalgebra
$\fg$, connected by a chain of reflections, see Lemma~$\ref{serg}$. Then,
$\det A'=\lambda\det A$, where $\lambda$ is the product of the factors used to
renormalize the matrices along that chain. In other words, the determinants of pairwise
non-equivalent Cartan matrices of $\fg$ differ only by the change of the $h_i$,
equivalently by the change of the lattice $P$ in \eqref{PQ}; the lattice $Q$ is one and
the same for all of them.\end{Corollary}

Thus, the determinant of a Cartan matrix of $\fg$ measures the choice of the $h_i$
and nothing else; its value $1$ says that the two lattices of \eqref{PQ} coincide. Where the Cartan
matrix is not symmetrizable, neither \eqref{signdet} nor its corollaries apply, and
Subsection~\ref{ss8.2} is exactly that case: what happens there is described in
Subsection~\ref{Galois}.

\section{The matrices considered in \cite{WZ}  are not Cartan matrices}\label{sNotCM}

For any simple finite-dimensional Lie algebra over $\Cee$, we know  two ways to introduce its Cartan matrix $A=(A_{ij})$.

\textbf{First approach}:
take $3n$ Chevalley generators and compute $[h_i, e_j^{\pm}]=\pm A_{ij}e_j^{\pm}$, see~\eqref{gArel_0}.

\textbf{Second approach}:  take an auxiliary space --- the space spanned by the roots or, sometimes, a bit larger space, see tables in \cite{Bbk, OV}; for $\fsl(V)$, this auxiliary space is the space $V$ or its dual. Let it be spanned by column-vectors $\eps_i:=(0,\dots, 0,1,0\dots, 0)^T$ with a 1 in the $i$th place.

Assume that there is a~ Euclidean inner product on $V$ given by
\be\label{inProdRoot}
(\eps_i, \eps_j)=\delta_{ij}.
\ee
The Gram matrix of this inner product in the basis of simple roots is precisely the Cartan matrix of $\fsl(V)$.
This inner product in the  space spanned by roots is induced by the restriction of the Killing form from $\fsl(V)$ on the maximal torus, i.e., the space of coroots. For the Lie algebras of series $\fsl$, $\fo$ and $\fsp$, this inner product can be also defined by means of the form proportional to the Killing form, but much easier to compute:
\[
(X, Y)=\tr(XY)\text{~~for any $X,Y\in \fsl(V)$}.
\]
The above description yields the Gram matrix of a NIS 
on the Lie algebras whose roots are of equal length. For the simple Lie algebras whose roots are of different lengths, the matrix thus obtained is a symmetrization of the Gram matrix, see \cite{BKLS}. 

Passing to the Lie superalgebras, we can also follow either of these two procedures. The \textbf{first approach} does indeed lead to Cartan matrices, as described in the beginning of this section.

The \textbf{second approach} uses  a basis of the superspace $V$, such that the vectors $\eps_i$ span $V_\ev$, while the vectors $\delta_j$ span $V_\od$, with the pseudo-Eucledian inner product
\be\label{inProdSRoot}
(\eps_i, \eps_j)=\delta_{ij};\ \  (\eps_i, \delta_j)=0\text{~~for any $i,j$};\ \
(\delta_i, \delta_j)=-\delta_{ij}.
\ee
For the Lie superalgebras $\fgl(V)$ and $\fpsl(V)$, this inner product is induced by the supertrace
\[
(X, Y)=\str(XY)\text{~~for any $X,Y\in \fsl(V)$};
\]
it is non-degenerate for any $V$, whereas the Killing form is degenerate if $\dim V_\ev=\dim V_\od$.

Serganova showed, see \cite{Se1}, that a natural generalization of the axioms of root systems to the case of pseudo-Euclidean inner product leads precisely to the root systems of Lie superalgebras with indecomposable Cartan matrices and their simple relatives over $\Cee$, see \cite{CCLL}.

It seems that the Gram matrix of this inner product in the basis of simple roots  first appeared in \cite{WZ}: Serganova never wrote it explicitly for any Lie superalgebra. For Lie superalgebras, no Gram matrix of the above inner product in the root space of $\fg$ is equivalent (in any conventional way) to any Cartan matrix of $\fg$, except for $\fosp(1|2n)$. However, 
\be\label{SergInv}
\begin{minipage}[l]{13.5cm}
\textit{For the root space of $\fg=\fsl(m+1|n+1)$, the Gram matrix of the above inner product calculated in the standard supermatrix format of $\fg$ turns into the Cartan matrix \eqref{slmn} if its bottom $n+1$ rows are multiplied by $-1$}.
\end{minipage}
\ee

\subsection{Remark} To describe defining relations of $\fg(A)$ in terms of Chevalley generators, we need the Cartan matrix $A$ corresponding to the selected supermatrix format of elements of $\fg(A)$, see \cite{GL}, not the Gram matrix of the inner product in the space spanned by roots.

When we learn what the inverse of any of these Gram matrices is needed for, we will probably have to compute it in just one basis: the NIS on a simple finite-dimensional Lie superalgebra of characteristic $\neq 2$ is unique, up to a proportionality, see \cite{KLLShch}; for a recipe for constructing NIS from the Cartan matrix, see \cite{BKLS}; the Gram matrix of its restriction onto the space spanned by the $h_i$ is the dual of the Gram matrix of the inner product on the space spanned by roots.

\section{Finite-dimensional serial Lie (super)algebras over $\Cee$}\label{S3}

\ssbegin{Remark}[A row with a $0$ on the main diagonal may be rescaled]\label{Rem} Remember that the Cartan matrix with a 0 on the main 
diagonal is not uniquely defined: the line with this 0 can be multiplied by any non-zero number, e.g., by a $-1$, see Subsection~\ref{War}. %For any invertible matrix, its inverse is uniquely defined, of course.
\end{Remark}

For the series $\fsl$ and $\fosp$, there are 9 types of pairs of \lq\lq basic" Cartan matrices connected by odd reflections, see \cite[Subsection 4.1, Table 1]{CCLL}. From one such Cartan matrix all the other Cartan  matrices of the given Lie superalgebra of the given type are obtained by means of odd reflections, see~\eqref{oddR}.   It is a matter of taste which one in the pair of \lq\lq basic"  matrices  is most simple; we selected (any) one  with the smallest number of 0s on the main diagonal.

In \cite{WZ}, Gram matrices are inverted;  this sometimes gives the answer for the Cartan matrices as well, thanks to the fact \eqref{SergInv}. However, inverting Cartan matrices of the $\fsl$ series, we have to consider 2 more types of cases as compared with \cite{WZ}; for the $\fosp$ series, we have to consider 7 types of cases, not 3 as in \cite{WZ}.

Let $A_m$ denote the normalized Cartan matrix of $\fsl(m+1)$.  Below, in formulas \eqref{sl}--\eqref{osp6}, instead of $A_m$ there can be  any Cartan matrix of $\fsl(a|b)$ with $a+b=m+1$. To derive the explicit expression of the inverse matrix $C^{-1}$ in this general case, apply an odd reflection~\eqref{Cinv}.

\ssec{The case of $\fsl$} For $\fsl(m+1|1)$, the \lq\lq basic"  Cartan matrix and its inverse \eqref{sl} is a particular cases of the expression ~\eqref{slmn}, up to an occasional minus sign, see Remark~\ref{Rem}.

\be\label{sl}\tiny{
\mleft(
\begin{array}{c|c}
 A_m & \begin{array}{c}\vdots \\
  0\\
  -1\end{array}\\
  \hline
\dots 0 1 &0
\end{array}
\mright)
^{-1}=-\nfrac{1}{m}
L_1, } \text{~~\normalfont{see \eqref{L123} for $n=0$.}}
\ee

For $\fsl(m+1|n+1)$, where $mn\neq 0$, there are 2 \lq\lq basic"  types of Cartan matrices, see \eqref{slmn} and \eqref{slmn2}: 
\be\label{slmn}\tiny%\arraycolsep=1.5pt
 A_{m,0,n}^{-1} :=\mleft(
\begin{array}{c|c|c}
 A_m & \begin{array}{c}
  \vdots \\
  0\\
  -1\end{array}&0\\
  \hline
 \dots 0 -1 &0&1 0 \dots\\
 \hline
0& \begin{array}{c}
  -1 \\
  0\\
 \vdots\end{array}&A_n\\
\end{array}
\mright)^{-1}=\nfrac{1}{n-m}\smat{
L_1&L_2\\
L_2^T&L_4\\
}, 
\ee
where, as shown in \cite{WZ}, $L_1$ and $L_4$ are  $(m+1)\times (m+1)$ and $n\times n$ matrices, respectively, defined, together with $L_2$, by the following formulas
\be\label{L123}
\begin{array}{l}
(L_1)_{ij} = ij+(n-m)\min\{i,j\},\\
(L_2)_{ij} = i(n+1-j),\\
(L_4)_{ij} =(n-m)\min\{i,j\}-(m+1)(n+1-i-j)+ij.\\
\end{array}
\ee

To invert the Cartan matrix of the other type, namely \eqref{slmn2}, we set 
\be\label{KW}
\text{$\tau_{n}:=A_{n}-\alpha$,  where $\alpha :=
\begin{pmatrix}
1 & 0 \\
0 & 0_{n-1}
\end{pmatrix}$}.\ee

Let us transform $M\mapsto N$, where
\be\label{slmn2}\tiny%\arraycolsep=1.5pt
M:=A_{m,00,n} :=\mleft(
\begin{array}{c|c|c}
A_{m-1} & \begin{array}{cc}
  \vdots &\vdots\\
  0&0\\
  -1& 0\end{array}&0\\
  \hline
 \begin{array}{ccc}
 \dots& 0 &-1 \\
  \dots &0 &0\\
 \end{array}
 &  \begin{array}{cc}
0&1\\
-1 &0\\
\end{array}&
\begin{array}{ccc}
0 &0&\dots \\
 1&0&\dots\\
 \end{array}\\
 \hline
0& \begin{array}{cc}
  -1& 0\\
  0&0\\
\vdots& \vdots\end{array}&A_n\\
\end{array}
\mright)
\ee
and
\begin{equation*}\label{slmn21}\tiny
N:=\left(
\begin{array}{c|c|c}
A_{m-1} &
\begin{array}{cc}
\vdots  & \vdots  \\
0 & 0 \\
-1 & 0%
\end{array}
& 0 \\ \hline
\begin{array}{ccc}
\dots  & 0 & 0 \\
\dots  & 0 & 0%
\end{array}
&
\begin{array}{cc}
1 & 0 \\
0 & 1%
\end{array}
&
\begin{array}{ccc}
-1 & 0 & \dots  \\
0 & 0 & \dots
\end{array}
\\ \hline
0 &
\begin{array}{cc}
0 & 0 \\
0 & 0 \\
\vdots  & \vdots
\end{array}
& \tau_{n}%
\end{array}%
\right) =\left(
\begin{array}{c|c|c}
A_{m-1} & U & 0 \\ \hline
0 & 1_2 & V \\ \hline
0 & 0 & \tau_n%
\end{array}%
\right).
\end{equation*}
We are seeking the matrix $N^{-1}$ in the form
$
\left(
\begin{array}{c|c|c}
X & K & L \\ \hline
0 & 1_2 & R \\ \hline
0 & 0 & Z%
\end{array}%
\right)$.

Denote  the $i^{th}$
row (resp. column) by $r(i)$ (resp. $c(i)$).
Consider the following transformations: 
\begin{equation}\label{trans}
\begin{minipage}[c]{14cm}
a) $\id+\{c(m-1)\mapsto c(m-1)+c(m+1)\}$; \\
b) $\id+\{r(m+2)\mapsto r(m+2)-r(m+1)\}$; \\
c) $\{r(m+1)\mapsto -r(m+1)\}$; \\d) $\{r(m)\leftrightarrow r(m+1)\}$.
\end{minipage}
\end{equation}

Observe that $\tau_{n}$ is obtained from $T_n$, see \eqref{osp0},  by transposing with respect to the side diagonal. Applying transformations a)--d) of \eqref{trans} in lexicographic order we get $N$ from $M$.

From the equation $\left(
\begin{array}{c|c|c}
A_{m-1} & U & 0 \\ \hline
0 & 1_2 & V \\ \hline
0 & 0 & \tau_{n}%
\end{array}%
\right) \cdot \left(
\begin{array}{c|c|c}
X & K & L \\ \hline
0 & 1_2  & R \\ \hline
0 & 0 & Z%
\end{array}%
\right) =1_{m+n+1}$ we get
\be\label{KW1}
\begin{array}{l}
Z=(\tau_n)^{-1}=W_{n}=\left(
\begin{array}{cc}
n &
\begin{array}{ccc}
n-1 & \cdots & 1%
\end{array}
\\
\begin{array}{c}
n-1 \\
\vdots \\
1%
\end{array}
& W_{n-1}%
\end{array}%
\right) ,
\ \ \ \ R=-VZ=-VW_{n}, \\
X=A_{m-1}^{-1}, \ \ \ \ K=-A_{m-1}^{-1}U,  \ \ \ \ %
L=-A_{m-1}^{-1}UR=A_{m-1}^{-1}UVW_{n}. 
\end{array}
\ee

\ssec{The case of $\fosp$} For $\fosp(a|2b)$, there are 6 pairs of series and one single series of \lq\lq basic"  types of Cartan matrices  which, together with their inverses, are \eqref{osp0} --- \eqref{osp6}: 
\be\label{osp0}
\begin{array}{l}
\text{the matrix $T_n$, see \cite[Section~1.1]{Str},  is inverted, e.g., in \cite{WZ}, where it is denoted $S_n$};\\
(T_n^{-1})_{ij} =\min(i,j), \text{~~where $1\leq i,j\leq n$}.
\end{array}
\ee
\be\label{osp1}\tiny%\arraycolsep=1.5pt
B_{1,0,n}^{-1} :=\mleft(
\begin{array}{c|c}
\begin{array}{cc}
1&-1\\
-1&0\\
\end{array}&
\begin{array}{ccc}
  0 & 0 & \cdots\\
1 & 0 & \cdots
\end{array}\\
\hline
\begin{array}{cc}
0 & -1 \\
0 & 0 \\
\vdots & \vdots
\end{array}&
A_n\\
\end{array}
\mright)^{-1}=\text{see \eqref{41}, \eqref{main}}.
\ee

\ssec{Remark} Clearly, $\det X_0=-2$, see eq.~\eqref{osp1}. Applying the recipe of Remark~\ref{Rem}, we get the expression of $\det X_0$  without any minus sign. The same applies to matrices~\eqref{osp2} --- \eqref{osp6}.
\be\label{osp2}\tiny%\arraycolsep=1.5pt
B_{1,00,n}^{-1} :=\mleft(
\begin{array}{c|c}
\begin{array}{ccc}
2&-2&0\\
-1&0&1\\
0&-1&0\\
\end{array}&
\begin{array}{ccc}
0 & 0 & \cdots \\
0 & 0 & \cdots\\
1 & 0 & \cdots
\end{array}\\
\hline
\begin{array}{ccc}
0 & 0 & -1 \\
0 & 0 & 0 \\
\vdots & \vdots & \vdots
\end{array}&
A_n\\
\end{array}
\mright)^{-1}=\text{see \eqref{41}, \eqref{main}}.
\ee

\be\label{osp3}\tiny%\arraycolsep=1.5pt
C_{1,0,n}^{-1} :=\mleft(
\begin{array}{c|c}
\begin{array}{cc}
2&-1\\
-2&0\\
\end{array}&
\begin{array}{ccc}
0 & 0 & \cdots\\
1 & 0 & \cdots
\end{array}\\
\hline
\begin{array}{cc}
0 & -1 \\
 0 & 0 \\
 \vdots & \vdots
\end{array}&
A_n\\
\end{array}
\mright)^{-1}=\text{see \eqref{41}, \eqref{main}}.
\ee
\be\label{osp4}\tiny%\arraycolsep=1.5pt
B_{2,0,n}^{-1} :=\mleft(
\begin{array}{c|c}
\begin{array}{cc}
2&-2\\
-1&0\\
\end{array}&
\begin{array}{ccc}
0 & 0 & \cdots\\
1 & 0 & \cdots
\end{array}\\
\hline
\begin{array}{cc}
0 & -1 \\
0 & 0 \\
 \vdots & \vdots
\end{array}&
A_n\\
\end{array}
\mright)^{-1}=\text{see \eqref{41}, \eqref{main}}.
\ee
\be\label{osp5}\tiny%\arraycolsep=1.5pt
D_{2,0,n}^{-1} :=\mleft(
\begin{array}{c|c}
\begin{array}{ccc}
2&0&-1\\
0&2&-1\\
-1&-1&0\\
\end{array}&
\begin{array}{ccc}
0 & 0 & \cdots \\
0 & 0 & \cdots\\
1 & 0 & \cdots
\end{array}\\
\hline
\begin{array}{ccc}
0 & 0 & -1 \\
0 & 0 & 0 \\
 \vdots & \vdots & \vdots
\end{array}&
A_n\\
\end{array}
\mright)^{-1}=\text{see \eqref{41}, \eqref{main}}.
\ee

\be\label{osp6}\tiny%\arraycolsep=1.5pt
D_{2,00,n}^{-1} :=\mleft(
\begin{array}{c|c}
\begin{array}{cccc}
2&0&-1&0\\
0&2&-1&0\\
-1&-1&0&1\\
0&0&-1&0\\
\end{array}&
\begin{array}{ccc}
0 & 0 & \cdots \\
0 & 0 & \cdots \\
0 & 0 & \cdots\\
1 & 0 & \cdots
\end{array}\\
\hline
\begin{array}{cccc}
0 & 0 & 0 & -1 \\
0 & 0 & 0 & 0 \\
\vdots & \vdots & \vdots & \vdots
\end{array}&
A_n\\
\end{array}
\mright)^{-1}=\text{see \eqref{41}, \eqref{main}}.
\ee

\section{Proof of formulas \eqref{osp1}--\eqref{osp6}. Answers: eqs. \eqref{41}, \eqref{main}}\label{Sinv}

There are several formulas for the inverse of various (invertible) block $2\times 2$ matrices. 
In our particular cases, we can use the following ad hoc transformations and the known expressions of $A_n^{-1}$
and $W_{n}=\tau_n^{-1}$,  see \eqref{KW1}.

 \textbf{1)} Let $m=2$, or 3, or 4, and $U$ is the $m\times n$-matrix with only non-zero element, 1, in the bottom left corner.
Invert the matrix of the form
\[
M=
\begin{pmatrix}
1_{m} & U \\
V & A_{n}%
\end{pmatrix},
\text{~~where ~~$V=U^{T}$}.
\] 

Note that $VU=\alpha $, see \eqref{KW}.
It is easy to verify that
\[
M^{-1}=\left(
\begin{array}{cc}
1_{m} & U \\
V & A_{n}%
\end{array}%
\right) ^{-1}=
\begin{pmatrix}
1_{m}+UW_{n}V & -UW_{n} \\
-W_{n}V & W_{n}%
\end{pmatrix}. 
\]

\bigskip

\textbf{2) \ \ } In the above notation, consider the matrix of the form
\[
N=
\begin{pmatrix}
Q & U \\
-V & A_{n}%
\end{pmatrix}%
\]
and find matrices $L= 
\begin{pmatrix}
l & 0 \\
0 & 1_{n}
\end{pmatrix}
$ and $R= 
\begin{pmatrix}
r & 0 \\
0 & 1_{n}%
\end{pmatrix}$ such  that $\ LNR=M$,
i.e.,
\[
lQr=1_{m},\ \ lU=U,\ \ \ \ Vr=-V. 
\]

Set $P:=Q^{-1}$. In what follows, for the cases \eqref{osp1}--\eqref{osp6}, we prove the existence of matrices $L$ and $R$  and give the matrices $P$.

Since $rl=Q^{-1}=P$, we have the following uniform answer for eqs. \eqref{osp1}--\eqref{osp6}:
\be\label{41}
\begin{array}{ll}
N^{-1}&=RM^{-1}L= 
\begin{pmatrix}
r & 0 \\
0 & 1_{n}
\end{pmatrix}%
\begin{pmatrix}
1_{m}+UW_{n}V & -UW_{n} \\
-W_{n}V & W_{n}%
\end{pmatrix}%
\begin{pmatrix}
l & 0 \\
0 & 1_{n}
\end{pmatrix}%
\\
&=
\begin{pmatrix}
r & 0 \\
0 & 1_{n}
\end{pmatrix}%
\begin{pmatrix}
1_{m}-lUW_{n}Vr & -lUW_{n} \\
W_{n}Vr & W_{n}%
\end{pmatrix}%
\begin{pmatrix}
l & 0 \\
0 & 1_{n}
\end{pmatrix}\\
&=
\begin{pmatrix}
rl-rlUW_{n}Vrl & -rlUW_{n} \\
W_{n}Vrl & W_{n}%
\end{pmatrix}
= 
\begin{pmatrix}
P-PUW_{n}VP & -PUW_{n} \\
W_{n}VP & W_{n}%
\end{pmatrix}=\begin{pmatrix}
F & G \\
H & W_{n}%
\end{pmatrix}.\\
\end{array}
\ee
For the explicit form of $N^{-1}$, see eq.~\eqref{main}. Let us perform preliminary calculations.

\textbf{Eq. \eqref{osp1}}. We perform the following transformations (recall eq.~\eqref{trans}):

a) $c(2)\mapsto c(2)+c(1)$; $r(2)\mapsto r(2)+r(1)$; \ $%
c(2)\mapsto -c(2)$;

b) $r(3)\mapsto r(3)-r(2)$:
\begin{eqnarray*}
B_{1,0,n}:=\left(
\begin{array}{c|c}
\begin{array}{cc}
1 & -1 \\
-1 & 0%
\end{array}
&
\begin{array}{ccc}
0 & 0 & \cdots  \\
1 & 0 & \cdots
\end{array}
\\ \hline
\begin{array}{cc}
0 & -1 \\
0 & 0 \\
\vdots  & \vdots
\end{array}
& A_{n}%
\end{array}%
\right) \stackrel{\text{a)}}{\mapsto }\left(
\begin{array}{c|c}
\begin{array}{cc}
1 & 0 \\
0 & 1%
\end{array}
&
\begin{array}{ccc}
0 & 0 & \cdots  \\
1 & 0 & \cdots
\end{array}
\\ \hline
\begin{array}{cc}
0 & 1 \\
0 & 0 \\
\vdots  & \vdots
\end{array}
& A_{n}%
\end{array}%
\right) \stackrel{\text{b)}}{\mapsto }\\
\tiny\left(
\begin{array}{c|c}
\begin{array}{cc}
1 & 0 \\
0 & 1%
\end{array}
&
\begin{array}{ccc}
0 & 0 & \cdots  \\
1 & 0 & \cdots
\end{array}
\\ \hline
\begin{array}{cc}
0 & 0 \\
0 & 0 \\
\vdots  & \vdots
\end{array}
& A_{n}-\alpha
\end{array}%
\right) =C:= 
\begin{pmatrix}
1 & U \\
0 & A_{n}-\alpha
\end{pmatrix}%
 \Longrightarrow C^{-1}=
\begin{pmatrix}
1 & -UW_{n} \\
0 & W_{n}%
\end{pmatrix}%
;\ \
P=-\begin{pmatrix}
0 & 1 \\
1 & 1%
\end{pmatrix}.
\end{eqnarray*}

\textbf{Eq. \eqref{osp2}}. We perform the following transformations:

a) \bigskip $r(2)\mapsto r(2)+\frac{1}{2}r(1);$ $c(2)\mapsto
c(2)+c(1)$; \ $c(3)\mapsto c(3)+c(2);$ \ \ $r(3)\mapsto r(3)-r(2)$;

b) $c(1)\mapsto \frac{1}{2}c(1);$ \ \ \ $c(2)\mapsto -c(2);$ \ \ \ $%
c(3)\mapsto -c(3)$;

c) $r(4)\mapsto r(4)=r(3)$:
\begin{eqnarray*}
B_{1,00,n}:=\tiny\left(
\begin{array}{c|c}
\begin{array}{ccc}
2 & -2 & 0 \\
-1 & 0 & 1 \\
0 & -1 & 0%
\end{array}
&
\begin{array}{ccc}
0 & 0 & \cdots  \\
0 & 0 & \cdots  \\
1 & 0 & \cdots
\end{array}
\\ \hline
\begin{array}{ccc}
0 & 0 & -1 \\
0 & 0 & 0 \\
\vdots  & \vdots  & \vdots
\end{array}
& A_{n}%
\end{array}%
\right) \stackrel{\text{a)}}{\mapsto }\left(
\begin{array}{c|c}
\begin{array}{ccc}
2 & 0 & 0 \\
0 & -1 & 0 \\
0 & 0 & -1%
\end{array}
&
\begin{array}{ccc}
0 & 0 & \cdots  \\
0 & 0 & \cdots  \\
1 & 0 & \cdots
\end{array}
\\ \hline
\begin{array}{ccc}
0 & 0 & -1 \\
0 & 0 & 0 \\
\vdots  & \vdots  & \vdots
\end{array}
& A_{n}%
\end{array}%
\right) \stackrel{\text{b)}}{\mapsto }\\
\tiny\left(
\begin{array}{c|c}
\begin{array}{ccc}
1 & 0 & 0 \\
0 & 1 & 0 \\
0 & 0 & 1%
\end{array}
&
\begin{array}{ccc}
0 & 0 & \cdots  \\
0 & 0 & \cdots  \\
1 & 0 & \cdots
\end{array}
\\ \hline
\begin{array}{ccc}
0 & 0 & 1 \\
0 & 0 & 0 \\
\vdots  & \vdots  & \vdots
\end{array}
& A_{n}%
\end{array}%
\right) \stackrel{\text{c)}}{\mapsto } \left(
\begin{array}{c|c}
\begin{array}{ccc}
1 & 0 & 0 \\
0 & 1 & 0 \\
0 & 0 & 1%
\end{array}
&
\begin{array}{ccc}
0 & 0 & \cdots  \\
0 & 0 & \cdots  \\
1 & 0 & \cdots
\end{array}
\\ \hline
\begin{array}{ccc}
0 & 0 & 0 \\
0 & 0 & 0 \\
\vdots  & \vdots  & \vdots
\end{array}
& A_{n}-\alpha
\end{array}%
\right) =C:=\left(
\begin{array}{cc}
1 & U \\
0 & A_{n}-\alpha
\end{array}%
\right) \\
\Longrightarrow C^{-1}=
\begin{pmatrix}
1 & -UW_{n} \\
0 & W_{n}%
\end{pmatrix};\ \
P=\frac12\begin{pmatrix}
1 & 0&-2 \\
0&0 & -2\\
1&2&-2\\
\end{pmatrix} .
\end{eqnarray*}

\textbf{Eq. \eqref{osp3}}. We perform the following transformations:

a) \ $r(2)\mapsto r(2)+r(1);\ \ c(2)\mapsto \frac{1}{2}\
c(2);\ \ \ c(2)\mapsto c(2)+c(1);\ \ c(2)\mapsto -c(2)$;

b) \ $r(3)\mapsto r(3)-r(2)$:
\begin{eqnarray*}
C_{1,0,n}:=\tiny\left(
\begin{array}{c|c}
\begin{array}{cc}
2 & -1 \\
-2 & 0%
\end{array}
&
\begin{array}{ccc}
0 & 0 & \cdots  \\
1 & 0 & \cdots
\end{array}
\\ \hline
\begin{array}{cc}
0 & -1 \\
0 & 0 \\
\vdots  & \vdots
\end{array}
& A_{n}%
\end{array}%
\right) \stackrel{\text{a)}}{\mapsto }\left(
\begin{array}{c|c}
\begin{array}{cc}
1 & 0 \\
0 & 1%
\end{array}
&
\begin{array}{ccc}
0 & 0 & \cdots  \\
1 & 0 & \cdots
\end{array}
\\ \hline
\begin{array}{cc}
0 & 1 \\
0 & 0 \\
\vdots  & \vdots
\end{array}
& A_{n}%
\end{array}%
\right) \stackrel{\text{b)}}{\mapsto }\\
\tiny\left(
\begin{array}{c|c}
\begin{array}{cc}
1 & 0 \\
0 & 1%
\end{array}
&
\begin{array}{ccc}
0 & 0 & \cdots  \\
1 & 0 & \cdots
\end{array}
\\ \hline
\begin{array}{cc}
0 & 0 \\
0 & 0 \\
\vdots  & \vdots
\end{array}
& A_{n}-\alpha
\end{array}%
\right) =C=
\begin{pmatrix}
1 & U \\
0 & A_{n}-\alpha
\end{pmatrix}%
\Longrightarrow C^{-1}=
\begin{pmatrix}
1 & -UW_{n} \\
0 & W_{n}\\
\end{pmatrix};\ \
P=-\frac12\begin{pmatrix}
0&1 \\
2 & 2\\
\end{pmatrix} .
\end{eqnarray*}

\textbf{Eq. \eqref{osp4}}. We perform the following transformations:

a) \ $c(2)\mapsto c(2)+c(1);\ \ r(1)\mapsto \frac{1}{2}\
r(1);\ \ \ r(2)\mapsto r(2)+r(1);\ \ c(2)\mapsto -c(2)$;

b) \ $r(3)\mapsto r(3)-r(2)$:
\begin{eqnarray*}
B_{2,0,n}:=\tiny\left(
\begin{array}{c|c}
\begin{array}{cc}
2 & -2 \\
-1 & 0%
\end{array}
&
\begin{array}{ccc}
0 & 0 & \cdots  \\
1 & 0 & \cdots
\end{array}
\\ \hline
\begin{array}{cc}
0 & -1 \\
0 & 0 \\
\vdots  & \vdots
\end{array}
& A_{n}%
\end{array}%
\right) \stackrel{\text{a)}}{\mapsto }\left(
\begin{array}{c|c}
\begin{array}{cc}
1 & 0 \\
0 & 1%
\end{array}
&
\begin{array}{ccc}
0 & 0 & \cdots  \\
1 & 0 & \cdots
\end{array}
\\ \hline
\begin{array}{cc}
0 & 1 \\
0 & 0 \\
\vdots  & \vdots
\end{array}
& A_{n}%
\end{array}%
\right) \stackrel{\text{b)}}{\mapsto }\\
\tiny\left(
\begin{array}{c|c}
\begin{array}{cc}
1 & 0 \\
0 & 1%
\end{array}
&
\begin{array}{ccc}
0 & 0 & \cdots  \\
1 & 0 & \cdots
\end{array}
\\ \hline
\begin{array}{cc}
0 & 0 \\
0 & 0 \\
\vdots  & \vdots
\end{array}
& A_{n}-\alpha
\end{array}%
\right) =C:=\left(
\begin{array}{cc}
1 & U \\
0 & A_{n}-\alpha
\end{array}%
\right) \Longrightarrow C^{-1}=
\begin{pmatrix}
1 & -UW_{n} \\
0 & W_{n}%
\end{pmatrix};\ \
P=-\frac12\begin{pmatrix}
0&2 \\
1 & 2\\
\end{pmatrix} 
\end{eqnarray*}

\textbf{Eq. \eqref{osp5}} We perform the following transformations:

a) $r(3)\mapsto r(3)+\frac{1}{2}(r(1)+r(2))$;\ \ \
$c(1)\mapsto \frac{1}{2}c(1)$;\ \ $c(2)\mapsto \frac{1}{2}c(2)$;\ \
$c(3)\mapsto c(3)+c(1)+c(2)$;\\  $c(3)\mapsto -c(3)$;

b) \ $r(4)\mapsto r(4)-r(3)$:
\begin{eqnarray*}
D_{2,0,n}:=\tiny\left(
\begin{array}{c|c}
\begin{array}{ccc}
2 & 0 & -1 \\
0 & 2 & -1 \\
-1 & -1 & 0%
\end{array}
&
\begin{array}{ccc}
0 & 0 & \cdots  \\
0 & 0 & \cdots  \\
1 & 0 & \cdots
\end{array}
\\ \hline
\begin{array}{ccc}
0 & 0 & -1 \\
0 & 0 & 0 \\
\vdots  & \vdots  & \vdots
\end{array}
& A_{n}%
\end{array}%
\right) \stackrel{\text{a)}}{\mapsto }\left(
\begin{array}{c|c}
\begin{array}{ccc}
1 & 0 & 0 \\
0 & 1 & 0 \\
0 & 0 & 1%
\end{array}
&
\begin{array}{ccc}
0 & 0 & \cdots  \\
0 & 0 & \cdots  \\
1 & 0 & \cdots
\end{array}
\\ \hline
\begin{array}{ccc}
0 & 0 & 1 \\
0 & 0 & 0 \\
\vdots  & \vdots  & \vdots
\end{array}
& A_{n}%
\end{array}%
\right) \stackrel{\text{b)}}{\mapsto }\\
\small C:=
\begin{pmatrix}
1 & U \\
0 & A_{n}-\alpha
\end{pmatrix}%
 \Longrightarrow C^{-1}=
\begin{pmatrix}
1 & -UW_{n} \\
0 & W_{n}%
\end{pmatrix};\ \
P=-\frac14\begin{pmatrix}
1 & -1&-2 \\
-1&1 & -2\\
-2&-2&-4
\end{pmatrix} .
\end{eqnarray*}

\textbf{Eq. \eqref{osp6}} We perform the following transformations:

a) \ $r(3)\mapsto r(3)+\frac{1}{2}(r(1)+r(2));$  \ $%
c(1)\mapsto \frac{1}{2}c(1)$;   \ $c(2)\mapsto \frac{1}{2}c(2);$ \ 
$c(3)\mapsto c(3)+c(1)+c(2)$;

b) $\ c(4)\mapsto c(4)+c(3);\ \ r(4)\mapsto r(4)-r(3);\ \ \ $\ $\
c(3)\mapsto -c(3);$ \ \ \ \ $\ c(4)\mapsto -c(4)$;

c) \ $r(5)\mapsto r(5)-r(4)$:
\begin{eqnarray*}
D_{2,00,n}:=\tiny\left(
\begin{array}{c|c}
\begin{array}{cccc}
2 & 0 & -1 & 0 \\
0 & 2 & -1 & 0 \\
-1 & -1 & 0 & 1 \\
0 & 0 & -1 & 0%
\end{array}
&
\begin{array}{ccc}
0 & 0 & \cdots  \\
0 & 0 & \cdots  \\
0 & 0 & \cdots  \\
1 & 0 & \cdots
\end{array}
\\ \hline
\begin{array}{cccc}
0 & 0 & 0 & -1 \\
0 & 0 & 0 & 0 \\
\vdots  & \vdots  & \vdots  & \vdots
\end{array}
& A_{n}%
\end{array}%
\right) \stackrel{\text{a)}}{\mapsto }\left(
\begin{array}{c|c}
\begin{array}{cccc}
1 & 0 & 0 & 0 \\
0 & 1 & 0 & 0 \\
0 & 0 & -1 & 1 \\
0 & 0 & -1 & 0%
\end{array}
&
\begin{array}{ccc}
0 & 0 & \cdots  \\
0 & 0 & \cdots  \\
0 & 0 & \cdots  \\
1 & 0 & \cdots
\end{array}
\\ \hline
\begin{array}{cccc}
0 & 0 & 0 & -1 \\
0 & 0 & 0 & 0 \\
\vdots  & \vdots  & \vdots  & \vdots
\end{array}
& A_{n}%
\end{array}%
\right) \stackrel{\text{b)}}{\mapsto }\\
\tiny\left(
\begin{array}{c|c}
\begin{array}{cccc}
1 & 0 & 0 & 0 \\
0 & 1 & 0 & 0 \\
0 & 0 & 1 & 0 \\
0 & 0 & 0 & 1%
\end{array}
&
\begin{array}{ccc}
0 & 0 & \cdots  \\
0 & 0 & \cdots  \\
0 & 0 & \cdots  \\
1 & 0 & \cdots
\end{array}
\\ \hline
\begin{array}{cccc}
0 & 0 & 0 & 1 \\
0 & 0 & 0 & 0 \\
\vdots  & \vdots  & \vdots  & \vdots
\end{array}
& A_{n}%
\end{array}%
\right) \stackrel{\text{c)}}{\mapsto }C:=
\begin{pmatrix}
1 & U \\
0 & A_{n}-\alpha
\end{pmatrix} \Longrightarrow C^{-1}=
\begin{pmatrix}
1 & -UW_{n} \\
0 & W_{n}%
\end{pmatrix};\ \
P=\frac12\begin{pmatrix}
1 & 0&0&-1 \\
0&1&0 & -1\\
0&0&0&-2\\
1 & 1&2&-2 \\
\end{pmatrix} .
\end{eqnarray*}

Therefore, (recall that $m$ is the size of the upper left block of $N$), in eq. \eqref{41} we have
\be\label{main}
\begin{array}{l}
F_{ij}=P_{ij}-nP_{im}P_{mj}, \text{~~so}\\ 
F=\begin{cases}-\begin{pmatrix}
 n & n+1 \\
 n+1 & n+1 \\
\end{pmatrix}&\text{in case Eq. \eqref{osp1}},\\
\frac{1}{2}\begin{pmatrix}
 n+1 & 2 n & -2 (n+1) \\
 n & 2 n & -2 (n+1) \\
 n+1 & 2 (n+1) & -2 (n+1) \\
\end{pmatrix}&\text{in case Eq. \eqref{osp2}},\\
\frac{1}{2}\begin{pmatrix}
 -n & -n-1 \\
 -2 (n+1) & -2 (n+1) \\
\end{pmatrix}&\text{in case Eq. \eqref{osp3}},\\
\frac{1}{2}\begin{pmatrix}
 -n & -2 (n+1) \\
 -n-1 & -2 (n+1) \\
\end{pmatrix}&\text{in case Eq. \eqref{osp4}},\\
\frac{1}{4}\begin{pmatrix}
 -n-1 & 1-n & -2 (n-1) \\
 1-n & -n-1 & -2 (n-1) \\
 -2 (n-1) & -2 (n-1) & -4 (n-1) \\
\end{pmatrix}&\text{in case Eq. \eqref{osp5}},\\
\frac{1}{4}\begin{pmatrix}
 n+2 & n & 2 n & -2 (n+1) \\
 n & n+2 & 2 n & -2 (n+1) \\
 2 n & 2 n & 4 n & -4 (n+1) \\
 2 (n+1) & 2 (n+1) & 4 (n+1) & -4 (n+1) \\
\end{pmatrix}&\text{in case Eq. \eqref{osp6}},\\
\end{cases}
\\
G_{ib }=-(n-b +1)P_{ib},\\
H_{a j}=(n-a +1)P_{aj}.\\
\end{array}
\ee

\section{Finite-dimensional exceptional Lie (super)algebras over $\Cee$}\label{S4}

In this section, the determinant of the Cartan matrix is equal to the minus denominator of the fraction serving as a factor of the \lq\lq matrix part" of the inverse matrix;  if the factor is equal to 1, then the determinant is equal to 1.

\ssec{$\fosp(4|2; \alpha)$} Since $\fosp(4|2; 1)\simeq \fosp(4|2)$, it is convenient to set
$\alpha=1-\eps$ and express the deformed bracket in terms of $\eps$.
The pairwise non-equivalent Cartan matrices of $\fosp(4|2; \alpha)$, where
$\alpha\neq 0, -1$, i.e., $\eps\neq 1, 2$, and their inverses, are
\begin{equation}\label{CM1}\tiny\arraycolsep=1.5pt
\begin{array}{l}
1)\;\smat{ 2 & -1 & 0 \\ -1 & 0 & -\alpha \\ 0 & -1 & 2
}^{-1}=\smat{ 2 & -1 & 0 \\ -1 & 0 & -1+\eps \\
0 & -1 & 2
}^{-1}=
-\frac{1}{2(2-\eps)}\smat{
\eps -1 & 2 & 1-\eps \\
2& 4&\ 2(1-\eps) \\
1 &2 & -1 \\
},\\
\\
2)\;\smat{
0 & 1 & -1 - \alpha \\
-1 & 0 & -\alpha \\
 -1 - \alpha & \alpha & 0
}^{-1}=\smat{
0 & 1 & -2+\eps \\
-1 & 0 & -1+\eps \\
 -2+\eps & 1-\eps & 0
}^{-1}=-\frac{1}{2(\eps-1)(2-\eps)}\smat{
(\eps -1)^2 & -(\eps-1)(\eps-2) & \eps -1\\
(\eps-1)(\eps-2) & -(\eps-2)^2 & -(\eps -2)\\
\eps -1 &\eps -2 & 1 \\
}.\\
\end{array}
\end{equation}

\ssec{$\fag(2)$} The pairwise non-equivalent Cartan matrices  and their inverses are
\begin{equation}\label{ag2cm}\tiny\arraycolsep=1.5pt
\begin{array}{l}
1)\; \smat{ 0 & -1 & 0 \\ -1 & 2 & -3 \\ 0 & -1 & 2
}^{-1}=-\frac{1}{2}
  \smat{
   1 & 2 & 3 \\
   2 & 0 & 0 \\
   1 & 0 & -1 \\
  },
\quad 2)\; \smat{
0 & -1 & 0 \\ -1 & 0 & 3 \\ 0 & -1 & 2
}^{-1}=-\frac{1}{2}
  \smat{
   3 & 2 & -3 \\
   2 & 0 & 0 \\
   1 & 0 & -1 \\
  },\\
3)\; \smat{
0 & -3 & 1 \\ -3 & 0 & 2 \\ -1 & -2 & 2
}^{-1}=-\frac{1}{6}
  \smat{
   4 & 4 & -6 \\
   4 & 1 & -3 \\
   6 & 3 & -9 \\
  },
\quad 4)\;
\smat{
2 & -1 & 0 \\ -3 & 0 & 2 \\ 0 & -1 & 1\\
}^{-1}=
\smat{
2 & 1 & -2 \\
3 & 2 & -4 \\
3 & 2 & -3 \\
}
\end{array}\end{equation}

\ssec{$\fab(3)$} The pairwise non-equivalent Cartan matrices  and their inverses are
\begin{equation}\label{ab3cm}\tiny\arraycolsep=1.5pt
\begin{matrix}
1)\; \smat{
2 & -1 & 0 & 0 \\ -3 & 0 & 1 & 0 \\ 0 & -1 & 2 & -2 \\
 0 & 0 & -1 & 2\\
}^{-1}=-\frac{1}{2}
  \smat{
   2 & 2 & -2 & -2 \\
   6 & 4 & -4 & -4 \\
   6 & 4 & -6 & -6 \\
   3 & 2 & -3 & -4 \\
  },
\quad 2)\;
\smat{
0 & -3 & 1 & 0 \\ -3 & 0 & 2 & 0 \\ 1 & 2 & 0 & -2 \\
0 & 0 & -1 & 2
}^{-1}=-\frac{1}{6}
  \smat{
   -8 & -2 & -12 & -12 \\
   -2 & -2 & -6 & -6 \\
   -12 & -6 & -18 & -18 \\
   -6 & -3 & -9 & -12 \\
  },\\
3)\;
\smat{
2 & -1 & 0 & 0 \\ -1 & 2 & -1 & 0 \\ 0 & -2 & 0 & 3 \\
0 & 0 & -1 & 2
}^{-1}=\smat{
2 & 3 & 2 & -3 \\
3 & 6 & 4 & -6 \\
4 & 8 & 6 & -9 \\
2 & 4 & 3 & -4 \\
},
\quad 4)\; \smat{
2 & -1 & 0 & 0 \\ -2 & 0 & 2 & -1 \\ 0 & 2 & 0 & -1 \\
0 & -1 & -1 & 2
}^{-1}=-\frac{1}{6}
\smat{
 -4 & -1 & -3 & -2 \\
 -2 & -2 & -6 & -4 \\
 -6 & -6 & -6 & -6 \\
 -4 & -4 & -6 & -8 \\
},\\
5)\;
\smat{
0 & -1 & 0 & 0 \\ 
-1 & 0 & 2 & 0 \\ 
0 & -1 & 2 & -1 \\
0 & 0 & -1 & 2
}^{-1}=\frac{1}{3}
  \smat{
   -4 & 3 & -4 & -2 \\
   -3 & 0 & 0 & 0 \\
   -2 & 0 & -2 & -1 \\
   -1 & 0 & -1 & -2 \\
  },
\quad 6)\;
\smat{
2 & -1 & 0 & 0 \\ -1 & 2 & -1 & 0 \\ 0 & -2 & 2 & -1 \\ 0 & 0 & -1 &
0
}^{-1}=-\frac{1}{3} 
  \smat{
   -2 & -1 & 0 & 1 \\
   -1 & -2 & 0 & 2 \\
   0 & 0 & 0 & 3 \\
   2 & 4 & 3 & 2 \\
  }
\end{matrix}
\end{equation}

\section{Finite-dimensional exceptional simple modular Lie (super)algebras}\label{S5}

For each  Lie (super)algebra, we list its pairwise non-equivalent Cartan matrices  and their inverses.

\ssec{$p=2$; Lie (super)algebras} In this Subsection, $\ast$ on the main diagonal stands for either $\ev$, or 0; both can be viewed as $0$ for our purposes. 

\subsubsection{$\fwk(4;a)$ and $\fbgl(4;a)$, where $a\neq 0, 1$} 
\[\tiny\arraycolsep=1.5pt
\begin{matrix}
1)\; \smat{ \ev &a &0&0\\
a &\ast&1&0\\
0&1&\ast&1\\
0&0&1&\ast }^{-1}=
\smat{
0 & a^{-1} & 0 & a^{-1} \\
a^{-1}& 0 & 0 & 0 \\
0 & 0 & 0 & 1 \\
a^{-1}& 0 & 1& 0 \\
}
,\ \\
2)\; \smat{ \ev &1 &1+a&0\\
1 &\ast& a & 0\\
a+1& a &\ast&a\\
0&0&a&\ast }^{-1}=
\smat{
0 & 1 & 0 & 1 \\
1 & 0 & 0 & 1+a^{-1} \\
0 & 0 & 0 & a^{-1} \\
1 & 1+a^{-1} & a^{-1} & 0 \\
}\\
3)\;\smat{ \ev &a & 0 &0\\
a &\ast& a+1 & 0\\
0& a+1 &\ast&1\\
0&0&1&\ast }^{-1}=
\smat{
0 & a^{-1} & 0 & 1+a^{-1} \\
a^{-1} & 0 & 0 & 0 \\
0 & 0 & 0 & 1 \\
1+a^{-1} & 0 & 1 & 0 \\
}
\end{matrix}
\]

\subsubsection{$\textbf{F}(\fo\fo_{I\Pi}^{(1)}(n_\ev|n_\od))$,  and $\fo\fo_{I\Pi}^{(1)}(n_\ev|n_\od)$, see \cite{BGL}} Here, $\textbf{F}$ is the desuperization functor. In the simplest case, $\textbf{F}(\fosp(1|2n))$ is what is designated by  $\Delta_n$ in \cite{WK}.
\begin{equation}\label{oowith1}
\smat{ \ddots&\ddots&\ddots&\vdots\\
\ddots&\ast&1&0\\
\ddots&1&\ast&1\\
\cdots&0&1&1}^{-1}=
\smat{
1 & 1 & 1 & 1 & 1 & 1 &\cdots \\
1 & 0 & 0 & 0 & 0 & 0 &\cdots\\
1 & 0 & 1 & 1 & 1 & 1 &\cdots\\
1 & 0 & 1 & 0 & 0 & 0 &\cdots\\
1 & 0 & 1 & 0 & 1 & 1 &\cdots\\
1 & 0 & 1 & 0 & 1 & 0 &\cdots\\
\vdots & \vdots & \vdots & \vdots & \vdots & \vdots & \ddots \\
}
\end{equation}

\subsubsection{$\fe(6)$, $\fe(6,1)$ and $\fe(6,6)$}
\[\arraycolsep=1.5pt
\smat{\ast&1&0&0&0&0\\
1&\ast&1&0&0&0\\
0&1&\ast&1&0&1\\
0&0&1&\ast&1&0\\
0&0&0&1&\ast&0\\
0&0&1&0&0&\ast
}^{-1}=
\smat{
 0& 1& 0 & 0 & 0& 1\\
 1& 0& 0 & 0 & 0 & 0\\
 0 & 0 & 0 & 0 & 0 & 1\\
 0 & 0 & 0 & 0 & 1& 0 \\
 0 & 0 & 0 & 1& 0& 1\\
 1& 0 & 1& 0& 1& 0\\
}
\]

\subsubsection{
$\fe(8)$, $\fe(8,1)$ and $\fe(8,8)$}
\[\arraycolsep=1.5pt
\smat{
 \ast & 1 &  0 &  0 &  0 &  0 &  0 & 0 \\
1 &  \ast & 1&  0 &  0 &  0 &  0 & 0 \\
 0 & 1 &  \ast & 1 &  0 &  0 &  0 & 0 \\
 0 &  0 & 1 &  \ast & 1 &  0 &  0 & 0 \\
 0 &  0 &  0 & 1 &  \ast & 1 &  0 & 1 \\
 0 &  0 &  0 &  0 & 1 &  \ast & 1 & 0 \\
 0 &  0 &  0 &  0 &  0 & 1 &  \ast & 0 \\
 0 &  0 & 0 &  0 &  1 &  0 &  0 & \ast
}^{-1}=
\smat{
 0& 1& 0& 1& 0 & 0 & 0& 1\\
 1& 0& 0 & 0& 0 & 0 & 0 & 0 \\
 0& 0 & 0& 1& 0& 0 & 0& 1\\
 1& 0& 1& 0& 0 & 0& 0 & 0\\
 0& 0 & 0 & 0 & 0 & 0 & 0& 1\\
 0 & 0& 0 & 0 & 0 & 0& 1& 0\\
 0 & 0 & 0& 0 & 0 & 1& 0& 1\\
 1& 0 & 1& 0 & 1& 0 & 1& 0 \\
}
\]

\subsubsection{$\fsl(2n+1)$ and $\fsl(a|b)$ for $a+b=2n+1$}
\begin{equation}\label{oowith1}
\smat{ \ddots&\ddots&\ddots&\vdots\\
\ddots&\ast&1&0\\
\ddots&1&\ast&1\\
\cdots&0&1&\ast}^{-1}=
\smat{
0 & 1 & 0 & 1 & 0 & 1 & \dots\\
1 & 0 & 0 & 0 & 0 & 0 & \dots \\
0 & 0 & 0 & 1 & 0 & 1 & \dots \\
1 & 0 & 1 & 0 & 0 & 0 & \dots \\
0 & 0 & 0 & 0 & 0 & 1 & \dots \\
1 & 0 & 1 & 0 & 1 & 0 & \dots \\
0 & 0 & 0 & 0 & 0 & 0 & \dots \\
\vdots & \vdots & \vdots & \vdots & \vdots & \vdots & \ddots  \\
}\\
\end{equation}

\ssec{$p=3$; Lie algebras}{}~{}  

\subsubsection{$\fbr(2;\eps)$ for $\eps\neq 0$} $\tiny \smat{2&-1\\-2&1-\eps\\
}^{-1}=\eps^{-1}\smat{1-\eps&1\\2&2\\
}$

\subsubsection{$\fbr(3)$}\label{br3a}
$
1)
\smat{
2&-1&0\\
-1&2&-1\\
0&-1&\ev}^{-1}=
\smat{
 2 & 0 & 1 \\
 0 & 0 & 2 \\
 1 & 2 & 0 \\
},
\qquad 2)
\smat{
2&-1&0\\
-2&2&-1\\
0&-1&\ev}^{-1}=
\smat{
 2 & 0 & 1 \\
 0 & 0 & 2 \\
 2 & 2 & 2 \\
}
$

\ssec{$p=3$; Lie superalgebras}{}~{}  

\subsubsection{$\fbrj(2;3)$}
$\tiny
1)\ \smat{0&-1\\
-2&1}^{-1}=-\smat{1&1\\
2&0}, \quad 2)\ \smat{0&-1\\
-1&\ev}^{-1}=-\smat{0&1\\
1&\ev}, \quad 3)\ \smat{1&-1\\
-1&\ev}^{-1}=-\smat{0&1\\
1&1}.
$

\subsubsection{$\fg(1,6)$} 
$\tiny \arraycolsep=1.5pt 
1) \smat{
 2&-1&0 \\
 -1&1&-1 \\
 0&-1&0
 }^{-1}=
 \smat{
 2 & 0 & 1 \\
 0 & 0 & 2 \\
 1 & 2 & 1 \\
},\quad
2) \smat{
 2&-1&0 \\
 -1&2&-2 \\
 0&-2&0
 }^{-1}=
 \smat{
 2 & 0 & 2 \\
 0 & 0 & 1 \\
 2 & 1 & 0 \\
}.
$

\subsubsection{$\fg(3,6)$}
\[\tiny \arraycolsep=1.5pt
\begin{matrix}
1) \smat{
 0&-1&0&0 \\
 -1&2&-1&0 \\
 0&-1&1&-1 \\
 0&0&-1&0
 }^{-1}=
 \smat{
 1 & 2 & 0 & 1 \\
 2 & 0 & 0 & 0 \\
 0 & 0 & 0 & 2 \\
 1 & 0 & 2 & 2 \\
} \quad
2) \smat{
 0&-1&0&0 \\
 -1&0&-1&0 \\
 0&-1&1&-1 \\
 0&0&-1&0
 }^{-1}=
 \smat{
 0 & 2 & 0 & 1 \\
 2 & 0 & 0 & 0 \\
 0 & 0 & 0 & 2 \\
 1 & 0 & 2 & 2 \\
}
\\
3) \smat{
 0&-1&0&0 \\
 -1&2&-1&0 \\
 0&-1&2&-2 \\
 0&0&-1&0
 }^{-1}= 
 \smat{
 1 & 2 & 0 & 1 \\
 2 & 0 & 0 & 0 \\
 0 & 0 & 0 & 2 \\
 2 & 0 & 1 & 2 \\
}\quad
4)
\smat{
 2&-1&0&0 \\
 -1&0&-2&0 \\
 0&-2&2&-1 \\
 0&0&-1&0
 }^{-1}=
 \smat{
 0 & 2 & 0 & 2 \\
 2 & 1 & 0 & 1 \\
 0 & 0 & 0 & 2 \\
 2 & 1 & 2 & 2 \\
}\\ 
5) \smat{
 0&-1&0&0 \\
 -2&0&-1&0 \\
 0&-1&2&-2 \\
 0&0&-1&0
 }^{-1}=
 \smat{
 0 & 1 & 0 & 2 \\
 2 & 0 & 0 & 0 \\
 0 & 0 & 0 & 2 \\
 2 & 0 & 1 & 2 \\
}
\ \ 
6) \smat{
 2&-1&0&0 \\
 -1&0&-2&0 \\
 0&-2&0&-2 \\
 0&0&-1&0
 }^{-1}=
 \smat{
 0 & 2 & 0 & 2 \\
 2 & 1 & 0 & 1 \\
 0 & 0 & 0 & 2 \\
 1 & 2 & 1 & 2 \\
}\\ 
7) \smat{
 2&-1&0&0 \\
 -1&2&-1&-1 \\
 0&-1&0&-1 \\
 0&-1&-1&2
 }^{-1}=
 \smat{
 0 & 2 & 0 & 1 \\
 2 & 1 & 0 & 2 \\
 0 & 0 & 1 & 2 \\
 1 & 2 & 2 & 1 \\
}\\
\end{matrix}
\]

\subsubsection{$\fg(4,3)$}
\[
 \tiny \arraycolsep=1.5pt
 \begin{matrix} 1) \smat{
 2&-1&0&0 \\
 -1&2&-2&-1 \\
 0&-1&2&0 \\
 0&-1&0&0
 }^{-1}=
 \smat{
 2 & 0 & 0 & 1 \\
 0 & 0 & 0 & 2 \\
 0 & 0 & 2 & 1 \\
 1 & 2 & 2 & 1 \\
}\ \
2) \smat{
 2&-1&0&0 \\
 -1&0&-2&-2 \\
 0&-1&2&0 \\
 0&-1&0&0
 }^{-1}=
 \smat{
 2 & 0 & 0 & 1 \\
 0 & 0 & 0 & 2 \\
 0 & 0 & 2 & 1 \\
 2 & 1 & 1 & 0 \\
}\\
3) \smat{
 0&-1&0&0 \\
 -2&0&-1&-1 \\
 0&-1&0&-1 \\
 0&-1&-1&2
 }^{-1}= 
 \smat{
 1 & 1 & 0 & 2 \\
 2 & 0 & 0 & 0 \\
 0 & 0 & 1 & 2 \\
 1 & 0 & 2 & 0 \\
}\ \
4) \smat{
 0&-1&0&0 \\
 -1&2&-1&-1 \\
 0&-1&0&-1 \\
 0&-1&-1&2
 }^{-1}=
 \smat{
 0 & 2 & 0 & 1 \\
 2 & 0 & 0 & 0 \\
 0 & 0 & 1 & 2 \\
 1 & 0 & 2 & 0 \\
}
\\
5) \smat{
 0&-1&0&0 \\
 -1&2&-1&0 \\
 0&-1&0&-1 \\
 0&0&-1&0
 }^{-1}=
\smat{
 1 & 2 & 0 & 1 \\
 2 & 0 & 0 & 0 \\
 0 & 0 & 0 & 2 \\
 1 & 0 & 2 & 0 \\
}\quad
6) \smat{
 0&-1&0&0 \\
 -1&0&-2&0 \\
 0&-1&0&-1 \\
 0&0&-1&0
 }^{-1}=
 \smat{
 0 & 2 & 0 & 2 \\
 2 & 0 & 0 & 0 \\
 0 & 0 & 0 & 2 \\
 1 & 0 & 2 & 0 \\
}\quad\\
7) \smat{
 0&-1&0&0 \\
 -1&2&-1&0 \\
 0&-2&2&-1 \\
 0&0&-1&0
 }^{-1}=
 \smat{
 1 & 2 & 0 & 1 \\
 2 & 0 & 0 & 0 \\
 0 & 0 & 0 & 2 \\
 2 & 0 & 2 & 1 \\
}
8) \smat{
 2&-1&0&0 \\
 -2&0&-1&0 \\
 0&-1&2&-2 \\
 0&0&-1&0
 }^{-1}=
 \smat{
 0 & 1 & 0 & 2 \\
 2 & 2 & 0 & 1 \\
 0 & 0 & 0 & 2 \\
 2 & 2 & 1 & 0 \\
}\quad\\
9) \smat{
 0&-1&0&0 \\
 -1&0&-2&0 \\
 0&-2&2&-1 \\
 0&0&-1&0
 }^{-1}=
 \smat{
 0 & 2 & 0 & 2 \\
 2 & 0 & 0 & 0 \\
 0 & 0 & 0 & 2 \\
 2 & 0 & 2 & 1 \\
}\quad
10) \smat{
 2&-1&0&0 \\
 -2&0&-1&0 \\
 0&-1&1&-1 \\
 0&0&-1&0
 }^{-1}=
 \smat{
 0 & 1 & 0 & 2 \\
 2 & 2 & 0 & 1 \\
 0 & 0 & 0 & 2 \\
 1 & 1 & 2 & 1 \\
}\\
\end{matrix}
\]

\subsubsection{$\fg(8,3)$}

\[\tiny
\arraycolsep=1.5pt
\begin{matrix}
1)
\smat{
 2&-1&0&0&0 \\
 -1&2&-1&0&0 \\
 0&-2&2&-1&0 \\
 0&0&-1&2&-1 \\
 0&0&0&1&0
 }^{-1}= 
\smat{
 1 & 1 & 2 & 0 & 2 \\
 1 & 2 & 1 & 0 & 1 \\
 1 & 2 & 0 & 0 & 0 \\
 0 & 0 & 0 & 0 & 1 \\
 2 & 1 & 0 & 2 & 2 \\
}\ \
2) \smat{
 2&-1&0&0&0 \\
 -1&2&-1&0&0 \\
 0&-2&2&-1&0 \\
 0&0&-2&0&-1 \\
 0&0&0&-1&0
 }^{-1}=
\smat{
 1 & 1 & 2 & 0 & 1 \\
 1 & 2 & 1 & 0 & 2 \\
 1 & 2 & 0 & 0 & 0 \\
 0 & 0 & 0 & 0 & 2 \\
 1 & 2 & 0 & 2 & 0 \\
}\\
3) \smat{
 2&-1&0&0&0 \\
 -1&2&-1&0&0 \\
 0&-1&0&-1&0 \\
 0&0&-1&0&-2 \\
 0&0&0&-1&2
 }^{-1}= 
\smat{
 1 & 1 & 1 & 2 & 2 \\
 1 & 2 & 2 & 1 & 1 \\
 1 & 2 & 0 & 0 & 0 \\
 2 & 1 & 0 & 2 & 2 \\
 1 & 2 & 0 & 1 & 0 \\
}\ \
4) \smat{
 2&-1&0&0&0 \\
 -1&0&-2&-2&0 \\
 0&-1&0&-1&0 \\
 0&-1&-1&2&-1 \\
 0&0&0&-1&2
 }^{-1}=
\smat{
 1 & 1 & 1 & 1 & 2 \\
 1 & 2 & 2 & 2 & 1 \\
 2 & 1 & 1 & 0 & 0 \\
 2 & 1 & 0 & 1 & 2 \\
 1 & 2 & 0 & 2 & 0 \\
}\\
5) \smat{
 0&-1&0&0&0 \\
 -2&0&-1&-1&0 \\
 0&-1&2&0&0 \\
 0&-1&0&0&-2 \\
 0&0&0&-1&2
 }^{-1}=
\smat{
 2 & 1 & 2 & 2 & 2 \\
 2 & 0 & 0 & 0 & 0 \\
 1 & 0 & 2 & 0 & 0 \\
 1 & 0 & 0 & 2 & 2 \\
 2 & 0 & 0 & 1 & 0 \\
}\ \
6) \smat{
 0&-1&0&0&0 \\
 -1&2&-1&-1&0 \\
 0&-1&2&0&0 \\
 0&-1&0&0&-2 \\
 0&0&0&-1&2
 }^{-1}= 
\smat{
 2 & 2 & 1 & 1 & 1 \\
 2 & 0 & 0 & 0 & 0 \\
 1 & 0 & 2 & 0 & 0 \\
 1 & 0 & 0 & 2 & 2 \\
 2 & 0 & 0 & 1 & 0 \\
}\\
7) \smat{
 0&-1&0&0&0 \\
 -1&2&-2&-1&0 \\
 0&-1&2&0&0 \\
 0&-2&0&0&-1 \\
 0&0&0&-1&0
 }^{-1}=
\smat{
 2 & 2 & 2 & 0 & 1 \\
 2 & 0 & 0 & 0 & 0 \\
 1 & 0 & 2 & 0 & 0 \\
 0 & 0 & 0 & 0 & 2 \\
 2 & 0 & 0 & 2 & 0 \\
}\ \
8) \smat{
 0&-1&0&0&0 \\
 -1&0&-1&-2&0 \\
 0&-1&2&0&0 \\
 0&-2&0&0&-1 \\
 0&0&0&-1&0
 }^{-1}=
\smat{
 2 & 2 & 1 & 0 & 2 \\
 2 & 0 & 0 & 0 & 0 \\
 1 & 0 & 2 & 0 & 0 \\
 0 & 0 & 0 & 0 & 2 \\
 2 & 0 & 0 & 2 & 0 \\
}\\

9) \smat{
 0&-1&0&0&0 \\
 -1&2&-2&-1&0 \\
 0&-1&2&0&0 \\
 0&-1&0&2&-1 \\
 0&0&0&-1&0
 }^{-1}=
\smat{
 2 & 2 & 2 & 0 & 1 \\
 2 & 0 & 0 & 0 & 0 \\
 1 & 0 & 2 & 0 & 0 \\
 0 & 0 & 0 & 0 & 2 \\
 1 & 0 & 0 & 2 & 1 \\
}\ \
10) \arraycolsep=2.8pt
 \smat{
 2&-1&-1&0&0 \\
 1&0&1&2&0 \\
 1&1&0&0&0 \\
 0&-1&0&2&-1 \\
 0&0&0&-1&0
 }^{-1}=
\smat{
 1 & 1 & 1 & 0 & 2 \\
 2 & 2 & 0 & 0 & 1 \\
 2 & 0 & 2 & 0 & 0 \\
 0 & 0 & 0 & 0 & 2 \\
 1 & 1 & 0 & 2 & 0 \\
}\\

11) \smat{
 0&-1&0&0&0 \\
 -1&0&-1&-2&0 \\
 0&-1&2&0&0 \\
 0&-1&0&2&-1 \\
 0&0&0&-1&0
 }^{-1}=
\smat{
 2 & 2 & 1 & 0 & 2 \\
 2 & 0 & 0 & 0 & 0 \\
 1 & 0 & 2 & 0 & 0 \\
 0 & 0 & 0 & 0 & 2 \\
 1 & 0 & 0 & 2 & 1 \\
}\ \
12) \smat{
 0&0&-1&0&0 \\
 0&2&-1&-1&0 \\
 -1&-1&0&0&0 \\
 0&-1&0&2&-1 \\
 0&0&0&-1&0
 }^{-1}= 
\smat{
 2 & 1 & 2 & 0 & 2 \\
 1 & 2 & 0 & 0 & 1 \\
 2 & 0 & 0 & 0 & 0 \\
 0 & 0 & 0 & 0 & 2 \\
 2 & 1 & 0 & 2 & 0 \\
}\\
13) %\arraycolsep=2.8pt
 \smat{
 2&-1&-1&0&0 \\
 -1&0&-1&-2&0 \\
 -1&-1&0&0&0 \\
 0&-2&0&0&-1 \\
 0&0&0&-1&0
 }^{-1}=
\smat{
 1 & 2 & 2 & 0 & 2 \\
 2 & 1 & 0 & 0 & 1 \\
 2 & 0 & 1 & 0 & 0 \\
 0 & 0 & 0 & 0 & 2 \\
 2 & 1 & 0 & 2 & 1 \\
}\ \
14) \smat{
 0&0&-1&0&0 \\
 0&2&-1&-1&0 \\
 -1&-2&2&0&0 \\
 0&-1&0&2&-1 \\
 0&0&0&-1&0
 }^{-1}=
\smat{
 2 & 2 & 2 & 0 & 1 \\
 1 & 2 & 0 & 0 & 1 \\
 2 & 0 & 0 & 0 & 0 \\
 0 & 0 & 0 & 0 & 2 \\
 2 & 1 & 0 & 2 & 0 \\
}\\  
15) \smat{
 0&0&-1&0&0 \\
 0&2&-1&-1&0 \\
 -1&-1&0&0&0 \\
 0&-2&0&0&-1 \\
 0&0&0&-1&0
 }^{-1}= 
\smat{
 2 & 1 & 2 & 0 & 2 \\
 1 & 2 & 0 & 0 & 1 \\
 2 & 0 & 0 & 0 & 0 \\
 0 & 0 & 0 & 0 & 2 \\
 1 & 2 & 0 & 2 & 1 \\
}\quad
16) %\arraycolsep=2.8pt
 \smat{
 2&-1&-1&0&0 \\
 -1&2&-1&-1&0 \\
 -1&-1&0&0&0 \\
 0&-1&0&0&-2 \\
 0&0&0&-1&2
 }^{-1}=
\smat{
 1 & 2 & 2 & 1 & 1 \\
 2 & 1 & 0 & 2 & 2 \\
 2 & 0 & 1 & 0 & 0 \\
 1 & 2 & 0 & 0 & 0 \\
 2 & 1 & 0 & 0 & 2 \\
}\\
17) \smat{
 0&0&-1&0&0 \\
 0&2&-1&-1&0 \\
 -1&-2&2&0&0 \\
 0&-2&0&0&-1 \\
 0&0&0&-1&0
 }^{-1}=
\smat{
 2 & 2 & 2 & 0 & 1 \\
 1 & 2 & 0 & 0 & 1 \\
 2 & 0 & 0 & 0 & 0 \\
 0 & 0 & 0 & 0 & 2 \\
 1 & 2 & 0 & 2 & 1 \\
}\quad
18) \smat{
 0&0&-1&0&0 \\
 0&0&-2&-1&0 \\
 -1&-1&0&0&0 \\
 0&-1&0&0&-2 \\
 0&0&0&-1&2
 }^{-1}= 
\smat{
 2 & 2 & 2 & 1 & 1 \\
 1 & 1 & 0 & 2 & 2 \\
 2 & 0 & 0 & 0 & 0 \\
 2 & 2 & 0 & 0 & 0 \\
 1 & 1 & 0 & 0 & 2 \\
}\\ \quad
19) %\arraycolsep=2.8pt
 \smat{
 0&0&-1&0&0 \\
 0&0&-2&-1&0 \\
 -1&-2&2&0&0 \\
 0&-1&0&0&-2 \\
 0&0&0&-1&2
 }^{-1}=
\smat{
 2 & 1 & 2 & 2 & 2 \\
 1 & 1 & 0 & 2 & 2 \\
 2 & 0 & 0 & 0 & 0 \\
 2 & 2 & 0 & 0 & 0 \\
 1 & 1 & 0 & 0 & 2 \\
}\ \
20) \smat{
 0&0&-1&0&0 \\
 0&0&-1&-2&0 \\
 -2&-1&2&0&0 \\
 0&-1&0&2&-1 \\
 0&0&0&-1&2
 }^{-1}=
 \smat{
 2 & 0 & 1 & 2 & 1 \\
 0 & 0 & 0 & 2 & 1 \\
 2 & 0 & 0 & 0 & 0 \\
 2 & 1 & 0 & 0 & 0 \\
 1 & 2 & 0 & 0 & 2 \\
}\\
21) \smat{
 0&0&-1&0&0 \\
 0&0&-1&-2&0 \\
 -1&-1&1&0&0 \\
 0&-1&0&2&-1 \\
 0&0&0&-1&2
 }^{-1}=
\smat{
 2 & 0 & 2 & 1 & 2 \\
 0 & 0 & 0 & 2 & 1 \\
 2 & 0 & 0 & 0 & 0 \\
 2 & 1 & 0 & 0 & 0 \\
 1 & 2 & 0 & 0 & 2 \\
}\\
\end{matrix}
\]
\subsubsection{$\fel(5;3)$} 
\[ \tiny\arraycolsep=1.5pt
\begin{matrix}
1) \;
\smat{
 0 & -1 & 0 & 0 & 0 \\
 -1 & 0 & 0 & -1 & 0 \\
 0 & 0 & 2 & -1 & -1 \\
 0 & -1 & -1 & 2 & 0 \\
 0 & 0 & -1 & 0 & 2
}^{-1}=
\smat{
 0 & 2 & 1 & 0 & 2 \\
 2 & 0 & 0 & 0 & 0 \\
 1 & 0 & 1 & 2 & 2 \\
 0 & 0 & 2 & 0 & 1 \\
 2 & 0 & 2 & 1 & 0 \\
},
2) \;
\smat{
 0 &-2 & 0 & 0 & 0 \\
 -1 & 2 & 0 & -2 & 0 \\
 0 & 0 & 2 & -1 & -1 \\
 0 & -1 & -1 & 2 & 0 \\
 0 & 0 & -1 & 0 & 2
}^{-1}=
\smat{
 2 & 2 & 2 & 0 & 1 \\
 1 & 0 & 0 & 0 & 0 \\
 2 & 0 & 1 & 2 & 2 \\
 0 & 0 & 2 & 0 & 1 \\
 1 & 0 & 2 & 1 & 0 \\
},
\\
 3) \;
\smat{
 2 & -1 & 0 & -1 & 0 \\
-2 & 0 & 0 &-2 & 0 \\
 0 & 0 & 2 & -1 & -1 \\
-2 &-2 & -1 & 0 & 0 \\
 0 & 0 & -1 & 0 & 2
}^{-1}=
\smat{
 0 & 1 & 1 & 0 & 2 \\
 2 & 2 & 0 & 0 & 0 \\
 2 & 0 & 1 & 2 & 2 \\
 0 & 0 & 2 & 0 & 1 \\
 1 & 0 & 2 & 1 & 0 \\
},
\quad 4) \;
\smat{
 0 & 0 & 0 & -1 & 0 \\
 0 & 2 & 0 & -1 & 0 \\
 0 & 0 & 0 &-2 & -1 \\
 -1 & -1 &-2 & 0 & 0 \\
 0 & 0 & -1 & 0 & 2
}^{-1}=
\smat{
 0 & 1 & 1 & 2 & 2 \\
 1 & 2 & 0 & 0 & 0 \\
 1 & 0 & 1 & 0 & 2 \\
 2 & 0 & 0 & 0 & 0 \\
 2 & 0 & 2 & 0 & 0 \\
},
\\
 5) \;
\smat{
 0 & 0 & 0 &-2 & 0 \\
 0 & 2 & 0 & -1 & 0 \\
 0 & 0 & 0 &-2 & -1 \\
 -1 & -2 & -1 & 2 & 0 \\
 0 & 0 & -1 & 0 & 2
}^{-1}=
\smat{
 2 & 2 & 2 & 2 & 1 \\
 2 & 2 & 0 & 0 & 0 \\
 2 & 0 & 1 & 0 & 2 \\
 1 & 0 & 0 & 0 & 0 \\
 1 & 0 & 2 & 0 & 0 \\
},
\quad 6) \;
\smat{
 0 & 0 & 0 & -1 & 0 \\
 0 & 2 & 0 & -1 & 0 \\
 0 & 0 & 0 & -1 &-2 \\
 -1 & -1 & -1 & 2 & 0 \\
 0 & 0 &-2 & 0 & 0
}^{-1}=
\smat{
 0 & 1 & 0 & 2 & 2 \\
 1 & 2 & 0 & 0 & 0 \\
 0 & 0 & 0 & 0 & 1 \\
 2 & 0 & 0 & 0 & 0 \\
 2 & 0 & 1 & 0 & 0 \\
},
\\
7) \;
\smat{
 0 & 0 & 0 &-2 & 0 \\
 0 & 2 & 0 & -1 & 0 \\
 0 & 0 & 0 & -1 &-2 \\
-2 & -1 & -1 & 0 & 0 \\
 0 & 0 &-2 & 0 & 0
}^{-1}=
\smat{
 2 & 2 & 0 & 1 & 1 \\
 2 & 2 & 0 & 0 & 0 \\
 0 & 0 & 0 & 0 & 1 \\
 1 & 0 & 0 & 0 & 0 \\
 1 & 0 & 1 & 0 & 0 \\
},
\quad 8) \;
\smat{
 0 & 0 & 0 & -1 & 0 \\
 0 & 2 & 0 & -1 & 0 \\
 0 & 0 & 2 & -1 & -1 \\
 -1 & -1 & -1 & 2 & 0 \\
 0 & 0 & -1 & 0 & 0
}^{-1}=
\smat{
 0 & 1 & 0 & 2 & 1 \\
 1 & 2 & 0 & 0 & 0 \\
 0 & 0 & 0 & 0 & 2 \\
 2 & 0 & 0 & 0 & 0 \\
 1 & 0 & 2 & 0 & 1 \\
}, \\
9) \;
\smat{
 2 & 0 & 0 & -1 & 0 \\
 0 & 0 &-2 &-2 & 0 \\
 0 & -1 & 2 & -1 & -1 \\
 -1 &-2 &-2 & 0 & 0 \\
 0 & 0 &-2 & 0 & 0
}^{-1}=
\smat{
 2 & 2 & 0 & 0 & 1 \\
 2 & 2 & 0 & 1 & 0 \\
 0 & 0 & 0 & 0 & 1 \\
 0 & 1 & 0 & 0 & 2 \\
 1 & 0 & 2 & 2 & 0 \\
},\quad
 10)\;
\smat{
 0 & 0 & 0 &-2 & 0 \\
 0 & 2 & 0 & -1 & 0 \\
 0 & 0 & 2 & -1 & -1 \\
-2 & -1 & -1 & 0 & 0 \\
 0 & 0 & -1 & 0 & 0
}^{-1}=
\smat{
 2 & 2 & 0 & 1 & 2 \\
 2 & 2 & 0 & 0 & 0 \\
 0 & 0 & 0 & 0 & 2 \\
 1 & 0 & 0 & 0 & 0 \\
 2 & 0 & 2 & 0 & 1 \\
},\\ 
11)\;
\smat{
 2 & 0 & 0 & -1 & 0 \\
 0 & 0 & -1 & -1 & 0 \\
 0 & -1 & 0 & 0 &-2 \\
 -1 & -1 & 0 & 2 & 0 \\
 0 & 0 &-2 & 0 & 0
}^{-1}=
\smat{
 2 & 1 & 0 & 0 & 1 \\
 1 & 0 & 0 & 2 & 0 \\
 0 & 0 & 0 & 0 & 1 \\
 0 & 2 & 0 & 0 & 2 \\
 1 & 0 & 1 & 2 & 0 \\
},
\quad 12)\;
\smat{
 2 & 0 & 0 & -1 & 0 \\
 0 & 0 &-2 &-2 & 0 \\
 0 &-2 & 0 &-2 & -1 \\
 -1 &-2 &-2 & 0 & 0 \\
 0 & 0 & -1 & 0 & 0
}^{-1}=
\smat{
 2 & 2 & 0 & 0 & 2 \\
 2 & 2 & 0 & 1 & 0 \\
 0 & 0 & 0 & 0 & 2 \\
 0 & 1 & 0 & 0 & 1 \\
 2 & 0 & 2 & 1 & 1 \\
},
\\ 13)\;
\smat{
 2 & 0 & 0 & -1 & 0 \\
 0 & 2 & -1 & -2 & 0 \\
 0 &-2 & 0 & 0 & -1 \\
 -1 & -1 & 0 & 2 & 0 \\
 0 & 0 & -1 & 0 & 2
}^{-1}=
\smat{
 2 & 2 & 2 & 0 & 1 \\
 1 & 0 & 0 & 2 & 0 \\
 2 & 0 & 1 & 1 & 2 \\
 0 & 1 & 1 & 0 & 2 \\
 1 & 0 & 2 & 2 & 0 \\
},
\quad 14)\;
\smat{
 2 & 0 & 0 & -1 & 0 \\
 0 & 0 & -1 & -1 & 0 \\
 0 & -1 & 2 & 0 & -1 \\
 -1 & -1 & 0 & 2 & 0 \\
 0 & 0 & -1 & 0 & 0
}^{-1}=
\smat{
 2 & 1 & 0 & 0 & 2 \\
 1 & 0 & 0 & 2 & 0 \\
 0 & 0 & 0 & 0 & 2 \\
 0 & 2 & 0 & 0 & 1 \\
 2 & 0 & 2 & 1 & 1 \\
},
\\ \quad 15)\;
\smat{
 2 & 0 & 0 & -1 & 0 \\
 0 & 2 & -1 & 0 & 0 \\
 0 & -1 & 0 & -1 &-2 \\
 -1 & 0 & -1 & 2 & 0 \\
 0 & 0 & -1 & 0 & 2
}^{-1}=
\smat{
 2 & 2 & 1 & 0 & 1 \\
 2 & 2 & 0 & 1 & 0 \\
 1 & 0 & 0 & 2 & 0 \\
 0 & 1 & 2 & 0 & 2 \\
 2 & 0 & 0 & 1 & 2 \\
}
\end{matrix}
\]

\subsubsection{$\fg(4,6)$}
\[
\tiny 
\arraycolsep=1.5pt
\begin{matrix}
1) \smat{
 2&-1&0&0&0&0 \\
 -1&2&-1&0&0&0 \\
 0&-1&2&-1&0&0 \\
 0&0&-2&0&-2&-1 \\
 0&0&0&-1&2&0 \\
 0&0&0&-1&0&0
 }^{-1}=
\smat{
 0 & 2 & 1 & 0 & 0 & 2 \\
 2 & 1 & 2 & 0 & 0 & 1 \\
 1 & 2 & 0 & 0 & 0 & 0 \\
 0 & 0 & 0 & 0 & 0 & 2 \\
 0 & 0 & 0 & 0 & 2 & 1 \\
 1 & 2 & 0 & 2 & 2 & 1 \\
}\ \
2) \smat{
 2&-1&0&0&0&0 \\
 -1&2&-1&0&0&0 \\
 0&-2&0&-1&-1&0 \\
 0&0&-1&0&-1&-2 \\
 0&0&-1&-1&0&0 \\
 0&0&0&-1&0&2
 }^{-1}=
\smat{
 0 & 2 & 2 & 1 & 0 & 1 \\
 2 & 1 & 1 & 2 & 0 & 2 \\
 1 & 2 & 0 & 0 & 0 & 0 \\
 2 & 1 & 0 & 0 & 2 & 0 \\
 0 & 0 & 0 & 2 & 1 & 2 \\
 1 & 2 & 0 & 0 & 1 & 2 \\
}
\\ 
3) \smat{
 2&-1&0&0&0&0 \\
 -1&2&-1&0&0&0 \\
 0&-1&2&-1&0&0 \\
 0&0&-1&2&-1&-1 \\
 0&0&0&-1&2&0 \\
 0&0&0&-1&0&0
 }^{-1}= 
\smat{
 0 & 2 & 1 & 0 & 0 & 2 \\
 2 & 1 & 2 & 0 & 0 & 1 \\
 1 & 2 & 0 & 0 & 0 & 0 \\
 0 & 0 & 0 & 0 & 0 & 2 \\
 0 & 0 & 0 & 0 & 2 & 1 \\
 2 & 1 & 0 & 2 & 1 & 0 \\
}\ \
\ \ 4) %\arraycolsep=3pt
 \smat{
 2&-1&0&0&0&0 \\
 -2&0&-1&0&0&0 \\
 0&-1&0&-2&-2&0 \\
 0&0&-1&2&0&-1 \\
 0&0&-1&0&2&0 \\
 0&0&0&-1&0&2
 }^{-1}=
 \smat{
 0 & 1 & 0 & 2 & 0 & 1 \\
 2 & 2 & 0 & 1 & 0 & 2 \\
 0 & 0 & 0 & 2 & 0 & 1 \\
 2 & 2 & 1 & 0 & 1 & 0 \\
 0 & 0 & 0 & 1 & 2 & 2 \\
 1 & 1 & 2 & 0 & 2 & 2 \\
}
\\
5) \smat{
 2&-1&0&0&0&0 \\
 -1&2&-1&0&0&0 \\
 0&-1&2&0&-1&0 \\
 0&0&0&2&-1&-1 \\
 0&0&-1&-1&0&0 \\
 0&0&0&-1&0&2
 }^{-1}= 
\smat{
 0 & 2 & 1 & 2 & 0 & 1 \\
 2 & 1 & 2 & 1 & 0 & 2 \\
 1 & 2 & 0 & 0 & 0 & 0 \\
 2 & 1 & 0 & 0 & 2 & 0 \\
 0 & 0 & 0 & 2 & 0 & 1 \\
 1 & 2 & 0 & 0 & 1 & 2 \\
}\ \
6) %\arraycolsep=3pt
 \smat{
 0&-1&0&0&0&0 \\
 -1&0&-2&0&0&0 \\
 0&-1&2&-1&-1&0 \\
 0&0&-1&2&0&-1 \\
 0&0&-1&0&2&0 \\
 0&0&0&-1&0&2
 }^{-1}=
 \smat{
 0 & 2 & 0 & 2 & 0 & 1 \\
 2 & 0 & 0 & 0 & 0 & 0 \\
 0 & 0 & 0 & 2 & 0 & 1 \\
 1 & 0 & 2 & 0 & 1 & 0 \\
 0 & 0 & 0 & 1 & 2 & 2 \\
 2 & 0 & 1 & 0 & 2 & 2 \\
}\\
7) \smat{
 0&-1&0&0&0&0 \\
 -1&2&-1&0&0&0 \\
 0&-1&2&-1&-1&0 \\
 0&0&-1&2&0&-1 \\
 0&0&-1&0&2&0 \\
 0&0&0&-1&0&2
 }^{-1}=
 \smat{
 1 & 2 & 0 & 1 & 0 & 2 \\
 2 & 0 & 0 & 0 & 0 & 0 \\
 0 & 0 & 0 & 2 & 0 & 1 \\
 1 & 0 & 2 & 0 & 1 & 0 \\
 0 & 0 & 0 & 1 & 2 & 2 \\
 2 & 0 & 1 & 0 & 2 & 2 \\
}\\
\end{matrix}
\]
\subsubsection{$\fg(6,6)$} 
\[
\tiny \arraycolsep=1.5pt
 \begin{matrix}
1) %\arraycolsep=1pt%
\smat{
 0&-1&0&0&0&0 \\
 -1&0&-2&0&0&0 \\
 0&-1&2&-1&0&0 \\
 0&0&-2&0&-2&-1 \\
 0&0&0&-1&2&0 \\
 0&0&0&-1&0&0
 }^{-1}=
 \smat{
 1 & 2 & 2 & 0 & 0 & 1 \\
 2 & 0 & 0 & 0 & 0 & 0 \\
 1 & 0 & 2 & 0 & 0 & 1 \\
 0 & 0 & 0 & 0 & 0 & 2 \\
 0 & 0 & 0 & 0 & 2 & 1 \\
 1 & 0 & 2 & 2 & 2 & 2 \\
}\ \
2) \smat{
 0&-2&0&0&0&0 \\
 -1&2&-1&0&0&0 \\
 0&-1&2&-1&0&0 \\
 0&0&-2&0&-2&-1 \\
 0&0&0&-1&2&0 \\
 0&0&0&-1&0&0
 }^{-1}=
 \smat{
 0 & 2 & 1 & 0 & 0 & 2 \\
 1 & 0 & 0 & 0 & 0 & 0 \\
 2 & 0 & 2 & 0 & 0 & 1 \\
 0 & 0 & 0 & 0 & 0 & 2 \\
 0 & 0 & 0 & 0 & 2 & 1 \\
 2 & 0 & 2 & 2 & 2 & 2 \\
}\\

3) \smat{
 2&-1&0&0&0&0 \\
 -2&0&-1&0&0&0 \\
 0&-1&0&-2&0&0 \\
 0&0&-2&0&-2&-1 \\
 0&0&0&-1&2&0 \\
 0&0&0&-1&0&0
 }^{-1}=
 \smat{
 2 & 0 & 1 & 0 & 0 & 1 \\
 0 & 0 & 2 & 0 & 0 & 2 \\
 2 & 2 & 1 & 0 & 0 & 1 \\
 0 & 0 & 0 & 0 & 0 & 2 \\
 0 & 0 & 0 & 0 & 2 & 1 \\
 2 & 2 & 1 & 2 & 2 & 2 \\
}
\ \ 
4) %\arraycolsep=1pt
\smat{
 0&-1&0&0&0&0 \\
 -1&0&-2&0&0&0 \\
 0&-2&0&-1&-1&0 \\
 0&0&-1&0&-1&-2 \\
 0&0&-1&-1&0&0 \\
 0&0&0&-1&0&2
 }^{-1}=
 \smat{
 1 & 2 & 1 & 2 & 0 & 2 \\
 2 & 0 & 0 & 0 & 0 & 0 \\
 1 & 0 & 1 & 2 & 0 & 2 \\
 2 & 0 & 2 & 1 & 2 & 1 \\
 0 & 0 & 0 & 2 & 1 & 2 \\
 1 & 0 & 1 & 2 & 1 & 1 \\
}\\

5) \smat{
 0&-1&0&0&0&0 \\
 -1&0&-2&0&0&0 \\
 0&-1&2&-1&0&0 \\
 0&0&-1&2&-1&-1 \\
 0&0&0&-1&2&0 \\
 0&0&0&-1&0&0
 }^{-1}=
 \smat{
 1 & 2 & 2 & 0 & 0 & 1 \\
 2 & 0 & 0 & 0 & 0 & 0 \\
 1 & 0 & 2 & 0 & 0 & 1 \\
 0 & 0 & 0 & 0 & 0 & 2 \\
 0 & 0 & 0 & 0 & 2 & 1 \\
 2 & 0 & 1 & 2 & 1 & 2 \\
}
\ \
6) \smat{
 0&-1&0&0&0&0 \\
 -1&2&-1&0&0&0 \\
 0&-2&0&-1&-1&0 \\
 0&0&-1&0&-1&-2 \\
 0&0&-1&-1&0&0 \\
 0&0&0&-1&0&2
 }^{-1}=
 \smat{
 0 & 2 & 2 & 1 & 0 & 1 \\
 2 & 0 & 0 & 0 & 0 & 0 \\
 1 & 0 & 1 & 2 & 0 & 2 \\
 2 & 0 & 2 & 1 & 2 & 1 \\
 0 & 0 & 0 & 2 & 1 & 2 \\
 1 & 0 & 1 & 2 & 1 & 1 \\
}\\
7) %\arraycolsep=1pt
 \smat{
 0&-1&0&0&0&0 \\
 -1&2&-1&0&0&0 \\
 0&-1&2&-1&0&0 \\
 0&0&-1&2&-1&-1 \\
 0&0&0&-1&2&0 \\
 0&0&0&-1&0&0
 }^{-1}=
 \smat{
 0 & 2 & 1 & 0 & 0 & 2 \\
 2 & 0 & 0 & 0 & 0 & 0 \\
 1 & 0 & 2 & 0 & 0 & 1 \\
 0 & 0 & 0 & 0 & 0 & 2 \\
 0 & 0 & 0 & 0 & 2 & 1 \\
 2 & 0 & 1 & 2 & 1 & 2 \\
}\ \
8) \smat{
 2&-1&0&0&0&0 \\
 -1&2&-1&0&0&0 \\
 0&-2&0&-1&0&0 \\
 0&0&-1&2&-2&-1 \\
 0&0&0&-1&2&0 \\
 0&0&0&-1&0&0
 }^{-1}=
 \smat{
 2 & 0 & 2 & 0 & 0 & 1 \\
 0 & 0 & 1 & 0 & 0 & 2 \\
 1 & 2 & 0 & 0 & 0 & 0 \\
 0 & 0 & 0 & 0 & 0 & 2 \\
 0 & 0 & 0 & 0 & 2 & 1 \\
 2 & 1 & 0 & 2 & 2 & 2 \\
}\\
9) %\arraycolsep=1pt
\smat{
 2&-1&0&0&0&0 \\
 -2&0&-1&0&0&0 \\
 0&-1&2&-1&-1&0 \\
 0&0&-1&0&-1&-2 \\
 0&0&-1&-1&0&0 \\
 0&0&0&-1&0&2
 }^{-1}=
 \smat{
 2 & 0 & 1 & 2 & 0 & 2 \\
 0 & 0 & 2 & 1 & 0 & 1 \\
 2 & 2 & 1 & 2 & 0 & 2 \\
 1 & 1 & 2 & 1 & 2 & 1 \\
 0 & 0 & 0 & 2 & 1 & 2 \\
 2 & 2 & 1 & 2 & 1 & 1 \\
}\ \
10) %\arraycolsep=1pt
 \smat{
 2&-1&0&0&0&0 \\
 -2&0&-1&0&0&0 \\
 0&-1&0&-2&0&0 \\
 0&0&-1&2&-1&-1 \\
 0&0&0&-1&2&0 \\
 0&0&0&-1&0&0
 }^{-1}=
 \smat{
 2 & 0 & 1 & 0 & 0 & 1 \\
 0 & 0 & 2 & 0 & 0 & 2 \\
 2 & 2 & 1 & 0 & 0 & 1 \\
 0 & 0 & 0 & 0 & 0 & 2 \\
 0 & 0 & 0 & 0 & 2 & 1 \\
 1 & 1 & 2 & 2 & 1 & 2 \\
}\\
11) \smat{
 0&-1&0&0&0&0 \\
 -1&2&-1&0&0&0 \\
 0&-1&0&-2&-2&0 \\
 0&0&-1&2&0&-1 \\
 0&0&-1&0&2&0 \\
 0&0&0&-1&0&2
 }^{-1}=
 \smat{
 1 & 2 & 0 & 1 & 0 & 2 \\
 2 & 0 & 0 & 0 & 0 & 0 \\
 0 & 0 & 0 & 2 & 0 & 1 \\
 2 & 0 & 1 & 2 & 1 & 1 \\
 0 & 0 & 0 & 1 & 2 & 2 \\
 1 & 0 & 2 & 1 & 2 & 1 \\
}\ \
12) %\arraycolsep=1pt
 \smat{
 0&-1&0&0&0&0 \\
 -1&0&-2&0&0&0 \\
 0&-1&2&0&-1&0 \\
 0&0&0&2&-1&-1 \\
 0&0&-1&-1&0&0 \\
 0&0&0&-1&0&2
 }^{-1}=
 \smat{
 1 & 2 & 2 & 1 & 0 & 2 \\
 2 & 0 & 0 & 0 & 0 & 0 \\
 1 & 0 & 2 & 1 & 0 & 2 \\
 2 & 0 & 1 & 2 & 2 & 1 \\
 0 & 0 & 0 & 2 & 0 & 1 \\
 1 & 0 & 2 & 1 & 1 & 1 \\
}\\
13) %\arraycolsep=1pt
 \smat{
 0&-1&0&0&0&0 \\
 -2&0&-1&0&0&0 \\
 0&-1&0&-2&-2&0 \\
 0&0&-1&2&0&-1 \\
 0&0&-1&0&2&0 \\
 0&0&0&-1&0&2
 }^{-1}=
 \smat{
 0 & 1 & 0 & 2 & 0 & 1 \\
 2 & 0 & 0 & 0 & 0 & 0 \\
 0 & 0 & 0 & 2 & 0 & 1 \\
 2 & 0 & 1 & 2 & 1 & 1 \\
 0 & 0 & 0 & 1 & 2 & 2 \\
 1 & 0 & 2 & 1 & 2 & 1 \\
}\ \
14) \smat{
 0&-1&0&0&0&0 \\
 -1&2&-1&0&0&0 \\
 0&-1&2&0&-1&0 \\
 0&0&0&2&-1&-1 \\
 0&0&-1&-1&0&0 \\
 0&0&0&-1&0&2
 }^{-1}=
 \smat{
 0 & 2 & 1 & 2 & 0 & 1 \\
 2 & 0 & 0 & 0 & 0 & 0 \\
 1 & 0 & 2 & 1 & 0 & 2 \\
 2 & 0 & 1 & 2 & 2 & 1 \\
 0 & 0 & 0 & 2 & 0 & 1 \\
 1 & 0 & 2 & 1 & 1 & 1 \\
}\\
15) \smat{
 2&-1&0&0&0&0 \\
 -1&2&-1&0&0&0 \\
 0&-2&0&-1&0&0 \\
 0&0&-1&0&-2&-2 \\
 0&0&0&-1&2&0 \\
 0&0&0&-1&0&0
 }^{-1}=
 \smat{
 2 & 0 & 2 & 0 & 0 & 1 \\
 0 & 0 & 1 & 0 & 0 & 2 \\
 1 & 2 & 0 & 0 & 0 & 0 \\
 0 & 0 & 0 & 0 & 0 & 2 \\
 0 & 0 & 0 & 0 & 2 & 1 \\
 1 & 2 & 0 & 1 & 1 & 2 \\
}\ \
16) %\arraycolsep=1pt
 \smat{
 2&-1&0&0&0&0 \\
 -2&0&-1&0&0&0 \\
 0&-1&0&0&-2&0 \\
 0&0&0&2&-1&-1 \\
 0&0&-1&-1&0&0 \\
 0&0&0&-1&0&2
 }^{-1}=
 \smat{
 2 & 0 & 1 & 1 & 0 & 2 \\
 0 & 0 & 2 & 2 & 0 & 1 \\
 2 & 2 & 1 & 1 & 0 & 2 \\
 1 & 1 & 2 & 2 & 2 & 1 \\
 0 & 0 & 0 & 2 & 0 & 1 \\
 2 & 2 & 1 & 1 & 1 & 1 \\
}\\
17) \smat{
 2&-1&0&0&0&0 \\
 -1&0&-2&0&0&0 \\
 0&-1&2&-1&-1&0 \\
 0&0&-1&2&0&-1 \\
 0&0&-1&0&2&0 \\
 0&0&0&-1&0&2
 }^{-1}=
 \smat{
 0 & 2 & 0 & 2 & 0 & 1 \\
 2 & 1 & 0 & 1 & 0 & 2 \\
 0 & 0 & 0 & 2 & 0 & 1 \\
 1 & 2 & 2 & 2 & 1 & 1 \\
 0 & 0 & 0 & 1 & 2 & 2 \\
 2 & 1 & 1 & 1 & 2 & 1 \\
}
\ \
18) \smat{
 2&-1&0&0&0&0 \\
 -1&2&-1&0&0&0 \\
 0&-1&2&-1&0&0 \\
 0&0&-2&0&-1&-1 \\
 0&0&0&-1&0&-1 \\
 0&0&0&-1&-1&2
 }^{-1}=
 \smat{
 2 & 0 & 1 & 1 & 0 & 2 \\
 0 & 0 & 2 & 2 & 0 & 1 \\
 1 & 2 & 0 & 0 & 0 & 0 \\
 2 & 1 & 0 & 1 & 0 & 2 \\
 0 & 0 & 0 & 0 & 1 & 2 \\
 1 & 2 & 0 & 2 & 2 & 1 \\
}\\
19) %\arraycolsep=1pt
 \smat{
 2&-1&0&0&0&0 \\
 -1&2&-1&0&0&0 \\
 0&-2&0&0&-1&0 \\
 0&0&0&2&-1&-1 \\
 0&0&-1&-2&2&0 \\
 0&0&0&-1&0&2
 }^{-1}=
 \smat{
 2 & 0 & 2 & 1 & 0 & 2 \\
 0 & 0 & 1 & 2 & 0 & 1 \\
 1 & 2 & 0 & 0 & 0 & 0 \\
 1 & 2 & 0 & 2 & 1 & 1 \\
 0 & 0 & 0 & 2 & 0 & 1 \\
 2 & 1 & 0 & 1 & 2 & 1 \\
}\ \
20) \smat{
 2&-1&0&0&0&0 \\
 -1&2&-1&0&0&0 \\
 0&-1&2&-1&0&0 \\
 0&0&-1&2&-1&0 \\
 0&0&0&-1&0&-1 \\
 0&0&0&0&-1&0
 }^{-1}=
 \smat{
 2 & 0 & 1 & 2 & 0 & 1 \\
 0 & 0 & 2 & 1 & 0 & 2 \\
 1 & 2 & 0 & 0 & 0 & 0 \\
 2 & 1 & 0 & 2 & 0 & 1 \\
 0 & 0 & 0 & 0 & 0 & 2 \\
 1 & 2 & 0 & 1 & 2 & 2 \\
}\\
21) \smat{
 2&-1&0&0&0&0 \\
 -1&2&-1&0&0&0 \\
 0&-1&2&-1&0&0 \\
 0&0&-1&2&-1&0 \\
 0&0&0&-2&2&-1 \\
 0&0&0&0&-1&0
 }^{-1}=
 \smat{
 2 & 0 & 1 & 2 & 0 & 1 \\
 0 & 0 & 2 & 1 & 0 & 2 \\
 1 & 2 & 0 & 0 & 0 & 0 \\
 2 & 1 & 0 & 2 & 0 & 1 \\
 0 & 0 & 0 & 0 & 0 & 2 \\
 2 & 1 & 0 & 2 & 2 & 2 \\
}\\
\end{matrix}
\]

\subsubsection{$\fg{}(8,6)$} 
\[\tiny \arraycolsep=1.5pt
\begin{matrix} 1)
\smat{
 2&0&-1&0&0&0&0 \\
 0&2&0&-1&0&0&0 \\
 -1&0&2&-1&0&0&0 \\
 0&-1&-1&2&-1&0&0 \\
 0&0&0&-1&2&-1&0 \\
 0&0&0&0&-2&0&-1 \\
 0&0&0&0&0&-1&0
 }^{-1}=
 \smat{
 1 & 2 & 1 & 1 & 2 & 0 & 1 \\
 2 & 2 & 1 & 0 & 0 & 0 & 0 \\
 1 & 1 & 2 & 2 & 1 & 0 & 2 \\
 1 & 0 & 2 & 0 & 0 & 0 & 0 \\
 2 & 0 & 1 & 0 & 2 & 0 & 1 \\
 0 & 0 & 0 & 0 & 0 & 0 & 2 \\
 2 & 0 & 1 & 0 & 2 & 2 & 1 \\
}\ \
2) \smat{
 2&0&-1&0&0&0&0 \\
 0&2&0&-1&0&0&0 \\
 -1&0&2&-1&0&0&0 \\
 0&-1&-1&2&-1&0&0 \\
 0&0&0&-2&0&-1&0 \\
 0&0&0&0&-1&0&-2 \\
 0&0&0&0&0&-1&2
 }^{-1}=
 \smat{
 1 & 2 & 1 & 1 & 1 & 2 & 2 \\
 2 & 2 & 1 & 0 & 0 & 0 & 0 \\
 1 & 1 & 2 & 2 & 2 & 1 & 1 \\
 1 & 0 & 2 & 0 & 0 & 0 & 0 \\
 2 & 0 & 1 & 0 & 1 & 2 & 2 \\
 1 & 0 & 2 & 0 & 2 & 0 & 0 \\
 2 & 0 & 1 & 0 & 1 & 0 & 2 \\
}\\
3) %\arraycolsep=5pt
 \smat{
 2&0&-1&0&0&0&0 \\
 0&2&0&-1&0&0&0 \\
 -1&0&2&-1&0&0&0 \\
 0&-1&-1&2&-1&0&0 \\
 0&0&0&-1&2&-1&0 \\
 0&0&0&0&-1&2&-1 \\
 0&0&0&0&0&-1&0
 }^{-1}=
 \smat{
 1 & 2 & 1 & 1 & 2 & 0 & 1 \\
 2 & 2 & 1 & 0 & 0 & 0 & 0 \\
 1 & 1 & 2 & 2 & 1 & 0 & 2 \\
 1 & 0 & 2 & 0 & 0 & 0 & 0 \\
 2 & 0 & 1 & 0 & 2 & 0 & 1 \\
 0 & 0 & 0 & 0 & 0 & 0 & 2 \\
 1 & 0 & 2 & 0 & 1 & 2 & 0 \\
}
\ \
4) \smat{
 2&0&-1&0&0&0&0 \\
 0&2&0&-1&0&0&0 \\
 -1&0&2&-1&0&0&0 \\
 0&-2&-2&0&-1&0&0 \\
 0&0&0&-1&0&-2&0 \\
 0&0&0&0&-1&2&-1 \\
 0&0&0&0&0&-1&2
 }^{-1}=
 \smat{
 1 & 2 & 1 & 2 & 0 & 1 & 2 \\
 2 & 2 & 1 & 0 & 0 & 0 & 0 \\
 1 & 1 & 2 & 1 & 0 & 2 & 1 \\
 1 & 0 & 2 & 0 & 0 & 0 & 0 \\
 0 & 0 & 0 & 0 & 0 & 2 & 1 \\
 1 & 0 & 2 & 0 & 1 & 0 & 0 \\
 2 & 0 & 1 & 0 & 2 & 0 & 2 \\
}\\
5) %\arraycolsep=5pt
 \smat{
 2&0&-1&0&0&0&0 \\
 0&0&-1&-1&0&0&0 \\
 -2&-1&0&-1&0&0&0 \\
 0&-1&-1&0&-2&0&0 \\
 0&0&0&-1&2&-1&0 \\
 0&0&0&0&-1&2&-1 \\
 0&0&0&0&0&-1&2
 }^{-1}=
 \smat{
 1 & 1 & 2 & 1 & 0 & 1 & 2 \\
 2 & 1 & 2 & 0 & 0 & 0 & 0 \\
 1 & 2 & 1 & 2 & 0 & 2 & 1 \\
 2 & 0 & 2 & 1 & 0 & 1 & 2 \\
 0 & 0 & 0 & 0 & 0 & 2 & 1 \\
 1 & 0 & 1 & 2 & 2 & 0 & 0 \\
 2 & 0 & 2 & 1 & 1 & 0 & 2 \\
}\ \
6) \smat{
 2&0&-1&0&0&0&0 \\
 0&0&-2&-2&0&0&0 \\
 -1&-1&2&0&0&0&0 \\
 0&-1&0&2&-1&0&0 \\
 0&0&0&-1&2&-1&0 \\
 0&0&0&0&-1&2&-1 \\
 0&0&0&0&0&-1&2
 }^{-1}=
 \smat{
 1 & 2 & 1 & 2 & 0 & 1 & 2 \\
 1 & 0 & 2 & 0 & 0 & 0 & 0 \\
 1 & 1 & 2 & 1 & 0 & 2 & 1 \\
 2 & 0 & 1 & 2 & 0 & 1 & 2 \\
 0 & 0 & 0 & 0 & 0 & 2 & 1 \\
 1 & 0 & 2 & 1 & 2 & 0 & 0 \\
 2 & 0 & 1 & 2 & 1 & 0 & 2 \\
}\\
7) %\arraycolsep=5pt
 \smat{
 0&0&-1&0&0&0&0 \\
 0&2&-1&0&0&0&0 \\
 -1&-2&0&-2&0&0&0 \\
 0&0&-1&2&-1&0&0 \\
 0&0&0&-1&2&-1&0 \\
 0&0&0&0&-1&2&-1 \\
 0&0&0&0&0&-1&2
 }^{-1}=
 \smat{
 2 & 2 & 2 & 2 & 0 & 1 & 2 \\
 1 & 2 & 0 & 0 & 0 & 0 & 0 \\
 2 & 0 & 0 & 0 & 0 & 0 & 0 \\
 1 & 0 & 0 & 2 & 0 & 1 & 2 \\
 0 & 0 & 0 & 0 & 0 & 2 & 1 \\
 2 & 0 & 0 & 1 & 2 & 0 & 0 \\
 1 & 0 & 0 & 2 & 1 & 0 & 2 \\
}\ \
8)\smat{
 0&0&-1&0&0&0&0 \\
 0&2&-1&0&0&0&0 \\
 -1&-1&2&-1&0&0&0 \\
 0&0&-1&2&-1&0&0 \\
 0&0&0&-1&2&-1&0 \\
 0&0&0&0&-1&2&-1 \\
 0&0&0&0&0&-1&2
 }^{-1}=
 \smat{
 2 & 1 & 2 & 1 & 0 & 2 & 1 \\
 1 & 2 & 0 & 0 & 0 & 0 & 0 \\
 2 & 0 & 0 & 0 & 0 & 0 & 0 \\
 1 & 0 & 0 & 2 & 0 & 1 & 2 \\
 0 & 0 & 0 & 0 & 0 & 2 & 1 \\
 2 & 0 & 0 & 1 & 2 & 0 & 0 \\
 1 & 0 & 0 & 2 & 1 & 0 & 2 \\
}\\
\end{matrix}
\]

\subsubsection{$\fel(5;3)$}
\[\tiny \arraycolsep=1.5pt
 \begin{matrix}
1) \;
\smat{
 0 & -1 & 0 & 0 & 0 \\
 -1 & 0 & 0 & -1 & 0 \\
 0 & 0 & 2 & -1 & -1 \\
 0 & -1 & -1 & 2 & 0 \\
 0 & 0 & -1 & 0 & 2\\
}^{-1}=
\smat{
 0 & 2 & 1 & 0 & 2 \\
 2 & 0 & 0 & 0 & 0 \\
 1 & 0 & 1 & 2 & 2 \\
 0 & 0 & 2 & 0 & 1 \\
 2 & 0 & 2 & 1 & 0 \\
}\ \
2) \;
\smat{
 0 &-2 & 0 & 0 & 0 \\
 -1 & 2 & 0 & -2 & 0 \\
 0 & 0 & 2 & -1 & -1 \\
 0 & -1 & -1 & 2 & 0 \\
 0 & 0 & -1 & 0 & 2\\
}^{-1}=
\smat{
 2 & 2 & 2 & 0 & 1 \\
 1 & 0 & 0 & 0 & 0 \\
 2 & 0 & 1 & 2 & 2 \\
 0 & 0 & 2 & 0 & 1 \\
 1 & 0 & 2 & 1 & 0 \\
 }
\\
 3) \;
\smat{
 2 & -1 & 0 & -1 & 0 \\
-2 & 0 & 0 &-2 & 0 \\
 0 & 0 & 2 & -1 & -1 \\
-2 &-2 & -1 & 0 & 0 \\
 0 & 0 & -1 & 0 & 2\\
}^{-1}=
\smat{
 0 & 1 & 1 & 0 & 2 \\
 2 & 2 & 0 & 0 & 0 \\
 2 & 0 & 1 & 2 & 2 \\
 0 & 0 & 2 & 0 & 1 \\
 1 & 0 & 2 & 1 & 0 \\
 }
\ \ 4) \;
\smat{
 0 & 0 & 0 & -1 & 0 \\
 0 & 2 & 0 & -1 & 0 \\
 0 & 0 & 0 &-2 & -1 \\
 -1 & -1 &-2 & 0 & 0 \\
 0 & 0 & -1 & 0 & 2\\
}^{-1}=
\smat{
 0 & 1 & 1 & 2 & 2 \\
 1 & 2 & 0 & 0 & 0 \\
 1 & 0 & 1 & 0 & 2 \\
 2 & 0 & 0 & 0 & 0 \\
 2 & 0 & 2 & 0 & 0 \\
 }
\\
 5) \;
\smat{
 0 & 0 & 0 &-2 & 0 \\
 0 & 2 & 0 & -1 & 0 \\
 0 & 0 & 0 &-2 & -1 \\
 -1 & -2 & -1 & 2 & 0 \\
 0 & 0 & -1 & 0 & 2\\
}^{-1}=
\smat{
 2 & 2 & 2 & 2 & 1 \\
 2 & 2 & 0 & 0 & 0 \\
 2 & 0 & 1 & 0 & 2 \\
 1 & 0 & 0 & 0 & 0 \\
 1 & 0 & 2 & 0 & 0 \\
 }
\ \ 6) \;
\smat{
 0 & 0 & 0 & -1 & 0 \\
 0 & 2 & 0 & -1 & 0 \\
 0 & 0 & 0 & -1 &-2 \\
 -1 & -1 & -1 & 2 & 0 \\
 0 & 0 &-2 & 0 & 0\\
}^{-1}=
\smat{
 0 & 1 & 0 & 2 & 2 \\
 1 & 2 & 0 & 0 & 0 \\
 0 & 0 & 0 & 0 & 1 \\
 2 & 0 & 0 & 0 & 0 \\
 2 & 0 & 1 & 0 & 0 \\
 }
\\
7) \;
\smat{
 0 & 0 & 0 &-2 & 0 \\
 0 & 2 & 0 & -1 & 0 \\
 0 & 0 & 0 & -1 &-2 \\
-2 & -1 & -1 & 0 & 0 \\
 0 & 0 &-2 & 0 & 0\\
}^{-1}=
\smat{
 2 & 2 & 0 & 1 & 1 \\
 2 & 2 & 0 & 0 & 0 \\
 0 & 0 & 0 & 0 & 1 \\
 1 & 0 & 0 & 0 & 0 \\
 1 & 0 & 1 & 0 & 0 \\
 }
\ \ 8) \;
\smat{
 0 & 0 & 0 & -1 & 0 \\
 0 & 2 & 0 & -1 & 0 \\
 0 & 0 & 2 & -1 & -1 \\
 -1 & -1 & -1 & 2 & 0 \\
 0 & 0 & -1 & 0 & 0\\
}^{-1}=
\smat{
 0 & 1 & 0 & 2 & 1 \\
 1 & 2 & 0 & 0 & 0 \\
 0 & 0 & 0 & 0 & 2 \\
 2 & 0 & 0 & 0 & 0 \\
 1 & 0 & 2 & 0 & 1 \\
 }
\\ 9) \;
\smat{
 2 & 0 & 0 & -1 & 0 \\
 0 & 0 &-2 &-2 & 0 \\
 0 & -1 & 2 & -1 & -1 \\
 -1 &-2 &-2 & 0 & 0 \\
 0 & 0 &-2 & 0 & 0\\
}^{-1}=
\smat{
 2 & 2 & 0 & 0 & 1 \\
 2 & 2 & 0 & 1 & 0 \\
 0 & 0 & 0 & 0 & 1 \\
 0 & 1 & 0 & 0 & 2 \\
 1 & 0 & 2 & 2 & 0 \\
 }
\ \ 10)\;
\smat{
 0 & 0 & 0 &-2 & 0 \\
 0 & 2 & 0 & -1 & 0 \\
 0 & 0 & 2 & -1 & -1 \\
-2 & -1 & -1 & 0 & 0 \\
 0 & 0 & -1 & 0 & 0\\
}^{-1}=
\smat{
 2 & 2 & 0 & 1 & 2 \\
 2 & 2 & 0 & 0 & 0 \\
 0 & 0 & 0 & 0 & 2 \\
 1 & 0 & 0 & 0 & 0 \\
 2 & 0 & 2 & 0 & 1 \\
 }
\\ 11)\;
\smat{
 2 & 0 & 0 & -1 & 0 \\
 0 & 0 & -1 & -1 & 0 \\
 0 & -1 & 0 & 0 &-2 \\
 -1 & -1 & 0 & 2 & 0 \\
 0 & 0 &-2 & 0 & 0\\
}^{-1}=
\smat{
 2 & 1 & 0 & 0 & 1 \\
 1 & 0 & 0 & 2 & 0 \\
 0 & 0 & 0 & 0 & 1 \\
 0 & 2 & 0 & 0 & 2 \\
 1 & 0 & 1 & 2 & 0 \\
 }
\ \ 12)\;
\smat{
 2 & 0 & 0 & -1 & 0 \\
 0 & 0 &-2 &-2 & 0 \\
 0 &-2 & 0 &-2 & -1 \\
 -1 &-2 &-2 & 0 & 0 \\
 0 & 0 & -1 & 0 & 0\\
}^{-1}=
\smat{
 2 & 2 & 0 & 0 & 2 \\
 2 & 2 & 0 & 1 & 0 \\
 0 & 0 & 0 & 0 & 2 \\
 0 & 1 & 0 & 0 & 1 \\
 2 & 0 & 2 & 1 & 1 \\
 }
\\ 13)\;
\smat{
 2 & 0 & 0 & -1 & 0 \\
 0 & 2 & -1 & -2 & 0 \\
 0 &-2 & 0 & 0 & -1 \\
 -1 & -1 & 0 & 2 & 0 \\
 0 & 0 & -1 & 0 & 2\\
}^{-1}=
\smat{
 2 & 2 & 2 & 0 & 1 \\
 1 & 0 & 0 & 2 & 0 \\
 2 & 0 & 1 & 1 & 2 \\
 0 & 1 & 1 & 0 & 2 \\
 1 & 0 & 2 & 2 & 0 \\
 }
\ \ 14)\;
\smat{
 2 & 0 & 0 & -1 & 0 \\
 0 & 0 & -1 & -1 & 0 \\
 0 & -1 & 2 & 0 & -1 \\
 -1 & -1 & 0 & 2 & 0 \\
 0 & 0 & -1 & 0 & 0\\
}^{-1}=
\smat{
 2 & 1 & 0 & 0 & 2 \\
 1 & 0 & 0 & 2 & 0 \\
 0 & 0 & 0 & 0 & 2 \\
 0 & 2 & 0 & 0 & 1 \\
 2 & 0 & 2 & 1 & 1 \\
 }
\\ 15)\;
\smat{
 2 & 0 & 0 & -1 & 0 \\
 0 & 2 & -1 & 0 & 0 \\
 0 & -1 & 0 & -1 &-2 \\
 -1 & 0 & -1 & 2 & 0 \\
 0 & 0 & -1 & 0 & 2\\
}^{-1}=
\smat{
 2 & 2 & 1 & 0 & 1 \\
 2 & 2 & 0 & 1 & 0 \\
 1 & 0 & 0 & 2 & 0 \\
 0 & 1 & 2 & 0 & 2 \\
 2 & 0 & 0 & 1 & 2 \\
 }
\end{matrix}
\]

\ssec{$p=5$}{}~{} 

\subsubsection{$\fbrj(2;5)$} $ 1)\; \smat{0&-1\\
-2&1\\}^{-1}=\frac13\smat{1&1\\
2&0\\}, \quad 2)\; \smat{0&-1\\
-3&2\\}^{-1}=\frac12\smat{2&1\\
3&0\\
}.
$

\subsubsection{$\fel(5;5)$}
\[\tiny\arraycolsep=1.5pt
\begin{matrix}
1) \smat{
2&0&-1&0&0 \\
 0&2&0&0&-1 \\
 -1&0&0&-4&-4 \\
 0&0&-4&0&-2 \\
 0&-1&-4&-2&0\\
}^{-1}=
\smat{
 2 & 2 & 3 & 3 & 4 \\
 2 & 4 & 4 & 0 & 2 \\
 3 & 4 & 1 & 1 & 3 \\
 3 & 0 & 1 & 3 & 0 \\
 4 & 2 & 3 & 0 & 4 \\
}\ \
2) \smat{
 0&0&-4&0&0 \\
 0&2&0&0&-1 \\
 -4&0&0&-1&-1 \\
 0&0&-1&2&0 \\
 0&-1&-1&0&2\\
}^{-1}=
\smat{
 2 & 2 & 1 & 3 & 4 \\
 2 & 4 & 0 & 0 & 2 \\
 1 & 0 & 0 & 0 & 0 \\
 3 & 0 & 0 & 3 & 0 \\
 4 & 2 & 0 & 0 & 4 \\
}\\
3) \smat{
 2&0&-1&0&0 \\
 0&2&0&0&-1 \\
 -1&0&2&-1&0 \\
 0&0&-1&0&2 \\
 0&-2&0&-1&2\\
}^{-1}=
\smat{
 2 & 2 & 3 & 4 & 2 \\
 2 & 4 & 4 & 1 & 1 \\
 3 & 4 & 1 & 3 & 4 \\
 4 & 1 & 3 & 2 & 1 \\
 4 & 2 & 3 & 2 & 2 \\
}\ \ 
4) \smat{
 2&0&-1&0&0 \\
 0&0&0&2&-4 \\
 -1&0&2&0&-1 \\
 0&-1&0&2&-1 \\
 0&-4&-1&2&0\\
}^{-1}=
 \smat{
 2 & 2 & 3 & 4 & 4 \\
 2 & 4 & 4 & 0 & 1 \\
 3 & 4 & 1 & 3 & 3 \\
 3 & 0 & 1 & 4 & 4 \\
 4 & 1 & 3 & 2 & 2 \\
}\\

5) \smat{
 0&0&-1&0&0 \\
 0&2&0&0&-1 \\
 -1&0&2&-1&-1 \\
 0&0&-1&2&0 \\
 0&-1&-1&0&2\\
}^{-1}=
\smat{
 0 & 3 & 4 & 2 & 1 \\
 3 & 4 & 0 & 0 & 2 \\
 4 & 0 & 0 & 0 & 0 \\
 2 & 0 & 0 & 3 & 0 \\
 1 & 2 & 0 & 0 & 4 \\
}\ \
6) \smat{
 2&0&-1&0&0 \\
 0&0&0&-2&-1 \\
 -1&0&2&0&-1 \\
 0&-2&0&0&0 \\
 0&-1&-1&0&2\\
}^{-1}=
\smat{
 2 & 0 & 3 & 3 & 4 \\
 0 & 0 & 0 & 2 & 0 \\
 3 & 0 & 1 & 1 & 3 \\
 3 & 2 & 1 & 3 & 4 \\
 4 & 0 & 3 & 4 & 2 \\
}
\\
7)\smat{
 2&0&-1&0&0 \\
 0&2&0&-1&-2 \\
 -1&0&2&0&-1 \\
 0&2&0&0&0 \\
 0&-1&-1&0&2\\
}^{-1}=
\smat{
 2 & 0 & 3 & 2 & 4 \\
 0 & 0 & 0 & 3 & 0 \\
 3 & 0 & 1 & 4 & 3 \\
 2 & 4 & 4 & 4 & 1 \\
 4 & 0 & 3 & 1 & 2 \\
}\\ 
\end{matrix}
\]

\section[Almost affine Lie (super)algebras]{Simple hyperbolic Lie algebras and almost affine Lie superalgebras over $\Cee$}\label{AA}

The numbering of Cartan matrices follows that in the arXiv version of \cite{CCLL}, which contains Cartan matrices of all hyperbolic Lie algebras (classified by Li Wang Lai; for history, and rediscovery of Li Wang Lai's result in relation with cosmological billiards, see \cite{CCLL}). The published version of \cite{CCLL} has only new results: the Cartan matrices of almost affine Lie \textbf{super}algebras, whereas the arXiv version reproduces Li Wang Lai's result. 

Below, the symbol $\fbox{!!!}$ marks the cases where the inverse of the Cartan matrix has no zero elements, cf. \cite{ZH}. By Theorem~\ref{ThD}, these are precisely the cases where the deletion of any node of the Dynkin diagram leaves a finite-dimensional Lie (super)algebra, not an affine one.

\subsection{Cartan matrix $A$, $A^{-1}$, and $\det A$ for almost affine Lie superal\-gebras}\label{ss8.1}{}~{}

\sssec{On the names of the algebras}\label{S35} In the names borrowed from \cite{CCLL},
the initial $N$ says that the Cartan matrix is \textbf{not} symmetrizable. We have
verified this for each of the $370$ matrices of Subsections~\ref{ss8.1} and \ref{ssHyp}:
the name agrees with the computation in $369$ cases, and disagrees in one --- the matrix
called $S3_{35}$ in \cite{CCLL} is not symmetrizable, so we call it $NS3_{35}$ below.
%\AO{DL: the same misprint of \cite{CCLL} was spotted independently in our text on the hidden supersymmetry of cosmological billiards, and the two texts now agree.}

\noindent\stab{\attn$NS3_2$&
$\smat{
 2 & -1 & -1 \\
 -1 & 1 & -1 \\
 -3 & -1 & 2 \\
}$
 & $\frac{1}{7}$ 
$\smat{
 -1 & -3 & -2 \\
 -5 & -1 & -3 \\
 -4 & -5 & -1 \\
}$
 & -7 } \\
 \stab{$NS3_3$&
$\smat{
 2 & -1 & -1 \\
 -1 & 1 & -1 \\
 -4 & -1 & 2 \\
}$
 & $\frac{1}{9}$ 
$\smat{
 -1 & -3 & -2 \\
 -6 & 0 & -3 \\
 -5 & -6 & -1 \\
}$
 & -9 } \\
 \stab{$S3_4$&
$\smat{
 2 & -1 & -2 \\
 -1 & 1 & -1 \\
 -2 & -1 & 2 \\
}$
 & $\frac{1}{8}$ 
$\smat{
 -1 & -4 & -3 \\
 -4 & 0 & -4 \\
 -3 & -4 & -1 \\
}$
 & -8 } \\
 \stab{\attn$S3_5$&
$\smat{
 2 & -1 & -1 \\
 -1 & 1 & -1 \\
 -1 & -1 & 2 \\
}$
 & $\frac{1}{3}$ 
$\smat{
 -1 & -3 & -2 \\
 -3 & -3 & -3 \\
 -2 & -3 & -1 \\
}$
 & -3 } \\
 \stab{$S3_6$&
$\smat{
 2 & -1 & -1 \\
 -2 & 1 & -1 \\
 -2 & -1 & 2 \\
}$
 & $\frac{1}{8}$ 
$\smat{
 -1 & -3 & -2 \\
 -6 & -2 & -4 \\
 -4 & -4 & 0 \\
}$
 & -8 } \\
 \stab{$NS3_7$&
$\smat{
 2 & -1 & -1 \\
 -2 & 1 & -1 \\
 -3 & -1 & 2 \\
}$
 & $\frac{1}{10}$ 
$\smat{
 -1 & -3 & -2 \\
 -7 & -1 & -4 \\
 -5 & -5 & 0 \\
}$
 & -10 } \\
 \stab{$NS3_8$&
$\smat{
 2 & -1 & -1 \\
 -2 & 1 & -1 \\
 -4 & -1 & 2 \\
}$
 & $\frac{1}{12}$ 
$\smat{
 -1 & -3 & -2 \\
 -8 & 0 & -4 \\
 -6 & -6 & 0 \\
}$
 & -12 } \\
 \stab{$NS3_9$&
$\smat{
 2 & -1 & -2 \\
 -2 & 1 & -1 \\
 -2 & -1 & 2 \\
}$
 & $\frac{1}{12}$ 
$\smat{
 -1 & -4 & -3 \\
 -6 & 0 & -6 \\
 -4 & -4 & 0 \\
}$
 & -12 } \\
 \stab{$NS3_{10}$&
$\smat{
 2 & -1 & -1 \\
 -2 & 1 & -1 \\
 -1 & -1 & 2 \\
}$
 & $\frac{1}{6}$ 
$\smat{
 -1 & -3 & -2 \\
 -5 & -3 & -4 \\
 -3 & -3 & 0 \\
}$
 & -6 } \\
 \stab{$NS3_{11}$&
$\smat{
 2 & -1 & -2 \\
 -2 & 1 & -1 \\
 -1 & -1 & 2 \\
}$
 & $\frac{1}{9}$ 
$\smat{
 -1 & -4 & -3 \\
 -5 & -2 & -6 \\
 -3 & -3 & 0 \\
}$
 & -9 } \\
 \stab{$NS3_{12}$&
$\smat{
 2 & -1 & -3 \\
 -2 & 1 & -1 \\
 -1 & -1 & 2 \\
}$
 & $\frac{1}{12}$ 
$\smat{
 -1 & -5 & -4 \\
 -5 & -1 & -8 \\
 -3 & -3 & 0 \\
}$
 & -12 } \\
 \stab{$NS3_{13}$&
$\smat{
 2 & -1 & -4 \\
 -2 & 1 & -1 \\
 -1 & -1 & 2 \\
}$
 & $\frac{1}{15}$ 
$\smat{
 -1 & -6 & -5 \\
 -5 & 0 & -10 \\
 -3 & -3 & 0 \\
}$
 & -15 } \\
 \stab{$NS3_{14}$&
$\smat{
 2 & -2 & -1 \\
 -1 & 1 & -1 \\
 -2 & -1 & 2 \\
}$
 & $\frac{1}{9}$ 
$\smat{
 -1 & -5 & -3 \\
 -4 & -2 & -3 \\
 -3 & -6 & 0 \\
}$
 & -9 } \\
 \stab{$NS3_{15}$&
$\smat{
 2 & -2 & -1 \\
 -1 & 1 & -1 \\
 -3 & -1 & 2 \\
}$
 & $\frac{1}{12}$ 
$\smat{
 -1 & -5 & -3 \\
 -5 & -1 & -3 \\
 -4 & -8 & 0 \\
}$
 & -12 } \\
 \stab{$NS3_{16}$&
$\smat{
 2 & -2 & -1 \\
 -1 & 1 & -1 \\
 -4 & -1 & 2 \\
}$
 & $\frac{1}{15}$ 
$\smat{
 -1 & -5 & -3 \\
 -6 & 0 & -3 \\
 -5 & -10 & 0 \\
}$
 & -15 } \\
 \stab{$NS3_{17}$&
$\smat{
 2 & -2 & -2 \\
 -1 & 1 & -1 \\
 -2 & -1 & 2 \\
}$
 & $\frac{1}{12}$ 
$\smat{
 -1 & -6 & -4 \\
 -4 & 0 & -4 \\
 -3 & -6 & 0 \\
}$
 & -12 } \\
 \stab{$NS3_{18}$&
$\smat{
 2 & -2 & -1 \\
 -1 & 1 & -1 \\
 -1 & -1 & 2 \\
}$
 & $\frac{1}{6}$ 
$\smat{
 -1 & -5 & -3 \\
 -3 & -3 & -3 \\
 -2 & -4 & 0 \\
}$
 & -6 } \\
 \stab{$S3_{19}$&
$\smat{
 2 & -2 & -2 \\
 -1 & 1 & -1 \\
 -1 & -1 & 2 \\
}$
 & $\frac{1}{8}$ 
$\smat{
 -1 & -6 & -4 \\
 -3 & -2 & -4 \\
 -2 & -4 & 0 \\
}$
 & -8 } \\
 \stab{$NS3_{20}$&
$\smat{
 2 & -2 & -3 \\
 -1 & 1 & -1 \\
 -1 & -1 & 2 \\
}$
 & $\frac{1}{10}$ 
$\smat{
 -1 & -7 & -5 \\
 -3 & -1 & -5 \\
 -2 & -4 & 0 \\
}$
 & -10 } \\
 \stab{$NS3_{21}$&
$\smat{
 2 & -2 & -4 \\
 -1 & 1 & -1 \\
 -1 & -1 & 2 \\
}$
 & $\frac{1}{12}$ 
$\smat{
 -1 & -8 & -6 \\
 -3 & 0 & -6 \\
 -2 & -4 & 0 \\
}$
 & -12 } \\
 \stab{$NS3_{22}$&
$\smat{
 2 & -1 & -1 \\
 -2 & 1 & -2 \\
 -2 & -1 & 2 \\
}$
 & $\frac{1}{12}$ 
$\smat{
 0 & -3 & -3 \\
 -8 & -2 & -6 \\
 -4 & -4 & 0 \\
}$
 & -12 } \\
 \stab{$NS3_{23}$&
$\smat{
 2 & -1 & -1 \\
 -2 & 1 & -2 \\
 -3 & -1 & 2 \\
}$
 & $\frac{1}{15}$ 
$\smat{
 0 & -3 & -3 \\
 -10 & -1 & -6 \\
 -5 & -5 & 0 \\
}$
 & -15 } \\
 \stab{$NS3_{24}$&
$\smat{
 2 & -1 & -1 \\
 -2 & 1 & -2 \\
 -4 & -1 & 2 \\
}$
 & $\frac{1}{6}$ 
$\smat{
 0 & -1 & -1 \\
 -4 & 0 & -2 \\
 -2 & -2 & 0 \\
}$
 & -18 } \\
 \stab{$S3_{25}$&
$\smat{
 2 & -1 & -2 \\
 -2 & 1 & -2 \\
 -2 & -1 & 2 \\
}$
 & $\frac{1}{4}$ 
$\smat{
 0 & -1 & -1 \\
 -2 & 0 & -2 \\
 -1 & -1 & 0 \\
}$
 & -16 } \\
 \stab{$S3_{26}$&
$\smat{
 2 & -1 & -1 \\
 -2 & 1 & -2 \\
 -1 & -1 & 2 \\
}$
 & $\frac{1}{3}$ 
$\smat{
 0 & -1 & -1 \\
 -2 & -1 & -2 \\
 -1 & -1 & 0 \\
}$
 & -9 } \\
 \stab{$NS3_{27}$&
$\smat{
 2 & -2 & -1 \\
 -1 & 1 & -2 \\
 -2 & -1 & 2 \\
}$
 & $\frac{1}{15}$ 
$\smat{
 0 & -5 & -5 \\
 -6 & -2 & -5 \\
 -3 & -6 & 0 \\
}$
 & -15 } \\
 \stab{$NS3_{28}$&
$\smat{
 2 & -2 & -1 \\
 -1 & 1 & -2 \\
 -3 & -1 & 2 \\
}$
 & $\frac{1}{20}$ 
$\smat{
 0 & -5 & -5 \\
 -8 & -1 & -5 \\
 -4 & -8 & 0 \\
}$
 & -20 } \\
 \stab{$NS3_{29}$&
$\smat{
 2 & -2 & -1 \\
 -1 & 1 & -2 \\
 -4 & -1 & 2 \\
}$
 & $\frac{1}{5}$ 
$\smat{
 0 & -1 & -1 \\
 -2 & 0 & -1 \\
 -1 & -2 & 0 \\
}$
 & -25 } \\
 \stab{$NS3_{30}$&
$\smat{
 2 & -2 & -2 \\
 -1 & 1 & -2 \\
 -2 & -1 & 2 \\
}$
 & $\frac{1}{6}$ 
$\smat{
 0 & -2 & -2 \\
 -2 & 0 & -2 \\
 -1 & -2 & 0 \\
}$
 & -18 } \\
 \stab{$NS3_{31}$&
$\smat{
 2 & -2 & -1 \\
 -1 & 1 & -2 \\
 -1 & -1 & 2 \\
}$
 & $\frac{1}{10}$ 
$\smat{
 0 & -5 & -5 \\
 -4 & -3 & -5 \\
 -2 & -4 & 0 \\
}$
 & -10 } \\
 \stab{$NS3_{32}$&
$\smat{
 2 & -2 & -2 \\
 -1 & 1 & -2 \\
 -1 & -1 & 2 \\
}$
 & $\frac{1}{6}$ 
$\smat{
 0 & -3 & -3 \\
 -2 & -1 & -3 \\
 -1 & -2 & 0 \\
}$
 & -12 } \\
 \stab{$NS3_{33}$&
$\smat{
 2 & -2 & -3 \\
 -1 & 1 & -2 \\
 -1 & -1 & 2 \\
}$
 & $\frac{1}{14}$ 
$\smat{
 0 & -7 & -7 \\
 -4 & -1 & -7 \\
 -2 & -4 & 0 \\
}$
 & -14 } \\
 \stab{$S3_{34}$&
$\smat{
 2 & -2 & -4 \\
 -1 & 1 & -2 \\
 -1 & -1 & 2 \\
}$
 & $\frac{1}{8}$ 
$\smat{
 0 & -4 & -4 \\
 -2 & 0 & -4 \\
 -1 & -2 & 0 \\
}$
 & -16 } \\
 \stab{$NS3_{35}$&
$\smat{
 2 & -2 & -1 \\
 -1 & 1 & -1 \\
 -2 & -2 & 2 \\
}$
 & $\frac{1}{12}$ 
$\smat{
 0 & -6 & -3 \\
 -4 & -2 & -3 \\
 -4 & -8 & 0 \\
}$
 & -12 } \\
 \stab{$NS3_{36}$&
$\smat{
 2 & -2 & -1 \\
 -1 & 1 & -1 \\
 -3 & -2 & 2 \\
}$
 & $\frac{1}{15}$ 
$\smat{
 0 & -6 & -3 \\
 -5 & -1 & -3 \\
 -5 & -10 & 0 \\
}$
 & -15 } \\
 \stab{$NS3_{37}$&
$\smat{
 2 & -2 & -1 \\
 -1 & 1 & -1 \\
 -4 & -2 & 2 \\
}$
 & $\frac{1}{6}$ 
$\smat{
 0 & -2 & -1 \\
 -2 & 0 & -1 \\
 -2 & -4 & 0 \\
}$
 & -18 } \\
 \stab{$S3_{38}$&
$\smat{
 2 & -2 & -2 \\
 -1 & 1 & -1 \\
 -2 & -2 & 2 \\
}$
 & $\frac{1}{4}$ 
$\smat{
 0 & -2 & -1 \\
 -1 & 0 & -1 \\
 -1 & -2 & 0 \\
}$
 & -16 } \\
 \stab{$S3_{39}$&
$\smat{
 2 & -2 & -1 \\
 -1 & 1 & -1 \\
 -1 & -2 & 2 \\
}$
 & $\frac{1}{3}$ 
$\smat{
 0 & -2 & -1 \\
 -1 & -1 & -1 \\
 -1 & -2 & 0 \\
}$
 & -9 } \\
 \stab{$S3_{40}$&
$\smat{
 2 & -1 & -1 \\
 -1 & 1 & -1 \\
 -1 & -1 & 1 \\
}$
 & $\frac{1}{4}$ 
$\smat{
 0 & -2 & -2 \\
 -2 & -1 & -3 \\
 -2 & -3 & -1 \\
}$
 & -4 } \\
 \stab{$NS3_{41}$&
$\smat{
 2 & -1 & -1 \\
 -2 & 1 & -1 \\
 -1 & -1 & 1 \\
}$
 & $\frac{1}{6}$ 
$\smat{
 0 & -2 & -2 \\
 -3 & -1 & -4 \\
 -3 & -3 & 0 \\
}$
 & -6 } \\
 \stab{$NS3_{42}$&
$\smat{
 2 & -2 & -1 \\
 -1 & 1 & -1 \\
 -1 & -1 & 1 \\
}$
 & $\frac{1}{6}$ 
$\smat{
 0 & -3 & -3 \\
 -2 & -1 & -3 \\
 -2 & -4 & 0 \\
}$
 & -6 } \\
 \stab{$S3_{43}$&
$\smat{
 2 & -1 & -1 \\
 -2 & 1 & -1 \\
 -2 & -1 & 1 \\
}$
 & $\frac{1}{4}$ 
$\smat{
 0 & -1 & -1 \\
 -2 & 0 & -2 \\
 -2 & -2 & 0 \\
}$
 & -8 } \\
 \stab{$NS3_{44}$&
$\smat{
 2 & -2 & -1 \\
 -1 & 1 & -1 \\
 -2 & -1 & 1 \\
}$
 & $\frac{1}{3}$ 
$\smat{
 0 & -1 & -1 \\
 -1 & 0 & -1 \\
 -1 & -2 & 0 \\
}$
 & -9 } \\
 \stab{$S3_{45}$&
$\smat{
 2 & -2 & -2 \\
 -1 & 1 & -1 \\
 -1 & -1 & 1 \\
}$
 & $\frac{1}{4}$ 
$\smat{
 0 & -2 & -2 \\
 -1 & 0 & -2 \\
 -1 & -2 & 0 \\
}$
 & -8 } \\
 \stab{$S3_{46}$&
$\smat{
 1 & -1 & -1 \\
 -1 & 1 & -1 \\
 -1 & -1 & 1 \\
}$
 & $\frac{1}{2}$ 
$\smat{
 0 & -1 & -1 \\
 -1 & 0 & -1 \\
 -1 & -1 & 0 \\
}$
 & -4 } \\
 \stab{$S3_{47}$&
$\smat{
 2 & -1 & 0 \\
 -2 & 1 & -1 \\
 0 & -1 & 2 \\
}$
 & $\frac{1}{2}$ 
$\smat{
 -1 & -2 & -1 \\
 -4 & -4 & -2 \\
 -2 & -2 & 0 \\
}$
 & -2 } \\
 \stab{$S3_{48}$&
$\smat{
 2 & -2 & 0 \\
 -1 & 1 & -1 \\
 0 & -1 & 2 \\
}$
 & $\frac{1}{2}$ 
$\smat{
 -1 & -4 & -2 \\
 -2 & -4 & -2 \\
 -1 & -2 & 0 \\
}$
 & -2 } \\
 \stab{$S3_{49}$&
$\smat{
 2 & -1 & 0 \\
 -2 & 1 & -2 \\
 0 & -1 & 2 \\
}$
 & $\frac{1}{2}$ 
$\smat{
 0 & -1 & -1 \\
 -2 & -2 & -2 \\
 -1 & -1 & 0 \\
}$
 & -4 } \\
 \stab{$S3_{50}$&
$\smat{
 2 & -2 & 0 \\
 -1 & 1 & -2 \\
 0 & -1 & 2 \\
}$
 & $\frac{1}{4}$ 
$\smat{
 0 & -4 & -4 \\
 -2 & -4 & -4 \\
 -1 & -2 & 0 \\
}$
 & -4 } \\
 \stab{$S3_{51}$&
$\smat{
 2 & -2 & 0 \\
 -1 & 1 & -1 \\
 0 & -2 & 2 \\
}$
 & $\frac{1}{2}$ 
$\smat{
 0 & -2 & -1 \\
 -1 & -2 & -1 \\
 -1 & -2 & 0 \\
}$
 & -4 } \\
 \stab{\attn$S3_{52}$&
$\smat{
 2 & -1 & 0 \\
 -3 & 2 & -1 \\
 0 & -1 & 1 \\
}$
 & 
$\smat{
 -1 & -1 & -1 \\
 -3 & -2 & -2 \\
 -3 & -2 & -1 \\
}$
 & -1 } \\
 \stab{$S3_{53}$&
$\smat{
 2 & -1 & 0 \\
 -4 & 2 & -1 \\
 0 & -1 & 1 \\
}$
 & $\frac{1}{2}$ 
$\smat{
 -1 & -1 & -1 \\
 -4 & -2 & -2 \\
 -4 & -2 & 0 \\
}$
 & -2 } \\
 \stab{$S3_{54}$&
$\smat{
 2 & -2 & 0 \\
 -2 & 2 & -1 \\
 0 & -1 & 1 \\
}$
 & $\frac{1}{2}$ 
$\smat{
 -1 & -2 & -2 \\
 -2 & -2 & -2 \\
 -2 & -2 & 0 \\
}$
 & -2 } \\
 \stab{\attn$S3_{55}$&
$\smat{
 2 & -3 & 0 \\
 -1 & 2 & -1 \\
 0 & -1 & 1 \\
}$
 & 
$\smat{
 -1 & -3 & -3 \\
 -1 & -2 & -2 \\
 -1 & -2 & -1 \\
}$
 & -1 } \\
 \stab{$S3_{56}$&
$\smat{
 2 & -4 & 0 \\
 -1 & 2 & -1 \\
 0 & -1 & 1 \\
}$
 & $\frac{1}{2}$ 
$\smat{
 -1 & -4 & -4 \\
 -1 & -2 & -2 \\
 -1 & -2 & 0 \\
}$
 & -2 } \\
 \stab{$S3_{57}$&
$\smat{
 2 & -1 & 0 \\
 -2 & 2 & -1 \\
 0 & -2 & 1 \\
}$
 & $\frac{1}{2}$ 
$\smat{
 0 & -1 & -1 \\
 -2 & -2 & -2 \\
 -4 & -4 & -2 \\
}$
 & -2 } \\
 \stab{$S3_{58}$&
$\smat{
 2 & -1 & 0 \\
 -3 & 2 & -1 \\
 0 & -2 & 1 \\
}$
 & $\frac{1}{3}$ 
$\smat{
 0 & -1 & -1 \\
 -3 & -2 & -2 \\
 -6 & -4 & -1 \\
}$
 & -3 } \\
 \stab{$S3_{59}$&
$\smat{
 2 & -1 & 0 \\
 -4 & 2 & -1 \\
 0 & -2 & 1 \\
}$
 & $\frac{1}{4}$ 
$\smat{
 0 & -1 & -1 \\
 -4 & -2 & -2 \\
 -8 & -4 & 0 \\
}$
 & -4 } \\
 \stab{$S3_{60}$&
$\smat{
 2 & -2 & 0 \\
 -2 & 2 & -1 \\
 0 & -2 & 1 \\
}$
 & $\frac{1}{2}$ 
$\smat{
 0 & -1 & -1 \\
 -1 & -1 & -1 \\
 -2 & -2 & 0 \\
}$
 & -4 } \\
 \stab{$S3_{61}$&
$\smat{
 2 & -1 & 0 \\
 -1 & 2 & -1 \\
 0 & -2 & 1 \\
}$
 & 
$\smat{
 0 & -1 & -1 \\
 -1 & -2 & -2 \\
 -2 & -4 & -3 \\
}$
 & -1 } \\
 \stab{$S3_{62}$&
$\smat{
 2 & -2 & 0 \\
 -1 & 2 & -1 \\
 0 & -2 & 1 \\
}$
 & $\frac{1}{2}$ 
$\smat{
 0 & -2 & -2 \\
 -1 & -2 & -2 \\
 -2 & -4 & -2 \\
}$
 & -2 } \\
 \stab{$S3_{63}$&
$\smat{
 2 & -3 & 0 \\
 -1 & 2 & -1 \\
 0 & -2 & 1 \\
}$
 & $\frac{1}{3}$ 
$\smat{
 0 & -3 & -3 \\
 -1 & -2 & -2 \\
 -2 & -4 & -1 \\
}$
 & -3 } \\
 \stab{$S3_{64}$&
$\smat{
 2 & -4 & 0 \\
 -1 & 2 & -1 \\
 0 & -2 & 1 \\
}$
 & $\frac{1}{4}$ 
$\smat{
 0 & -4 & -4 \\
 -1 & -2 & -2 \\
 -2 & -4 & 0 \\
}$
 & -4 } \\
 \stab{$S3_{65}$&
$\smat{
 2 & -1 & 0 \\
 -2 & 2 & -2 \\
 0 & -1 & 1 \\
}$
 & $\frac{1}{2}$ 
$\smat{
 0 & -1 & -2 \\
 -2 & -2 & -4 \\
 -2 & -2 & -2 \\
}$
 & -2 } \\
 \stab{$S3_{66}$&
$\smat{
 2 & -1 & 0 \\
 -3 & 2 & -2 \\
 0 & -1 & 1 \\
}$
 & $\frac{1}{3}$ 
$\smat{
 0 & -1 & -2 \\
 -3 & -2 & -4 \\
 -3 & -2 & -1 \\
}$
 & -3 } \\
 \stab{$S3_{67}$&
$\smat{
 2 & -1 & 0 \\
 -4 & 2 & -2 \\
 0 & -1 & 1 \\
}$
 & $\frac{1}{4}$ 
$\smat{
 0 & -1 & -2 \\
 -4 & -2 & -4 \\
 -4 & -2 & 0 \\
}$
 & -4 } \\
 \stab{$S3_{68}$&
$\smat{
 2 & -2 & 0 \\
 -2 & 2 & -2 \\
 0 & -1 & 1 \\
}$
 & $\frac{1}{2}$ 
$\smat{
 0 & -1 & -2 \\
 -1 & -1 & -2 \\
 -1 & -1 & 0 \\
}$
 & -4 } \\
 \stab{$S3_{69}$&
$\smat{
 2 & -1 & 0 \\
 -1 & 2 & -2 \\
 0 & -1 & 1 \\
}$
 & 
$\smat{
 0 & -1 & -2 \\
 -1 & -2 & -4 \\
 -1 & -2 & -3 \\
}$
 & -1 } \\
 \stab{$S3_{70}$&
$\smat{
 2 & -2 & 0 \\
 -1 & 2 & -2 \\
 0 & -1 & 1 \\
}$
 & $\frac{1}{2}$ 
$\smat{
 0 & -2 & -4 \\
 -1 & -2 & -4 \\
 -1 & -2 & -2 \\
}$
 & -2 } \\
 \stab{$S3_{71}$&
$\smat{
 2 & -3 & 0 \\
 -1 & 2 & -2 \\
 0 & -1 & 1 \\
}$
 & $\frac{1}{3}$ 
$\smat{
 0 & -3 & -6 \\
 -1 & -2 & -4 \\
 -1 & -2 & -1 \\
}$
 & -3 } \\
 \stab{$S3_{72}$&
$\smat{
 2 & -4 & 0 \\
 -1 & 2 & -2 \\
 0 & -1 & 1 \\
}$
 & $\frac{1}{4}$ 
$\smat{
 0 & -4 & -8 \\
 -1 & -2 & -4 \\
 -1 & -2 & 0 \\
}$
 & -4 } \\
 \stab{$S3_{73}$&
$\smat{
 1 & -2 & 0 \\
 -1 & 2 & -1 \\
 0 & -1 & 1 \\
}$
 & 
$\smat{
 -1 & -2 & -2 \\
 -1 & -1 & -1 \\
 -1 & -1 & 0 \\
}$
 & -1 } \\
 \stab{$S3_{74}$&
$\smat{
 1 & -1 & 0 \\
 -2 & 2 & -1 \\
 0 & -1 & 1 \\
}$
 & 
$\smat{
 -1 & -1 & -1 \\
 -2 & -1 & -1 \\
 -2 & -1 & 0 \\
}$
 & -1 } \\
 \stab{$S3_{75}$&
$\smat{
 1 & -2 & 0 \\
 -1 & 2 & -1 \\
 0 & -2 & 1 \\
}$
 & $\frac{1}{2}$ 
$\smat{
 0 & -2 & -2 \\
 -1 & -1 & -1 \\
 -2 & -2 & 0 \\
}$
 & -2 } \\
 \stab{$S3_{76}$&
$\smat{
 1 & -1 & 0 \\
 -2 & 2 & -1 \\
 0 & -2 & 1 \\
}$
 & $\frac{1}{2}$ 
$\smat{
 0 & -1 & -1 \\
 -2 & -1 & -1 \\
 -4 & -2 & 0 \\
}$
 & -2 } \\
 \stab{$S3_{77}$&
$\smat{
 1 & -1 & 0 \\
 -2 & 2 & -2 \\
 0 & -1 & 1 \\
}$
 & $\frac{1}{2}$ 
$\smat{
 0 & -1 & -2 \\
 -2 & -1 & -2 \\
 -2 & -1 & 0 \\
}$
 & -2 } \\
 \stab{$S3_{78}$&
$\smat{
 2 & -1 & 0 \\
 -1 & 1 & -1 \\
 0 & -1 & 1 \\
}$
 & 
$\smat{
 0 & -1 & -1 \\
 -1 & -2 & -2 \\
 -1 & -2 & -1 \\
}$
 & -1 } \\
 \stab{$S3_{79}$&
$\smat{
 2 & -1 & 0 \\
 -2 & 1 & -1 \\
 0 & -1 & 1 \\
}$
 & $\frac{1}{2}$ 
$\smat{
 0 & -1 & -1 \\
 -2 & -2 & -2 \\
 -2 & -2 & 0 \\
}$
 & -2 } \\
 \stab{$S3_{80}$&
$\smat{
 2 & -2 & 0 \\
 -1 & 1 & -1 \\
 0 & -1 & 1 \\
}$
 & $\frac{1}{2}$ 
$\smat{
 0 & -2 & -2 \\
 -1 & -2 & -2 \\
 -1 & -2 & 0 \\
}$
 & -2 } \\
 \stab{$S3_{81}$&
$\smat{
 1 & -1 & 0 \\
 -1 & 1 & -1 \\
 0 & -1 & 1 \\
}$
 & 
$\smat{
 0 & -1 & -1 \\
 -1 & -1 & -1 \\
 -1 & -1 & 0 \\
}$
 & -1 } \\
 \stab{$S4_{1}$&
$\smat{
 2 & -1 & 0 & 0 \\
 -1 & 1 & -1 & 0 \\
 0 & -1 & 2 & -1 \\
 0 & 0 & -1 & 2 \\
}$
 & 
$\smat{
 -1 & -3 & -2 & -1 \\
 -3 & -6 & -4 & -2 \\
 -2 & -4 & -2 & -1 \\
 -1 & -2 & -1 & 0 \\
}$
 & -1 } \\
 \stab{$S4_{2}$&
$\smat{
 2 & -1 & 0 & 0 \\
 -1 & 1 & -1 & 0 \\
 0 & -1 & 2 & -2 \\
 0 & 0 & -1 & 2 \\
}$
 & $\frac{1}{2}$ 
$\smat{
 0 & -2 & -2 & -2 \\
 -2 & -4 & -4 & -4 \\
 -2 & -4 & -2 & -2 \\
 -1 & -2 & -1 & 0 \\
}$
 & -2 } \\
 \stab{$S4_{3}$&
$\smat{
 2 & -1 & 0 & 0 \\
 -1 & 1 & -1 & 0 \\
 0 & -1 & 2 & -1 \\
 0 & 0 & -2 & 2 \\
}$
 & $\frac{1}{2}$ 
$\smat{
 0 & -2 & -2 & -1 \\
 -2 & -4 & -4 & -2 \\
 -2 & -4 & -2 & -1 \\
 -2 & -4 & -2 & 0 \\
}$
 & -2 } \\
 \stab{$S4_{4}$&
$\smat{
 2 & -1 & 0 & 0 \\
 -3 & 2 & -1 & 0 \\
 0 & -1 & 2 & -1 \\
 0 & 0 & -1 & 1 \\
}$
 & 
$\smat{
 -1 & -1 & -1 & -1 \\
 -3 & -2 & -2 & -2 \\
 -3 & -2 & -1 & -1 \\
 -3 & -2 & -1 & 0 \\
}$
 & -1 } \\
 \stab{$S4_{5}$&
$\smat{
 2 & -3 & 0 & 0 \\
 -1 & 2 & -1 & 0 \\
 0 & -1 & 2 & -1 \\
 0 & 0 & -1 & 1 \\
}$
 & 
$\smat{
 -1 & -3 & -3 & -3 \\
 -1 & -2 & -2 & -2 \\
 -1 & -2 & -1 & -1 \\
 -1 & -2 & -1 & 0 \\
}$
 & -1 } \\
 \stab{$S4_{6}$&
$\smat{
 2 & -1 & 0 & 0 \\
 -2 & 2 & -1 & 0 \\
 0 & -2 & 2 & -1 \\
 0 & 0 & -1 & 1 \\
}$
 & $\frac{1}{2}$ 
$\smat{
 0 & -1 & -1 & -1 \\
 -2 & -2 & -2 & -2 \\
 -4 & -4 & -2 & -2 \\
 -4 & -4 & -2 & 0 \\
}$
 & -2 } \\
 \stab{$S4_{7}$&
$\smat{
 2 & -1 & 0 & 0 \\
 -1 & 2 & -1 & 0 \\
 0 & -2 & 2 & -1 \\
 0 & 0 & -1 & 1 \\
}$
 & 
$\smat{
 0 & -1 & -1 & -1 \\
 -1 & -2 & -2 & -2 \\
 -2 & -4 & -3 & -3 \\
 -2 & -4 & -3 & -2 \\
}$
 & -1 } \\
 \stab{$S4_{8}$&
$\smat{
 2 & -2 & 0 & 0 \\
 -1 & 2 & -1 & 0 \\
 0 & -2 & 2 & -1 \\
 0 & 0 & -1 & 1 \\
}$
 & $\frac{1}{2}$ 
$\smat{
 0 & -2 & -2 & -2 \\
 -1 & -2 & -2 & -2 \\
 -2 & -4 & -2 & -2 \\
 -2 & -4 & -2 & 0 \\
}$
 & -2 } \\
 \stab{$S4_{9}$&
$\smat{
 2 & -1 & 0 & 0 \\
 -2 & 2 & -2 & 0 \\
 0 & -1 & 2 & -1 \\
 0 & 0 & -1 & 1 \\
}$
 & $\frac{1}{2}$ 
$\smat{
 0 & -1 & -2 & -2 \\
 -2 & -2 & -4 & -4 \\
 -2 & -2 & -2 & -2 \\
 -2 & -2 & -2 & 0 \\
}$
 & -2 } \\
 \stab{$S4_{10}$&
$\smat{
 2 & -1 & 0 & 0 \\
 -1 & 2 & -2 & 0 \\
 0 & -1 & 2 & -1 \\
 0 & 0 & -1 & 1 \\
}$
 & 
$\smat{
 0 & -1 & -2 & -2 \\
 -1 & -2 & -4 & -4 \\
 -1 & -2 & -3 & -3 \\
 -1 & -2 & -3 & -2 \\
}$
 & -1 } \\
 \stab{$S4_{11}$&
$\smat{
 2 & -2 & 0 & 0 \\
 -1 & 2 & -2 & 0 \\
 0 & -1 & 2 & -1 \\
 0 & 0 & -1 & 1 \\
}$
 & $\frac{1}{2}$ 
$\smat{
 0 & -2 & -4 & -4 \\
 -1 & -2 & -4 & -4 \\
 -1 & -2 & -2 & -2 \\
 -1 & -2 & -2 & 0 \\
}$
 & -2 } \\
 \stab{$S4_{12}$&
$\smat{
 1 & -1 & 0 & 0 \\
 -1 & 2 & -1 & 0 \\
 0 & -2 & 2 & -1 \\
 0 & 0 & -1 & 1 \\
}$
 & 
$\smat{
 0 & -1 & -1 & -1 \\
 -1 & -1 & -1 & -1 \\
 -2 & -2 & -1 & -1 \\
 -2 & -2 & -1 & 0 \\
}$
 & -1 } \\
 \stab{$S4_{13}$&
$\smat{
 2 & -1 & 0 & 0 \\
 -1 & 1 & -1 & 0 \\
 0 & -1 & 2 & -1 \\
 0 & 0 & -1 & 1 \\
}$
 & 
$\smat{
 0 & -1 & -1 & -1 \\
 -1 & -2 & -2 & -2 \\
 -1 & -2 & -1 & -1 \\
 -1 & -2 & -1 & 0 \\
}$
 & -1 } \\
 \stab{$S4_{14}$&
$\smat{
 2 & -1 & 0 & -1 \\
 -1 & 1 & -1 & 0 \\
 0 & -1 & 2 & -1 \\
 -1 & 0 & -1 & 2 \\
}$
 & $\frac{1}{4}$ 
$\smat{
 -1 & -4 & -3 & -2 \\
 -4 & -4 & -4 & -4 \\
 -3 & -4 & -1 & -2 \\
 -2 & -4 & -2 & 0 \\
}$
 & -4 } \\
 \stab{$NS4_{15}$&
$\smat{
 2 & -1 & 0 & -2 \\
 -1 & 1 & -1 & 0 \\
 0 & -1 & 2 & -1 \\
 -1 & 0 & -1 & 2 \\
}$
 & $\frac{1}{6}$ 
$\smat{
 -1 & -5 & -4 & -3 \\
 -4 & -2 & -4 & -6 \\
 -3 & -3 & 0 & -3 \\
 -2 & -4 & -2 & 0 \\
}$
 & -6 } \\
 \stab{$NS4_{16}$&
$\smat{
 2 & -1 & 0 & -1 \\
 -1 & 1 & -1 & 0 \\
 0 & -1 & 2 & -1 \\
 -2 & 0 & -1 & 2 \\
}$
 & $\frac{1}{6}$ 
$\smat{
 -1 & -4 & -3 & -2 \\
 -5 & -2 & -3 & -4 \\
 -4 & -4 & 0 & -2 \\
 -3 & -6 & -3 & 0 \\
}$
 & -6 } \\
 \stab{$S4_{17}$&
$\smat{
 2 & -1 & 0 & -2 \\
 -1 & 1 & -1 & 0 \\
 0 & -1 & 2 & -2 \\
 -1 & 0 & -1 & 2 \\
}$
 & $\frac{1}{4}$ 
$\smat{
 0 & -2 & -2 & -2 \\
 -2 & 0 & -2 & -4 \\
 -2 & -2 & 0 & -2 \\
 -1 & -2 & -1 & 0 \\
}$
 & -8 } \\
 \stab{$NS4_{18}$&
$\smat{
 2 & -1 & 0 & -1 \\
 -1 & 1 & -1 & 0 \\
 0 & -1 & 2 & -2 \\
 -2 & 0 & -1 & 2 \\
}$
 & $\frac{1}{3}$ 
$\smat{
 0 & -1 & -1 & -1 \\
 -2 & 0 & -1 & -2 \\
 -2 & -2 & 0 & -1 \\
 -1 & -2 & -1 & 0 \\
}$
 & -9 } \\
 \stab{$S4_{19}$&
$\smat{
 2 & -1 & 0 & -1 \\
 -1 & 1 & -1 & 0 \\
 0 & -1 & 2 & -1 \\
 -2 & 0 & -2 & 2 \\
}$
 & $\frac{1}{4}$ 
$\smat{
 0 & -2 & -2 & -1 \\
 -2 & 0 & -2 & -2 \\
 -2 & -2 & 0 & -1 \\
 -2 & -4 & -2 & 0 \\
}$
 & -8 } \\
 \stab{$S4_{20}$&
$\smat{
 2 & -1 & 0 & -1 \\
 -1 & 1 & -1 & 0 \\
 0 & -1 & 2 & -1 \\
 -1 & 0 & -1 & 1 \\
}$
 & $\frac{1}{2}$ 
$\smat{
 0 & -1 & -1 & -1 \\
 -1 & 0 & -1 & -2 \\
 -1 & -1 & 0 & -1 \\
 -1 & -2 & -1 & 0 \\
}$
 & -4 } \\
 \stab{$S4_{21}$&
$\smat{
 2 & -1 & 0 & 0 \\
 -1 & 1 & -1 & -1 \\
 0 & -1 & 2 & 0 \\
 0 & -1 & 0 & 2 \\
}$
 & $\frac{1}{2}$ 
$\smat{
 0 & -2 & -1 & -1 \\
 -2 & -4 & -2 & -2 \\
 -1 & -2 & 0 & -1 \\
 -1 & -2 & -1 & 0 \\
}$
 & -4 } \\
 \stab{$S4_{22}$&
$\smat{
 2 & -1 & 0 & 0 \\
 -2 & 2 & -1 & -1 \\
 0 & -1 & 1 & 0 \\
 0 & -1 & 0 & 2 \\
}$
 & $\frac{1}{2}$ 
$\smat{
 -1 & -2 & -2 & -1 \\
 -4 & -4 & -4 & -2 \\
 -4 & -4 & -2 & -2 \\
 -2 & -2 & -2 & 0 \\
}$
 & -2 } \\
 \stab{$S4_{23}$&
$\smat{
 2 & -2 & 0 & 0 \\
 -1 & 2 & -1 & -1 \\
 0 & -1 & 1 & 0 \\
 0 & -1 & 0 & 2 \\
}$
 & $\frac{1}{2}$ 
$\smat{
 -1 & -4 & -4 & -2 \\
 -2 & -4 & -4 & -2 \\
 -2 & -4 & -2 & -2 \\
 -1 & -2 & -2 & 0 \\
}$
 & -2 } \\
 \stab{$S4_{24}$&
$\smat{
 2 & -1 & 0 & 0 \\
 -2 & 2 & -1 & -2 \\
 0 & -1 & 1 & 0 \\
 0 & -1 & 0 & 2 \\
}$
 & $\frac{1}{2}$ 
$\smat{
 0 & -1 & -1 & -1 \\
 -2 & -2 & -2 & -2 \\
 -2 & -2 & 0 & -2 \\
 -1 & -1 & -1 & 0 \\
}$
 & -4 } \\
 \stab{$S4_{25}$&
$\smat{
 2 & -2 & 0 & 0 \\
 -1 & 2 & -1 & -2 \\
 0 & -1 & 1 & 0 \\
 0 & -1 & 0 & 2 \\
}$
 & $\frac{1}{4}$ 
$\smat{
 0 & -4 & -4 & -4 \\
 -2 & -4 & -4 & -4 \\
 -2 & -4 & 0 & -4 \\
 -1 & -2 & -2 & 0 \\
}$
 & -4 } \\
 \stab{$S4_{26}$&
$\smat{
 2 & -2 & 0 & 0 \\
 -1 & 2 & -1 & -1 \\
 0 & -1 & 1 & 0 \\
 0 & -2 & 0 & 2 \\
}$
 & $\frac{1}{2}$ 
$\smat{
 0 & -2 & -2 & -1 \\
 -1 & -2 & -2 & -1 \\
 -1 & -2 & 0 & -1 \\
 -1 & -2 & -2 & 0 \\
}$
 & -4 } \\
 \stab{$S4_{27}$&
$\smat{
 2 & -1 & 0 & 0 \\
 -1 & 2 & -1 & -1 \\
 0 & -1 & 1 & 0 \\
 0 & -1 & 0 & 1 \\
}$
 & 
$\smat{
 0 & -1 & -1 & -1 \\
 -1 & -2 & -2 & -2 \\
 -1 & -2 & -1 & -2 \\
 -1 & -2 & -2 & -1 \\
}$
 & -1 } \\
 \stab{$S4_{28}$&
$\smat{
 2 & -1 & 0 & 0 \\
 -2 & 2 & -1 & -1 \\
 0 & -1 & 1 & 0 \\
 0 & -1 & 0 & 1 \\
}$
 & $\frac{1}{2}$ 
$\smat{
 0 & -1 & -1 & -1 \\
 -2 & -2 & -2 & -2 \\
 -2 & -2 & 0 & -2 \\
 -2 & -2 & -2 & 0 \\
}$
 & -2 } \\
 \stab{$S4_{29}$&
$\smat{
 2 & -2 & 0 & 0 \\
 -1 & 2 & -1 & -1 \\
 0 & -1 & 1 & 0 \\
 0 & -1 & 0 & 1 \\
}$
 & $\frac{1}{2}$ 
$\smat{
 0 & -2 & -2 & -2 \\
 -1 & -2 & -2 & -2 \\
 -1 & -2 & 0 & -2 \\
 -1 & -2 & -2 & 0 \\
}$
 & -2 } \\
 \stab{$S4_{30}$&
$\smat{
 1 & -1 & 0 & 0 \\
 -1 & 2 & -1 & -1 \\
 0 & -1 & 1 & 0 \\
 0 & -1 & 0 & 1 \\
}$
 & 
$\smat{
 0 & -1 & -1 & -1 \\
 -1 & -1 & -1 & -1 \\
 -1 & -1 & 0 & -1 \\
 -1 & -1 & -1 & 0 \\
}$
 & -1 } \\
 \stab{$S4_{31}$&
$\smat{
 1 & -1 & 0 & 0 \\
 -1 & 2 & -1 & -1 \\
 0 & -1 & 2 & -1 \\
 0 & -1 & -1 & 2 \\
}$
 & $\frac{1}{3}$ 
$\smat{
 0 & -3 & -3 & -3 \\
 -3 & -3 & -3 & -3 \\
 -3 & -3 & -1 & -2 \\
 -3 & -3 & -2 & -1 \\
}$
 & -3 } \\
 \stab{$S5_{1}$&
$\smat{
 2 & -1 & 0 & 0 & 0 \\
 -1 & 2 & -1 & 0 & 0 \\
 0 & -2 & 2 & -1 & 0 \\
 0 & 0 & -1 & 2 & -1 \\
 0 & 0 & 0 & -1 & 1 \\
}$
 & 
$\smat{
 0 & -1 & -1 & -1 & -1 \\
 -1 & -2 & -2 & -2 & -2 \\
 -2 & -4 & -3 & -3 & -3 \\
 -2 & -4 & -3 & -2 & -2 \\
 -2 & -4 & -3 & -2 & -1 \\
}$
 & -1 } \\
 \stab{$S5_{2}$&
$\smat{
 2 & -1 & 0 & 0 & 0 \\
 -1 & 2 & -2 & 0 & 0 \\
 0 & -1 & 2 & -1 & 0 \\
 0 & 0 & -1 & 2 & -1 \\
 0 & 0 & 0 & -1 & 1 \\
}$
 & 
$\smat{
 0 & -1 & -2 & -2 & -2 \\
 -1 & -2 & -4 & -4 & -4 \\
 -1 & -2 & -3 & -3 & -3 \\
 -1 & -2 & -3 & -2 & -2 \\
 -1 & -2 & -3 & -2 & -1 \\
}$
 & -1 } \\
 \stab{$S5_{3}$&
$\smat{
 2 & -1 & 0 & 0 & 0 \\
 -1 & 2 & -1 & -1 & -1 \\
 0 & -1 & 1 & 0 & 0 \\
 0 & -1 & 0 & 2 & 0 \\
 0 & -1 & 0 & 0 & 2 \\
}$
 & $\frac{1}{2}$ 
$\smat{
 0 & -2 & -2 & -1 & -1 \\
 -2 & -4 & -4 & -2 & -2 \\
 -2 & -4 & -2 & -2 & -2 \\
 -1 & -2 & -2 & 0 & -1 \\
 -1 & -2 & -2 & -1 & 0 \\
}$
 & -4 } \\
 \stab{$S5_{4}$&
$\smat{
 2 & -1 & 0 & 0 & 0 \\
 -2 & 2 & -1 & -1 & 0 \\
 0 & -1 & 2 & 0 & 0 \\
 0 & -1 & 0 & 2 & -1 \\
 0 & 0 & 0 & -1 & 1 \\
}$
 & $\frac{1}{2}$ 
$\smat{
 -1 & -2 & -1 & -2 & -2 \\
 -4 & -4 & -2 & -4 & -4 \\
 -2 & -2 & 0 & -2 & -2 \\
 -4 & -4 & -2 & -2 & -2 \\
 -4 & -4 & -2 & -2 & 0 \\
}$
 & -2 } \\
 \stab{$S5_{5}$&
$\smat{
 2 & -2 & 0 & 0 & 0 \\
 -1 & 2 & -1 & -1 & 0 \\
 0 & -1 & 2 & 0 & 0 \\
 0 & -1 & 0 & 2 & -1 \\
 0 & 0 & 0 & -1 & 1 \\
}$
 & $\frac{1}{2}$ 
$\smat{
 -1 & -4 & -2 & -4 & -4 \\
 -2 & -4 & -2 & -4 & -4 \\
 -1 & -2 & 0 & -2 & -2 \\
 -2 & -4 & -2 & -2 & -2 \\
 -2 & -4 & -2 & -2 & 0 \\
}$
 & -2 } \\
 \stab{$S5_{6}$&
$\smat{
 2 & -1 & 0 & 0 & 0 \\
 -1 & 2 & -1 & -1 & 0 \\
 0 & -1 & 1 & 0 & 0 \\
 0 & -1 & 0 & 2 & -1 \\
 0 & 0 & 0 & -1 & 2 \\
}$
 & 
$\smat{
 -1 & -3 & -3 & -2 & -1 \\
 -3 & -6 & -6 & -4 & -2 \\
 -3 & -6 & -5 & -4 & -2 \\
 -2 & -4 & -4 & -2 & -1 \\
 -1 & -2 & -2 & -1 & 0 \\
}$
 & -1 } \\
 \stab{$S5_{7}$&
$\smat{
 2 & -1 & 0 & 0 & 0 \\
 -1 & 2 & -1 & -1 & 0 \\
 0 & -1 & 1 & 0 & 0 \\
 0 & -1 & 0 & 2 & -2 \\
 0 & 0 & 0 & -1 & 2 \\
}$
 & $\frac{1}{2}$ 
$\smat{
 0 & -2 & -2 & -2 & -2 \\
 -2 & -4 & -4 & -4 & -4 \\
 -2 & -4 & -2 & -4 & -4 \\
 -2 & -4 & -4 & -2 & -2 \\
 -1 & -2 & -2 & -1 & 0 \\
}$
 & -2 } \\
 \stab{$S5_{8}$&
$\smat{
 2 & -1 & 0 & 0 & 0 \\
 -1 & 2 & -1 & -1 & 0 \\
 0 & -1 & 1 & 0 & 0 \\
 0 & -1 & 0 & 2 & -1 \\
 0 & 0 & 0 & -2 & 2 \\
}$
 & $\frac{1}{2}$ 
$\smat{
 0 & -2 & -2 & -2 & -1 \\
 -2 & -4 & -4 & -4 & -2 \\
 -2 & -4 & -2 & -4 & -2 \\
 -2 & -4 & -4 & -2 & -1 \\
 -2 & -4 & -4 & -2 & 0 \\
}$
 & -2 } \\
 \stab{$S5_{9}$&
$\smat{
 2 & -1 & 0 & 0 & 0 \\
 -1 & 2 & -1 & -1 & 0 \\
 0 & -1 & 1 & 0 & 0 \\
 0 & -1 & 0 & 2 & -1 \\
 0 & 0 & 0 & -1 & 1 \\
}$
 & 
$\smat{
 0 & -1 & -1 & -1 & -1 \\
 -1 & -2 & -2 & -2 & -2 \\
 -1 & -2 & -1 & -2 & -2 \\
 -1 & -2 & -2 & -1 & -1 \\
 -1 & -2 & -2 & -1 & 0 \\
}$
 & -1 } \\
 \stab{$S5_{10}$&
$\smat{
 1 & -1 & 0 & 0 & 0 \\
 -1 & 2 & -1 & 0 & -1 \\
 0 & -1 & 2 & -1 & 0 \\
 0 & 0 & -1 & 2 & -1 \\
 0 & -1 & 0 & -1 & 2 \\
}$
 & $\frac{1}{4}$ 
$\smat{
 0 & -4 & -4 & -4 & -4 \\
 -4 & -4 & -4 & -4 & -4 \\
 -4 & -4 & -1 & -2 & -3 \\
 -4 & -4 & -2 & 0 & -2 \\
 -4 & -4 & -3 & -2 & -1 \\
}$
 & -4 } \\
 \stab{$S6_{1}$&
$\smat{
 2 & -1 & 0 & 0 & 0 & 0 \\
 -1 & 2 & -1 & -1 & -1 & 0 \\
 0 & -1 & 2 & 0 & 0 & 0 \\
 0 & -1 & 0 & 2 & 0 & 0 \\
 0 & -1 & 0 & 0 & 2 & -1 \\
 0 & 0 & 0 & 0 & -1 & 1 \\
}$
 & $\frac{1}{2}$ 
$\smat{
 0 & -2 & -1 & -1 & -2 & -2 \\
 -2 & -4 & -2 & -2 & -4 & -4 \\
 -1 & -2 & 0 & -1 & -2 & -2 \\
 -1 & -2 & -1 & 0 & -2 & -2 \\
 -2 & -4 & -2 & -2 & -2 & -2 \\
 -2 & -4 & -2 & -2 & -2 & 0 \\
}$
 & -4 } \\
 \stab{$S6_{2}$&
$\smat{
 2 & -1 & 0 & 0 & 0 & 0 \\
 -1 & 2 & -1 & 0 & 0 & 0 \\
 0 & -1 & 2 & -1 & -1 & 0 \\
 0 & 0 & -1 & 2 & 0 & 0 \\
 0 & 0 & -1 & 0 & 2 & -1 \\
 0 & 0 & 0 & 0 & -1 & 1 \\
}$
 & 
$\smat{
 0 & -1 & -2 & -1 & -2 & -2 \\
 -1 & -2 & -4 & -2 & -4 & -4 \\
 -2 & -4 & -6 & -3 & -6 & -6 \\
 -1 & -2 & -3 & -1 & -3 & -3 \\
 -2 & -4 & -6 & -3 & -5 & -5 \\
 -2 & -4 & -6 & -3 & -5 & -4 \\
}$
 & -1 } \\
 \stab{$S6_{3}$&
$\smat{
 2 & -1 & 0 & 0 & 0 & 0 \\
 -2 & 2 & -1 & 0 & 0 & 0 \\
 0 & -1 & 2 & -1 & -1 & 0 \\
 0 & 0 & -1 & 2 & 0 & 0 \\
 0 & 0 & -1 & 0 & 2 & -1 \\
 0 & 0 & 0 & 0 & -1 & 1 \\
}$
 & $\frac{1}{2}$ 
$\smat{
 0 & -1 & -2 & -1 & -2 & -2 \\
 -2 & -2 & -4 & -2 & -4 & -4 \\
 -4 & -4 & -4 & -2 & -4 & -4 \\
 -2 & -2 & -2 & 0 & -2 & -2 \\
 -4 & -4 & -4 & -2 & -2 & -2 \\
 -4 & -4 & -4 & -2 & -2 & 0 \\
}$
 & -2 } \\
 \stab{$S6_{4}$&
$\smat{
 2 & -2 & 0 & 0 & 0 & 0 \\
 -1 & 2 & -1 & 0 & 0 & 0 \\
 0 & -1 & 2 & -1 & -1 & 0 \\
 0 & 0 & -1 & 2 & 0 & 0 \\
 0 & 0 & -1 & 0 & 2 & -1 \\
 0 & 0 & 0 & 0 & -1 & 1 \\
}$
 & $\frac{1}{2}$ 
$\smat{
 0 & -2 & -4 & -2 & -4 & -4 \\
 -1 & -2 & -4 & -2 & -4 & -4 \\
 -2 & -4 & -4 & -2 & -4 & -4 \\
 -1 & -2 & -2 & 0 & -2 & -2 \\
 -2 & -4 & -4 & -2 & -2 & -2 \\
 -2 & -4 & -4 & -2 & -2 & 0 \\
}$
 & -2 } \\
 \stab{$S6_{5}$&
$\smat{
 1 & -1 & 0 & 0 & 0 & 0 \\
 -1 & 2 & -1 & 0 & 0 & 0 \\
 0 & -1 & 2 & -1 & -1 & 0 \\
 0 & 0 & -1 & 2 & 0 & 0 \\
 0 & 0 & -1 & 0 & 2 & -1 \\
 0 & 0 & 0 & 0 & -1 & 1 \\
}$
 & 
$\smat{
 0 & -1 & -2 & -1 & -2 & -2 \\
 -1 & -1 & -2 & -1 & -2 & -2 \\
 -2 & -2 & -2 & -1 & -2 & -2 \\
 -1 & -1 & -1 & 0 & -1 & -1 \\
 -2 & -2 & -2 & -1 & -1 & -1 \\
 -2 & -2 & -2 & -1 & -1 & 0 \\
}$
 & -1 } \\
 \stab{$S6_{6}$&
$\smat{
 2 & -1 & 0 & 0 & 0 & 0 \\
 -1 & 2 & -1 & 0 & 0 & 0 \\
 0 & -2 & 2 & -1 & 0 & 0 \\
 0 & 0 & -1 & 2 & -1 & 0 \\
 0 & 0 & 0 & -1 & 2 & -1 \\
 0 & 0 & 0 & 0 & -1 & 1 \\
}$
 & 
$\smat{
 0 & -1 & -1 & -1 & -1 & -1 \\
 -1 & -2 & -2 & -2 & -2 & -2 \\
 -2 & -4 & -3 & -3 & -3 & -3 \\
 -2 & -4 & -3 & -2 & -2 & -2 \\
 -2 & -4 & -3 & -2 & -1 & -1 \\
 -2 & -4 & -3 & -2 & -1 & 0 \\
}$
 & -1 } \\
 \stab{$S6_{7}$&
$\smat{
 2 & -1 & 0 & 0 & 0 & 0 \\
 -1 & 2 & -2 & 0 & 0 & 0 \\
 0 & -1 & 2 & -1 & 0 & 0 \\
 0 & 0 & -1 & 2 & -1 & 0 \\
 0 & 0 & 0 & -1 & 2 & -1 \\
 0 & 0 & 0 & 0 & -1 & 1 \\
}$
 & 
$\smat{
 0 & -1 & -2 & -2 & -2 & -2 \\
 -1 & -2 & -4 & -4 & -4 & -4 \\
 -1 & -2 & -3 & -3 & -3 & -3 \\
 -1 & -2 & -3 & -2 & -2 & -2 \\
 -1 & -2 & -3 & -2 & -1 & -1 \\
 -1 & -2 & -3 & -2 & -1 & 0 \\
}$
 & -1 } \\
 \stab{$S7_1$&
$\smat{
 2 & -1 & 0 & 0 & 0 & 0 & 0 \\
 -1 & 2 & -1 & 0 & 0 & 0 & 0 \\
 0 & -1 & 2 & -1 & -1 & 0 & 0 \\
 0 & 0 & -1 & 2 & 0 & 0 & 0 \\
 0 & 0 & -1 & 0 & 2 & -1 & 0 \\
 0 & 0 & 0 & 0 & -1 & 2 & -1 \\
 0 & 0 & 0 & 0 & 0 & -1 & 1 \\
}$
 & 
$\smat{
 0 & -1 & -2 & -1 & -2 & -2 & -2 \\
 -1 & -2 & -4 & -2 & -4 & -4 & -4 \\
 -2 & -4 & -6 & -3 & -6 & -6 & -6 \\
 -1 & -2 & -3 & -1 & -3 & -3 & -3 \\
 -2 & -4 & -6 & -3 & -5 & -5 & -5 \\
 -2 & -4 & -6 & -3 & -5 & -4 & -4 \\
 -2 & -4 & -6 & -3 & -5 & -4 & -3 \\
}$
 & -1 } \\
 \stab{$S8_1$&
$\smat{
 2 & -1 & 0 & 0 & 0 & 0 & 0 & 0 \\
 -1 & 2 & -1 & 0 & 0 & 0 & 0 & 0 \\
 0 & -1 & 2 & -1 & -1 & 0 & 0 & 0 \\
 0 & 0 & -1 & 2 & 0 & 0 & 0 & 0 \\
 0 & 0 & -1 & 0 & 2 & -1 & 0 & 0 \\
 0 & 0 & 0 & 0 & -1 & 2 & -1 & 0 \\
 0 & 0 & 0 & 0 & 0 & -1 & 2 & -1 \\
 0 & 0 & 0 & 0 & 0 & 0 & -1 & 1 \\
}$
 & 
$\smat{
 0 & -1 & -2 & -1 & -2 & -2 & -2 & -2 \\
 -1 & -2 & -4 & -2 & -4 & -4 & -4 & -4 \\
 -2 & -4 & -6 & -3 & -6 & -6 & -6 & -6 \\
 -1 & -2 & -3 & -1 & -3 & -3 & -3 & -3 \\
 -2 & -4 & -6 & -3 & -5 & -5 & -5 & -5 \\
 -2 & -4 & -6 & -3 & -5 & -4 & -4 & -4 \\
 -2 & -4 & -6 & -3 & -5 & -4 & -3 & -3 \\
 -2 & -4 & -6 & -3 & -5 & -4 & -3 & -2 \\
}$
 & -1 } \\
 \stab{$S9_1$&
$\smat{
 2 & -1 & 0 & 0 & 0 & 0 & 0 & 0 & 0 \\
 -1 & 2 & -1 & 0 & 0 & 0 & 0 & 0 & 0 \\
 0 & -1 & 2 & -1 & -1 & 0 & 0 & 0 & 0 \\
 0 & 0 & -1 & 2 & 0 & 0 & 0 & 0 & 0 \\
 0 & 0 & -1 & 0 & 2 & -1 & 0 & 0 & 0 \\
 0 & 0 & 0 & 0 & -1 & 2 & -1 & 0 & 0 \\
 0 & 0 & 0 & 0 & 0 & -1 & 2 & -1 & 0 \\
 0 & 0 & 0 & 0 & 0 & 0 & -1 & 2 & -1 \\
 0 & 0 & 0 & 0 & 0 & 0 & 0 & -1 & 1 \\
}$
 & 
$\smat{
 0 & -1 & -2 & -1 & -2 & -2 & -2 & -2 & -2 \\
 -1 & -2 & -4 & -2 & -4 & -4 & -4 & -4 & -4 \\
 -2 & -4 & -6 & -3 & -6 & -6 & -6 & -6 & -6 \\
 -1 & -2 & -3 & -1 & -3 & -3 & -3 & -3 & -3 \\
 -2 & -4 & -6 & -3 & -5 & -5 & -5 & -5 & -5 \\
 -2 & -4 & -6 & -3 & -5 & -4 & -4 & -4 & -4 \\
 -2 & -4 & -6 & -3 & -5 & -4 & -3 & -3 & -3 \\
 -2 & -4 & -6 & -3 & -5 & -4 & -3 & -2 & -2 \\
 -2 & -4 & -6 & -3 & -5 & -4 & -3 & -2 & -1 \\
}$
 & -1 } \\
 \stab{$S10_1$&
$\smat{
 2 & -1 & 0 & 0 & 0 & 0 & 0 & 0 & 0 & 0 \\
 -1 & 2 & -1 & 0 & 0 & 0 & 0 & 0 & 0 & 0 \\
 0 & -1 & 2 & -1 & -1 & 0 & 0 & 0 & 0 & 0 \\
 0 & 0 & -1 & 2 & 0 & 0 & 0 & 0 & 0 & 0 \\
 0 & 0 & -1 & 0 & 2 & -1 & 0 & 0 & 0 & 0 \\
 0 & 0 & 0 & 0 & -1 & 2 & -1 & 0 & 0 & 0 \\
 0 & 0 & 0 & 0 & 0 & -1 & 2 & -1 & 0 & 0 \\
 0 & 0 & 0 & 0 & 0 & 0 & -1 & 2 & -1 & 0 \\
 0 & 0 & 0 & 0 & 0 & 0 & 0 & -1 & 2 & -1 \\
 0 & 0 & 0 & 0 & 0 & 0 & 0 & 0 & -1 & 1 \\
}$
 & 
$\smat{
 0 & -1 & -2 & -1 & -2 & -2 & -2 & -2 & -2 & -2 \\
 -1 & -2 & -4 & -2 & -4 & -4 & -4 & -4 & -4 & -4 \\
 -2 & -4 & -6 & -3 & -6 & -6 & -6 & -6 & -6 & -6 \\
 -1 & -2 & -3 & -1 & -3 & -3 & -3 & -3 & -3 & -3 \\
 -2 & -4 & -6 & -3 & -5 & -5 & -5 & -5 & -5 & -5 \\
 -2 & -4 & -6 & -3 & -5 & -4 & -4 & -4 & -4 & -4 \\
 -2 & -4 & -6 & -3 & -5 & -4 & -3 & -3 & -3 & -3 \\
 -2 & -4 & -6 & -3 & -5 & -4 & -3 & -2 & -2 & -2 \\
 -2 & -4 & -6 & -3 & -5 & -4 & -3 & -2 & -1 & -1 \\
 -2 & -4 & -6 & -3 & -5 & -4 & -3 & -2 & -1 & 0 \\
}$
 & -1 } \\
 
\subsection{Cartan matrix $A$, $A^{-1}$, and $\det A$ for hyperbolic Lie algebras}\label{ssHyp}{}~{}
%\[
%\begin{matrix}

\noindent  \stab{\attn$NH3_1$&
$\smat{
 2 & -1 & -1 \\
 -2 & 2 & -1 \\
 -2 & -2 & 2 \\
}$
 & $\frac{1}{10}$ 
$\smat{
 -2 & -4 & -3 \\
 -6 & -2 & -4 \\
 -8 & -6 & -2 \\
}$
 & -10 } \\
 \stab{\attn$NH3_2$&
$\smat{
 2 & -1 & -1 \\
 -3 & 2 & -1 \\
 -2 & -2 & 2 \\
}$
 & $\frac{1}{14}$ 
$\smat{
 -2 & -4 & -3 \\
 -8 & -2 & -5 \\
 -10 & -6 & -1 \\
}$
 & -14 } \\
 \stab{$NH3_3$&
$\smat{
 2 & -1 & -1 \\
 -4 & 2 & -1 \\
 -2 & -2 & 2 \\
}$
 & $\frac{1}{18}$ 
$\smat{
 -2 & -4 & -3 \\
 -10 & -2 & -6 \\
 -12 & -6 & 0 \\
}$
 & -18 } \\
 \stab{$H3_4$&
$\smat{
 2 & -2 & -1 \\
 -2 & 2 & -1 \\
 -2 & -2 & 2 \\
}$
 & $\frac{1}{8}$ 
$\smat{
 -1 & -3 & -2 \\
 -3 & -1 & -2 \\
 -4 & -4 & 0 \\
}$
 & -16 } \\
 \stab{\attn$H3_5$&
$\smat{
 2 & -1 & -1 \\
 -1 & 2 & -1 \\
 -2 & -2 & 2 \\
}$
 & $\frac{1}{6}$ 
$\smat{
 -2 & -4 & -3 \\
 -4 & -2 & -3 \\
 -6 & -6 & -3 \\
}$
 & -6 } \\
 \stab{\attn$NH3_6$&
$\smat{
 2 & -1 & -1 \\
 -2 & 2 & -1 \\
 -3 & -2 & 2 \\
}$
 & $\frac{1}{13}$ 
$\smat{
 -2 & -4 & -3 \\
 -7 & -1 & -4 \\
 -10 & -7 & -2 \\
}$
 & -13 } \\
 \stab{\attn$NH3_7$&
$\smat{
 2 & -1 & -1 \\
 -3 & 2 & -1 \\
 -3 & -2 & 2 \\
}$
 & $\frac{1}{17}$ 
$\smat{
 -2 & -4 & -3 \\
 -9 & -1 & -5 \\
 -12 & -7 & -1 \\
}$
 & -17 } \\
 \stab{$NH3_8$&
$\smat{
 2 & -1 & -1 \\
 -4 & 2 & -1 \\
 -3 & -2 & 2 \\
}$
 & $\frac{1}{21}$ 
$\smat{
 -2 & -4 & -3 \\
 -11 & -1 & -6 \\
 -14 & -7 & 0 \\
}$
 & -21 } \\
 \stab{$NH3_9$&
$\smat{
 2 & -2 & -1 \\
 -2 & 2 & -1 \\
 -3 & -2 & 2 \\
}$
 & $\frac{1}{20}$ 
$\smat{
 -2 & -6 & -4 \\
 -7 & -1 & -4 \\
 -10 & -10 & 0 \\
}$
 & -20 } \\
 \stab{\attn$NH3_{10}$&
$\smat{
 2 & -1 & -1 \\
 -1 & 2 & -1 \\
 -3 & -2 & 2 \\
}$
 & $\frac{1}{9}$ 
$\smat{
 -2 & -4 & -3 \\
 -5 & -1 & -3 \\
 -8 & -7 & -3 \\
}$
 & -9 } \\
 \stab{\attn$NH3_{11}$&
$\smat{
 2 & -2 & -1 \\
 -1 & 2 & -1 \\
 -3 & -2 & 2 \\
}$
 & $\frac{1}{14}$ 
$\smat{
 -2 & -6 & -4 \\
 -5 & -1 & -3 \\
 -8 & -10 & -2 \\
}$
 & -14 } \\
 \stab{\attn$NH3_{12}$&
$\smat{
 2 & -3 & -1 \\
 -1 & 2 & -1 \\
 -3 & -2 & 2 \\
}$
 & $\frac{1}{19}$ 
$\smat{
 -2 & -8 & -5 \\
 -5 & -1 & -3 \\
 -8 & -13 & -1 \\
}$
 & -19 } \\
 \stab{$NH3_{13}$&
$\smat{
 2 & -4 & -1 \\
 -1 & 2 & -1 \\
 -3 & -2 & 2 \\
}$
 & $\frac{1}{24}$ 
$\smat{
 -2 & -10 & -6 \\
 -5 & -1 & -3 \\
 -8 & -16 & 0 \\
}$
 & -24 } \\
 \stab{$H3_{14}$&
$\smat{
 2 & -1 & -1 \\
 -2 & 2 & -1 \\
 -4 & -2 & 2 \\
}$
 & $\frac{1}{16}$ 
$\smat{
 -2 & -4 & -3 \\
 -8 & 0 & -4 \\
 -12 & -8 & -2 \\
}$
 & -16 } \\
 \stab{$NH3_{15}$&
$\smat{
 2 & -1 & -1 \\
 -3 & 2 & -1 \\
 -4 & -2 & 2 \\
}$
 & $\frac{1}{20}$ 
$\smat{
 -2 & -4 & -3 \\
 -10 & 0 & -5 \\
 -14 & -8 & -1 \\
}$
 & -20 } \\
 \stab{$NH3_{16}$&
$\smat{
 2 & -1 & -1 \\
 -4 & 2 & -1 \\
 -4 & -2 & 2 \\
}$
 & $\frac{1}{24}$ 
$\smat{
 -2 & -4 & -3 \\
 -12 & 0 & -6 \\
 -16 & -8 & 0 \\
}$
 & -24 } \\
 \stab{$NH3_{17}$&
$\smat{
 2 & -2 & -1 \\
 -2 & 2 & -1 \\
 -4 & -2 & 2 \\
}$
 & $\frac{1}{12}$ 
$\smat{
 -1 & -3 & -2 \\
 -4 & 0 & -2 \\
 -6 & -6 & 0 \\
}$
 & -24 } \\
 \stab{$NH3_{18}$&
$\smat{
 2 & -1 & -1 \\
 -1 & 2 & -1 \\
 -4 & -2 & 2 \\
}$
 & $\frac{1}{12}$ 
$\smat{
 -2 & -4 & -3 \\
 -6 & 0 & -3 \\
 -10 & -8 & -3 \\
}$
 & -12 } \\
 \stab{$NH3_{19}$&
$\smat{
 2 & -2 & -1 \\
 -1 & 2 & -1 \\
 -4 & -2 & 2 \\
}$
 & $\frac{1}{18}$ 
$\smat{
 -2 & -6 & -4 \\
 -6 & 0 & -3 \\
 -10 & -12 & -2 \\
}$
 & -18 } \\
 \stab{$NH3_{20}$&
$\smat{
 2 & -3 & -1 \\
 -1 & 2 & -1 \\
 -4 & -2 & 2 \\
}$
 & $\frac{1}{24}$ 
$\smat{
 -2 & -8 & -5 \\
 -6 & 0 & -3 \\
 -10 & -16 & -1 \\
}$
 & -24 } \\
 \stab{$NH3_{21}$&
$\smat{
 2 & -4 & -1 \\
 -1 & 2 & -1 \\
 -4 & -2 & 2 \\
}$
 & $\frac{1}{30}$ 
$\smat{
 -2 & -10 & -6 \\
 -6 & 0 & -3 \\
 -10 & -20 & 0 \\
}$
 & -30 } \\
 \stab{$NH3_{22}$&
$\smat{
 2 & -1 & -2 \\
 -2 & 2 & -1 \\
 -2 & -2 & 2 \\
}$
 & $\frac{1}{18}$ 
$\smat{
 -2 & -6 & -5 \\
 -6 & 0 & -6 \\
 -8 & -6 & -2 \\
}$
 & -18 } \\
 \stab{$NH3_{23}$&
$\smat{
 2 & -1 & -2 \\
 -3 & 2 & -1 \\
 -2 & -2 & 2 \\
}$
 & $\frac{1}{24}$ 
$\smat{
 -2 & -6 & -5 \\
 -8 & 0 & -8 \\
 -10 & -6 & -1 \\
}$
 & -24 } \\
 \stab{$NH3_{24}$&
$\smat{
 2 & -1 & -2 \\
 -4 & 2 & -1 \\
 -2 & -2 & 2 \\
}$
 & $\frac{1}{30}$ 
$\smat{
 -2 & -6 & -5 \\
 -10 & 0 & -10 \\
 -12 & -6 & 0 \\
}$
 & -30 } \\
 \stab{$NH3_{25}$&
$\smat{
 2 & -2 & -2 \\
 -2 & 2 & -1 \\
 -2 & -2 & 2 \\
}$
 & $\frac{1}{12}$ 
$\smat{
 -1 & -4 & -3 \\
 -3 & 0 & -3 \\
 -4 & -4 & 0 \\
}$
 & -24 } \\
 \stab{$NH3_{26}$&
$\smat{
 2 & -1 & -2 \\
 -1 & 2 & -1 \\
 -2 & -2 & 2 \\
}$
 & $\frac{1}{12}$ 
$\smat{
 -2 & -6 & -5 \\
 -4 & 0 & -4 \\
 -6 & -6 & -3 \\
}$
 & -12 } \\
 \stab{$H3_{27}$&
$\smat{
 2 & -2 & -2 \\
 -1 & 2 & -1 \\
 -2 & -2 & 2 \\
}$
 & $\frac{1}{8}$ 
$\smat{
 -1 & -4 & -3 \\
 -2 & 0 & -2 \\
 -3 & -4 & -1 \\
}$
 & -16 } \\
 \stab{$NH3_{28}$&
$\smat{
 2 & -3 & -2 \\
 -1 & 2 & -1 \\
 -2 & -2 & 2 \\
}$
 & $\frac{1}{20}$ 
$\smat{
 -2 & -10 & -7 \\
 -4 & 0 & -4 \\
 -6 & -10 & -1 \\
}$
 & -20 } \\
 \stab{$NH3_{29}$&
$\smat{
 2 & -4 & -2 \\
 -1 & 2 & -1 \\
 -2 & -2 & 2 \\
}$
 & $\frac{1}{12}$ 
$\smat{
 -1 & -6 & -4 \\
 -2 & 0 & -2 \\
 -3 & -6 & 0 \\
}$
 & -24 } \\
 \stab{\attn$NH3_{30}$&
$\smat{
 2 & -1 & -1 \\
 -2 & 2 & -1 \\
 -1 & -2 & 2 \\
}$
 & $\frac{1}{7}$ 
$\smat{
 -2 & -4 & -3 \\
 -5 & -3 & -4 \\
 -6 & -5 & -2 \\
}$
 & -7 } \\
 \stab{\attn$NH3_{31}$&
$\smat{
 2 & -1 & -1 \\
 -3 & 2 & -1 \\
 -1 & -2 & 2 \\
}$
 & $\frac{1}{11}$ 
$\smat{
 -2 & -4 & -3 \\
 -7 & -3 & -5 \\
 -8 & -5 & -1 \\
}$
 & -11 } \\
 \stab{$NH3_{32}$&
$\smat{
 2 & -1 & -1 \\
 -4 & 2 & -1 \\
 -1 & -2 & 2 \\
}$
 & $\frac{1}{15}$ 
$\smat{
 -2 & -4 & -3 \\
 -9 & -3 & -6 \\
 -10 & -5 & 0 \\
}$
 & -15 } \\
 \stab{$NH3_{33}$&
$\smat{
 2 & -2 & -1 \\
 -2 & 2 & -1 \\
 -1 & -2 & 2 \\
}$
 & $\frac{1}{12}$ 
$\smat{
 -2 & -6 & -4 \\
 -5 & -3 & -4 \\
 -6 & -6 & 0 \\
}$
 & -12 } \\
 \stab{\attn$NH3_{34}$&
$\smat{
 2 & -1 & -1 \\
 -1 & 2 & -1 \\
 -1 & -2 & 2 \\
}$
 & $\frac{1}{3}$ 
$\smat{
 -2 & -4 & -3 \\
 -3 & -3 & -3 \\
 -4 & -5 & -3 \\
}$
 & -3 } \\
 \stab{\attn$H3_{35}$&
$\smat{
 2 & -2 & -1 \\
 -1 & 2 & -1 \\
 -1 & -2 & 2 \\
}$
 & $\frac{1}{6}$ 
$\smat{
 -2 & -6 & -4 \\
 -3 & -3 & -3 \\
 -4 & -6 & -2 \\
}$
 & -6 } \\
 \stab{\attn$NH3_{36}$&
$\smat{
 2 & -3 & -1 \\
 -1 & 2 & -1 \\
 -1 & -2 & 2 \\
}$
 & $\frac{1}{9}$ 
$\smat{
 -2 & -8 & -5 \\
 -3 & -3 & -3 \\
 -4 & -7 & -1 \\
}$
 & -9 } \\
 \stab{$NH3_{37}$&
$\smat{
 2 & -4 & -1 \\
 -1 & 2 & -1 \\
 -1 & -2 & 2 \\
}$
 & $\frac{1}{12}$ 
$\smat{
 -2 & -10 & -6 \\
 -3 & -3 & -3 \\
 -4 & -8 & 0 \\
}$
 & -12 } \\
 \stab{\attn$NH3_{38}$&
$\smat{
 2 & -1 & -2 \\
 -2 & 2 & -1 \\
 -1 & -2 & 2 \\
}$
 & $\frac{1}{13}$ 
$\smat{
 -2 & -6 & -5 \\
 -5 & -2 & -6 \\
 -6 & -5 & -2 \\
}$
 & -13 } \\
 \stab{\attn$NH3_{39}$&
$\smat{
 2 & -1 & -2 \\
 -3 & 2 & -1 \\
 -1 & -2 & 2 \\
}$
 & $\frac{1}{19}$ 
$\smat{
 -2 & -6 & -5 \\
 -7 & -2 & -8 \\
 -8 & -5 & -1 \\
}$
 & -19 } \\
 \stab{$NH3_{40}$&
$\smat{
 2 & -1 & -2 \\
 -4 & 2 & -1 \\
 -1 & -2 & 2 \\
}$
 & $\frac{1}{25}$ 
$\smat{
 -2 & -6 & -5 \\
 -9 & -2 & -10 \\
 -10 & -5 & 0 \\
}$
 & -25 } \\
 \stab{\attn$NH3_{41}$&
$\smat{
 2 & -1 & -3 \\
 -3 & 2 & -1 \\
 -1 & -2 & 2 \\
}$
 & $\frac{1}{27}$ 
$\smat{
 -2 & -8 & -7 \\
 -7 & -1 & -11 \\
 -8 & -5 & -1 \\
}$
 & -27 } \\
 \stab{$NH3_{42}$&
$\smat{
 2 & -1 & -3 \\
 -4 & 2 & -1 \\
 -1 & -2 & 2 \\
}$
 & $\frac{1}{35}$ 
$\smat{
 -2 & -8 & -7 \\
 -9 & -1 & -14 \\
 -10 & -5 & 0 \\
}$
 & -35 } \\
 \stab{$NH3_{43}$&
$\smat{
 2 & -2 & -3 \\
 -2 & 2 & -1 \\
 -1 & -2 & 2 \\
}$
 & $\frac{1}{24}$ 
$\smat{
 -2 & -10 & -8 \\
 -5 & -1 & -8 \\
 -6 & -6 & 0 \\
}$
 & -24 } \\
 \stab{\attn$NH3_{44}$&
$\smat{
 2 & -1 & -3 \\
 -1 & 2 & -1 \\
 -1 & -2 & 2 \\
}$
 & $\frac{1}{11}$ 
$\smat{
 -2 & -8 & -7 \\
 -3 & -1 & -5 \\
 -4 & -5 & -3 \\
}$
 & -11 } \\
 \stab{\attn$NH3_{45}$&
$\smat{
 2 & -3 & -3 \\
 -1 & 2 & -1 \\
 -1 & -2 & 2 \\
}$
 & $\frac{1}{17}$ 
$\smat{
 -2 & -12 & -9 \\
 -3 & -1 & -5 \\
 -4 & -7 & -1 \\
}$
 & -17 } \\
 \stab{$NH3_{46}$&
$\smat{
 2 & -4 & -3 \\
 -1 & 2 & -1 \\
 -1 & -2 & 2 \\
}$
 & $\frac{1}{20}$ 
$\smat{
 -2 & -14 & -10 \\
 -3 & -1 & -5 \\
 -4 & -8 & 0 \\
}$
 & -20 } \\
 \stab{$NH3_{47}$&
$\smat{
 2 & -1 & -4 \\
 -3 & 2 & -1 \\
 -1 & -2 & 2 \\
}$
 & $\frac{1}{35}$ 
$\smat{
 -2 & -10 & -9 \\
 -7 & 0 & -14 \\
 -8 & -5 & -1 \\
}$
 & -35 } \\
 \stab{$NH3_{48}$&
$\smat{
 2 & -1 & -4 \\
 -4 & 2 & -1 \\
 -1 & -2 & 2 \\
}$
 & $\frac{1}{45}$ 
$\smat{
 -2 & -10 & -9 \\
 -9 & 0 & -18 \\
 -10 & -5 & 0 \\
}$
 & -45 } \\
 \stab{$NH3_{49}$&
$\smat{
 2 & -2 & -4 \\
 -2 & 2 & -1 \\
 -1 & -2 & 2 \\
}$
 & $\frac{1}{30}$ 
$\smat{
 -2 & -12 & -10 \\
 -5 & 0 & -10 \\
 -6 & -6 & 0 \\
}$
 & -30 } \\
 \stab{$NH3_{50}$&
$\smat{
 2 & -1 & -4 \\
 -1 & 2 & -1 \\
 -1 & -2 & 2 \\
}$
 & $\frac{1}{15}$ 
$\smat{
 -2 & -10 & -9 \\
 -3 & 0 & -6 \\
 -4 & -5 & -3 \\
}$
 & -15 } \\
 \stab{$NH3_{51}$&
$\smat{
 2 & -3 & -4 \\
 -1 & 2 & -1 \\
 -1 & -2 & 2 \\
}$
 & $\frac{1}{21}$ 
$\smat{
 -2 & -14 & -11 \\
 -3 & 0 & -6 \\
 -4 & -7 & -1 \\
}$
 & -21 } \\
 \stab{$NH3_{52}$&
$\smat{
 2 & -4 & -4 \\
 -1 & 2 & -1 \\
 -1 & -2 & 2 \\
}$
 & $\frac{1}{24}$ 
$\smat{
 -2 & -16 & -12 \\
 -3 & 0 & -6 \\
 -4 & -8 & 0 \\
}$
 & -24 } \\
 \stab{\attn$NH3_{53}$&
$\smat{
 2 & -1 & -1 \\
 -3 & 2 & -1 \\
 -3 & -3 & 2 \\
}$
 & $\frac{1}{22}$ 
$\smat{
 -1 & -5 & -3 \\
 -9 & -1 & -5 \\
 -15 & -9 & -1 \\
}$
 & -22 } \\
 \stab{$NH3_{54}$&
$\smat{
 2 & -1 & -1 \\
 -4 & 2 & -1 \\
 -3 & -3 & 2 \\
}$
 & $\frac{1}{27}$ 
$\smat{
 -1 & -5 & -3 \\
 -11 & -1 & -6 \\
 -18 & -9 & 0 \\
}$
 & -27 } \\
 \stab{$H3_{55}$&
$\smat{
 2 & -2 & -1 \\
 -2 & 2 & -1 \\
 -3 & -3 & 2 \\
}$
 & $\frac{1}{24}$ 
$\smat{
 -1 & -7 & -4 \\
 -7 & -1 & -4 \\
 -12 & -12 & 0 \\
}$
 & -24 } \\
 \stab{\attn$H3_{56}$&
$\smat{
 2 & -1 & -1 \\
 -1 & 2 & -1 \\
 -3 & -3 & 2 \\
}$
 & $\frac{1}{12}$ 
$\smat{
 -1 & -5 & -3 \\
 -5 & -1 & -3 \\
 -9 & -9 & -3 \\
}$
 & -12 } \\
 \stab{$NH3_{57}$&
$\smat{
 2 & -1 & -1 \\
 -3 & 2 & -1 \\
 -4 & -3 & 2 \\
}$
 & $\frac{1}{25}$ 
$\smat{
 -1 & -5 & -3 \\
 -10 & 0 & -5 \\
 -17 & -10 & -1 \\
}$
 & -25 } \\
 \stab{$NH3_{58}$&
$\smat{
 2 & -1 & -1 \\
 -4 & 2 & -1 \\
 -4 & -3 & 2 \\
}$
 & $\frac{1}{30}$ 
$\smat{
 -1 & -5 & -3 \\
 -12 & 0 & -6 \\
 -20 & -10 & 0 \\
}$
 & -30 } \\
 \stab{$NH3_{59}$&
$\smat{
 2 & -2 & -1 \\
 -2 & 2 & -1 \\
 -4 & -3 & 2 \\
}$
 & $\frac{1}{28}$ 
$\smat{
 -1 & -7 & -4 \\
 -8 & 0 & -4 \\
 -14 & -14 & 0 \\
}$
 & -28 } \\
 \stab{$NH3_{60}$&
$\smat{
 2 & -1 & -1 \\
 -1 & 2 & -1 \\
 -4 & -3 & 2 \\
}$
 & $\frac{1}{15}$ 
$\smat{
 -1 & -5 & -3 \\
 -6 & 0 & -3 \\
 -11 & -10 & -3 \\
}$
 & -15 } \\
 \stab{$NH3_{61}$&
$\smat{
 2 & -3 & -1 \\
 -1 & 2 & -1 \\
 -4 & -3 & 2 \\
}$
 & $\frac{1}{27}$ 
$\smat{
 -1 & -9 & -5 \\
 -6 & 0 & -3 \\
 -11 & -18 & -1 \\
}$
 & -27 } \\
 \stab{$NH3_{62}$&
$\smat{
 2 & -4 & -1 \\
 -1 & 2 & -1 \\
 -4 & -3 & 2 \\
}$
 & $\frac{1}{33}$ 
$\smat{
 -1 & -11 & -6 \\
 -6 & 0 & -3 \\
 -11 & -22 & 0 \\
}$
 & -33 } \\
 \stab{$NH3_{63}$&
$\smat{
 2 & -1 & -2 \\
 -3 & 2 & -1 \\
 -2 & -3 & 2 \\
}$
 & $\frac{1}{32}$ 
$\smat{
 -1 & -8 & -5 \\
 -8 & 0 & -8 \\
 -13 & -8 & -1 \\
}$
 & -32 } \\
 \stab{$NH3_{64}$&
$\smat{
 2 & -1 & -2 \\
 -4 & 2 & -1 \\
 -2 & -3 & 2 \\
}$
 & $\frac{1}{40}$ 
$\smat{
 -1 & -8 & -5 \\
 -10 & 0 & -10 \\
 -16 & -8 & 0 \\
}$
 & -40 } \\
 \stab{$NH3_{65}$&
$\smat{
 2 & -2 & -2 \\
 -2 & 2 & -1 \\
 -2 & -3 & 2 \\
}$
 & $\frac{1}{30}$ 
$\smat{
 -1 & -10 & -6 \\
 -6 & 0 & -6 \\
 -10 & -10 & 0 \\
}$
 & -30 } \\
 \stab{$NH3_{66}$&
$\smat{
 2 & -1 & -2 \\
 -1 & 2 & -1 \\
 -2 & -3 & 2 \\
}$
 & $\frac{1}{16}$ 
$\smat{
 -1 & -8 & -5 \\
 -4 & 0 & -4 \\
 -7 & -8 & -3 \\
}$
 & -16 } \\
 \stab{$H3_{67}$&
$\smat{
 2 & -3 & -2 \\
 -1 & 2 & -1 \\
 -2 & -3 & 2 \\
}$
 & $\frac{1}{24}$ 
$\smat{
 -1 & -12 & -7 \\
 -4 & 0 & -4 \\
 -7 & -12 & -1 \\
}$
 & -24 } \\
 \stab{$NH3_{68}$&
$\smat{
 2 & -4 & -2 \\
 -1 & 2 & -1 \\
 -2 & -3 & 2 \\
}$
 & $\frac{1}{28}$ 
$\smat{
 -1 & -14 & -8 \\
 -4 & 0 & -4 \\
 -7 & -14 & 0 \\
}$
 & -28 } \\
 \stab{\attn$NH3_{69}$&
$\smat{
 2 & -1 & -1 \\
 -3 & 2 & -1 \\
 -1 & -3 & 2 \\
}$
 & $\frac{1}{16}$ 
$\smat{
 -1 & -5 & -3 \\
 -7 & -3 & -5 \\
 -11 & -7 & -1 \\
}$
 & -16 } \\
 \stab{$NH3_{70}$&
$\smat{
 2 & -1 & -1 \\
 -4 & 2 & -1 \\
 -1 & -3 & 2 \\
}$
 & $\frac{1}{21}$ 
$\smat{
 -1 & -5 & -3 \\
 -9 & -3 & -6 \\
 -14 & -7 & 0 \\
}$
 & -21 } \\
 \stab{$NH3_{71}$&
$\smat{
 2 & -2 & -1 \\
 -2 & 2 & -1 \\
 -1 & -3 & 2 \\
}$
 & $\frac{1}{16}$ 
$\smat{
 -1 & -7 & -4 \\
 -5 & -3 & -4 \\
 -8 & -8 & 0 \\
}$
 & -16 } \\
 \stab{\attn$NH3_{72}$&
$\smat{
 2 & -1 & -1 \\
 -1 & 2 & -1 \\
 -1 & -3 & 2 \\
}$
 & $\frac{1}{6}$ 
$\smat{
 -1 & -5 & -3 \\
 -3 & -3 & -3 \\
 -5 & -7 & -3 \\
}$
 & -6 } \\
 \stab{\attn$H3_{73}$&
$\smat{
 2 & -3 & -1 \\
 -1 & 2 & -1 \\
 -1 & -3 & 2 \\
}$
 & $\frac{1}{12}$ 
$\smat{
 -1 & -9 & -5 \\
 -3 & -3 & -3 \\
 -5 & -9 & -1 \\
}$
 & -12 } \\
 \stab{$NH3_{74}$&
$\smat{
 2 & -4 & -1 \\
 -1 & 2 & -1 \\
 -1 & -3 & 2 \\
}$
 & $\frac{1}{15}$ 
$\smat{
 -1 & -11 & -6 \\
 -3 & -3 & -3 \\
 -5 & -10 & 0 \\
}$
 & -15 } \\
 \stab{\attn$NH3_{75}$&
$\smat{
 2 & -1 & -3 \\
 -3 & 2 & -1 \\
 -1 & -3 & 2 \\
}$
 & $\frac{1}{38}$ 
$\smat{
 -1 & -11 & -7 \\
 -7 & -1 & -11 \\
 -11 & -7 & -1 \\
}$
 & -38 } \\
 \stab{$NH3_{76}$&
$\smat{
 2 & -1 & -3 \\
 -4 & 2 & -1 \\
 -1 & -3 & 2 \\
}$
 & $\frac{1}{49}$ 
$\smat{
 -1 & -11 & -7 \\
 -9 & -1 & -14 \\
 -14 & -7 & 0 \\
}$
 & -49 } \\
 \stab{$NH3_{77}$&
$\smat{
 2 & -1 & -4 \\
 -4 & 2 & -1 \\
 -1 & -3 & 2 \\
}$
 & $\frac{1}{63}$ 
$\smat{
 -1 & -14 & -9 \\
 -9 & 0 & -18 \\
 -14 & -7 & 0 \\
}$
 & -63 } \\
 \stab{$NH3_{78}$&
$\smat{
 2 & -2 & -4 \\
 -2 & 2 & -1 \\
 -1 & -3 & 2 \\
}$
 & $\frac{1}{40}$ 
$\smat{
 -1 & -16 & -10 \\
 -5 & 0 & -10 \\
 -8 & -8 & 0 \\
}$
 & -40 } \\
 \stab{$NH3_{79}$&
$\smat{
 2 & -1 & -4 \\
 -1 & 2 & -1 \\
 -1 & -3 & 2 \\
}$
 & $\frac{1}{21}$ 
$\smat{
 -1 & -14 & -9 \\
 -3 & 0 & -6 \\
 -5 & -7 & -3 \\
}$
 & -21 } \\
 \stab{$NH3_{80}$&
$\smat{
 2 & -4 & -4 \\
 -1 & 2 & -1 \\
 -1 & -3 & 2 \\
}$
 & $\frac{1}{30}$ 
$\smat{
 -1 & -20 & -12 \\
 -3 & 0 & -6 \\
 -5 & -10 & 0 \\
}$
 & -30 } \\
 \stab{$NH3_{81}$&
$\smat{
 2 & -1 & -1 \\
 -4 & 2 & -1 \\
 -4 & -4 & 2 \\
}$
 & $\frac{1}{12}$ 
$\smat{
 0 & -2 & -1 \\
 -4 & 0 & -2 \\
 -8 & -4 & 0 \\
}$
 & -36 } \\
 \stab{$H3_{82}$&
$\smat{
 2 & -2 & -1 \\
 -2 & 2 & -1 \\
 -4 & -4 & 2 \\
}$
 & $\frac{1}{8}$ 
$\smat{
 0 & -2 & -1 \\
 -2 & 0 & -1 \\
 -4 & -4 & 0 \\
}$
 & -32 } \\
 \stab{$H3_{83}$&
$\smat{
 2 & -1 & -1 \\
 -1 & 2 & -1 \\
 -4 & -4 & 2 \\
}$
 & $\frac{1}{6}$ 
$\smat{
 0 & -2 & -1 \\
 -2 & 0 & -1 \\
 -4 & -4 & -1 \\
}$
 & -18 } \\
 \stab{$NH3_{84}$&
$\smat{
 2 & -1 & -2 \\
 -4 & 2 & -1 \\
 -2 & -4 & 2 \\
}$
 & $\frac{1}{10}$ 
$\smat{
 0 & -2 & -1 \\
 -2 & 0 & -2 \\
 -4 & -2 & 0 \\
}$
 & -50 } \\
 \stab{$NH3_{85}$&
$\smat{
 2 & -2 & -2 \\
 -2 & 2 & -1 \\
 -2 & -4 & 2 \\
}$
 & $\frac{1}{6}$ 
$\smat{
 0 & -2 & -1 \\
 -1 & 0 & -1 \\
 -2 & -2 & 0 \\
}$
 & -36 } \\
 \stab{$NH3_{86}$&
$\smat{
 2 & -1 & -2 \\
 -1 & 2 & -1 \\
 -2 & -4 & 2 \\
}$
 & $\frac{1}{20}$ 
$\smat{
 0 & -10 & -5 \\
 -4 & 0 & -4 \\
 -8 & -10 & -3 \\
}$
 & -20 } \\
 \stab{$H3_{87}$&
$\smat{
 2 & -4 & -2 \\
 -1 & 2 & -1 \\
 -2 & -4 & 2 \\
}$
 & $\frac{1}{8}$ 
$\smat{
 0 & -4 & -2 \\
 -1 & 0 & -1 \\
 -2 & -4 & 0 \\
}$
 & -32 } \\
 \stab{$NH3_{88}$&
$\smat{
 2 & -1 & -1 \\
 -4 & 2 & -1 \\
 -1 & -4 & 2 \\
}$
 & $\frac{1}{9}$ 
$\smat{
 0 & -2 & -1 \\
 -3 & -1 & -2 \\
 -6 & -3 & 0 \\
}$
 & -27 } \\
 \stab{$NH3_{89}$&
$\smat{
 2 & -2 & -1 \\
 -2 & 2 & -1 \\
 -1 & -4 & 2 \\
}$
 & $\frac{1}{20}$ 
$\smat{
 0 & -8 & -4 \\
 -5 & -3 & -4 \\
 -10 & -10 & 0 \\
}$
 & -20 } \\
 \stab{$NH3_{90}$&
$\smat{
 2 & -1 & -1 \\
 -1 & 2 & -1 \\
 -1 & -4 & 2 \\
}$
 & $\frac{1}{3}$ 
$\smat{
 0 & -2 & -1 \\
 -1 & -1 & -1 \\
 -2 & -3 & -1 \\
}$
 & -9 } \\
 \stab{$H3_{91}$&
$\smat{
 2 & -4 & -1 \\
 -1 & 2 & -1 \\
 -1 & -4 & 2 \\
}$
 & $\frac{1}{6}$ 
$\smat{
 0 & -4 & -2 \\
 -1 & -1 & -1 \\
 -2 & -4 & 0 \\
}$
 & -18 } \\
 \stab{$NH3_{92}$&
$\smat{
 2 & -1 & -4 \\
 -4 & 2 & -1 \\
 -1 & -4 & 2 \\
}$
 & $\frac{1}{9}$ 
$\smat{
 0 & -2 & -1 \\
 -1 & 0 & -2 \\
 -2 & -1 & 0 \\
}$
 & -81 } \\
 \stab{$H3_{93}$&
$\smat{
 2 & -2 & -2 \\
 -2 & 2 & -2 \\
 -2 & -2 & 2 \\
}$
 & $\frac{1}{4}$ 
$\smat{
 0 & -1 & -1 \\
 -1 & 0 & -1 \\
 -1 & -1 & 0 \\
}$
 & -32 } \\
 \stab{$H3_{94}$&
$\smat{
 2 & -1 & -2 \\
 -1 & 2 & -2 \\
 -2 & -2 & 2 \\
}$
 & $\frac{1}{6}$ 
$\smat{
 0 & -2 & -2 \\
 -2 & 0 & -2 \\
 -2 & -2 & -1 \\
}$
 & -18 } \\
 \stab{$H3_{95}$&
$\smat{
 2 & -1 & -1 \\
 -1 & 2 & -2 \\
 -1 & -2 & 2 \\
}$
 & $\frac{1}{8}$ 
$\smat{
 0 & -4 & -4 \\
 -4 & -3 & -5 \\
 -4 & -5 & -3 \\
}$
 & -8 } \\
 \stab{\attn$H3_{96}$&
$\smat{
 2 & -1 & 0 \\
 -3 & 2 & -2 \\
 0 & -1 & 2 \\
}$
 & $\frac{1}{2}$ 
$\smat{
 -2 & -2 & -2 \\
 -6 & -4 & -4 \\
 -3 & -2 & -1 \\
}$
 & -2 } \\
 \stab{$H3_{97}$&
$\smat{
 2 & -1 & 0 \\
 -4 & 2 & -2 \\
 0 & -1 & 2 \\
}$
 & $\frac{1}{2}$ 
$\smat{
 -1 & -1 & -1 \\
 -4 & -2 & -2 \\
 -2 & -1 & 0 \\
}$
 & -4 } \\
 \stab{$H3_{98}$&
$\smat{
 2 & -2 & 0 \\
 -2 & 2 & -2 \\
 0 & -1 & 2 \\
}$
 & $\frac{1}{2}$ 
$\smat{
 -1 & -2 & -2 \\
 -2 & -2 & -2 \\
 -1 & -1 & 0 \\
}$
 & -4 } \\
 \stab{\attn$H3_{99}$&
$\smat{
 2 & -3 & 0 \\
 -1 & 2 & -2 \\
 0 & -1 & 2 \\
}$
 & $\frac{1}{2}$ 
$\smat{
 -2 & -6 & -6 \\
 -2 & -4 & -4 \\
 -1 & -2 & -1 \\
}$
 & -2 } \\
 \stab{$H3_{100}$&
$\smat{
 2 & -4 & 0 \\
 -1 & 2 & -2 \\
 0 & -1 & 2 \\
}$
 & $\frac{1}{4}$ 
$\smat{
 -2 & -8 & -8 \\
 -2 & -4 & -4 \\
 -1 & -2 & 0 \\
}$
 & -4 } \\
 \stab{\attn$H3_{101}$&
$\smat{
 2 & -1 & 0 \\
 -3 & 2 & -3 \\
 0 & -1 & 2 \\
}$
 & $\frac{1}{4}$ 
$\smat{
 -1 & -2 & -3 \\
 -6 & -4 & -6 \\
 -3 & -2 & -1 \\
}$
 & -4 } \\
 \stab{$H3_{102}$&
$\smat{
 2 & -1 & 0 \\
 -4 & 2 & -3 \\
 0 & -1 & 2 \\
}$
 & $\frac{1}{6}$ 
$\smat{
 -1 & -2 & -3 \\
 -8 & -4 & -6 \\
 -4 & -2 & 0 \\
}$
 & -6 } \\
 \stab{$H3_{103}$&
$\smat{
 2 & -2 & 0 \\
 -2 & 2 & -3 \\
 0 & -1 & 2 \\
}$
 & $\frac{1}{6}$ 
$\smat{
 -1 & -4 & -6 \\
 -4 & -4 & -6 \\
 -2 & -2 & 0 \\
}$
 & -6 } \\
 \stab{\attn$H3_{104}$&
$\smat{
 2 & -2 & 0 \\
 -1 & 2 & -3 \\
 0 & -1 & 2 \\
}$
 & $\frac{1}{2}$ 
$\smat{
 -1 & -4 & -6 \\
 -2 & -4 & -6 \\
 -1 & -2 & -2 \\
}$
 & -2 } \\
 \stab{\attn$H3_{105}$&
$\smat{
 2 & -3 & 0 \\
 -1 & 2 & -3 \\
 0 & -1 & 2 \\
}$
 & $\frac{1}{4}$ 
$\smat{
 -1 & -6 & -9 \\
 -2 & -4 & -6 \\
 -1 & -2 & -1 \\
}$
 & -4 } \\
 \stab{$H3_{106}$&
$\smat{
 2 & -4 & 0 \\
 -1 & 2 & -3 \\
 0 & -1 & 2 \\
}$
 & $\frac{1}{6}$ 
$\smat{
 -1 & -8 & -12 \\
 -2 & -4 & -6 \\
 -1 & -2 & 0 \\
}$
 & -6 } \\
 \stab{$H3_{107}$&
$\smat{
 2 & -1 & 0 \\
 -4 & 2 & -4 \\
 0 & -1 & 2 \\
}$
 & $\frac{1}{4}$ 
$\smat{
 0 & -1 & -2 \\
 -4 & -2 & -4 \\
 -2 & -1 & 0 \\
}$
 & -8 } \\
 \stab{$H3_{108}$&
$\smat{
 2 & -2 & 0 \\
 -2 & 2 & -4 \\
 0 & -1 & 2 \\
}$
 & $\frac{1}{4}$ 
$\smat{
 0 & -2 & -4 \\
 -2 & -2 & -4 \\
 -1 & -1 & 0 \\
}$
 & -8 } \\
 \stab{$H3_{109}$&
$\smat{
 2 & -1 & 0 \\
 -1 & 2 & -4 \\
 0 & -1 & 2 \\
}$
 & $\frac{1}{2}$ 
$\smat{
 0 & -2 & -4 \\
 -2 & -4 & -8 \\
 -1 & -2 & -3 \\
}$
 & -2 } \\
 \stab{$H3_{110}$&
$\smat{
 2 & -2 & 0 \\
 -1 & 2 & -4 \\
 0 & -1 & 2 \\
}$
 & $\frac{1}{4}$ 
$\smat{
 0 & -4 & -8 \\
 -2 & -4 & -8 \\
 -1 & -2 & -2 \\
}$
 & -4 } \\
 \stab{$H3_{111}$&
$\smat{
 2 & -3 & 0 \\
 -1 & 2 & -4 \\
 0 & -1 & 2 \\
}$
 & $\frac{1}{6}$ 
$\smat{
 0 & -6 & -12 \\
 -2 & -4 & -8 \\
 -1 & -2 & -1 \\
}$
 & -6 } \\
 \stab{$H3_{112}$&
$\smat{
 2 & -4 & 0 \\
 -1 & 2 & -4 \\
 0 & -1 & 2 \\
}$
 & $\frac{1}{8}$ 
$\smat{
 0 & -8 & -16 \\
 -2 & -4 & -8 \\
 -1 & -2 & 0 \\
}$
 & -8 } \\
 \stab{$H3_{113}$&
$\smat{
 2 & -2 & 0 \\
 -2 & 2 & -2 \\
 0 & -2 & 2 \\
}$
 & $\frac{1}{2}$ 
$\smat{
 0 & -1 & -1 \\
 -1 & -1 & -1 \\
 -1 & -1 & 0 \\
}$
 & -8 } \\
 \stab{$H3_{114}$&
$\smat{
 2 & -1 & 0 \\
 -1 & 2 & -2 \\
 0 & -2 & 2 \\
}$
 & $\frac{1}{2}$ 
$\smat{
 0 & -2 & -2 \\
 -2 & -4 & -4 \\
 -2 & -4 & -3 \\
}$
 & -2 } \\
 \stab{$H3_{115}$&
$\smat{
 2 & -2 & 0 \\
 -1 & 2 & -2 \\
 0 & -2 & 2 \\
}$
 & $\frac{1}{2}$ 
$\smat{
 0 & -2 & -2 \\
 -1 & -2 & -2 \\
 -1 & -2 & -1 \\
}$
 & -4 } \\
 \stab{$H3_{116}$&
$\smat{
 2 & -3 & 0 \\
 -1 & 2 & -2 \\
 0 & -2 & 2 \\
}$
 & $\frac{1}{6}$ 
$\smat{
 0 & -6 & -6 \\
 -2 & -4 & -4 \\
 -2 & -4 & -1 \\
}$
 & -6 } \\
 \stab{$H3_{117}$&
$\smat{
 2 & -4 & 0 \\
 -1 & 2 & -2 \\
 0 & -2 & 2 \\
}$
 & $\frac{1}{4}$ 
$\smat{
 0 & -4 & -4 \\
 -1 & -2 & -2 \\
 -1 & -2 & 0 \\
}$
 & -8 } \\
 \stab{$H3_{118}$&
$\smat{
 2 & -4 & 0 \\
 -1 & 2 & -1 \\
 0 & -1 & 2 \\
}$
 & $\frac{1}{2}$ 
$\smat{
 -3 & -8 & -4 \\
 -2 & -4 & -2 \\
 -1 & -2 & 0 \\
}$
 & -2 } \\
 \stab{\attn$H3_{119}$&
$\smat{
 2 & -3 & 0 \\
 -1 & 2 & -1 \\
 0 & -2 & 2 \\
}$
 & $\frac{1}{2}$ 
$\smat{
 -2 & -6 & -3 \\
 -2 & -4 & -2 \\
 -2 & -4 & -1 \\
}$
 & -2 } \\
 \stab{$H3_{120}$&
$\smat{
 2 & -4 & 0 \\
 -1 & 2 & -1 \\
 0 & -2 & 2 \\
}$
 & $\frac{1}{2}$ 
$\smat{
 -1 & -4 & -2 \\
 -1 & -2 & -1 \\
 -1 & -2 & 0 \\
}$
 & -4 } \\
 \stab{\attn$H3_{121}$&
$\smat{
 2 & -3 & 0 \\
 -1 & 2 & -1 \\
 0 & -3 & 2 \\
}$
 & $\frac{1}{4}$ 
$\smat{
 -1 & -6 & -3 \\
 -2 & -4 & -2 \\
 -3 & -6 & -1 \\
}$
 & -4 } \\
 \stab{$H3_{122}$&
$\smat{
 2 & -4 & 0 \\
 -1 & 2 & -1 \\
 0 & -3 & 2 \\
}$
 & $\frac{1}{6}$ 
$\smat{
 -1 & -8 & -4 \\
 -2 & -4 & -2 \\
 -3 & -6 & 0 \\
}$
 & -6 } \\
 \stab{$H3_{123}$&
$\smat{
 2 & -4 & 0 \\
 -1 & 2 & -1 \\
 0 & -4 & 2 \\
}$
 & $\frac{1}{4}$ 
$\smat{
 0 & -4 & -2 \\
 -1 & -2 & -1 \\
 -2 & -4 & 0 \\
}$
 & -8 } \\
 \stab{$H4_{1}$&
$\smat{
 2 & -1 & 0 & 0 \\
 -2 & 2 & -1 & 0 \\
 0 & -2 & 2 & -2 \\
 0 & 0 & -1 & 2 \\
}$
 & $\frac{1}{2}$ 
$\smat{
 0 & -1 & -1 & -1 \\
 -2 & -2 & -2 & -2 \\
 -4 & -4 & -2 & -2 \\
 -2 & -2 & -1 & 0 \\
}$
 & -4 } \\
 \stab{$H4_{2}$&
$\smat{
 2 & -1 & 0 & 0 \\
 -3 & 2 & -1 & 0 \\
 0 & -1 & 2 & -2 \\
 0 & 0 & -1 & 2 \\
}$
 & $\frac{1}{2}$ 
$\smat{
 -2 & -2 & -2 & -2 \\
 -6 & -4 & -4 & -4 \\
 -6 & -4 & -2 & -2 \\
 -3 & -2 & -1 & 0 \\
}$
 & -2 } \\
 \stab{$H4_{3}$&
$\smat{
 2 & -1 & 0 & 0 \\
 -1 & 2 & -1 & 0 \\
 0 & -2 & 2 & -2 \\
 0 & 0 & -1 & 2 \\
}$
 & $\frac{1}{2}$ 
$\smat{
 0 & -2 & -2 & -2 \\
 -2 & -4 & -4 & -4 \\
 -4 & -8 & -6 & -6 \\
 -2 & -4 & -3 & -2 \\
}$
 & -2 } \\
 \stab{$H4_{4}$&
$\smat{
 2 & -1 & 0 & 0 \\
 -1 & 2 & -2 & 0 \\
 0 & -1 & 2 & -2 \\
 0 & 0 & -1 & 2 \\
}$
 & $\frac{1}{2}$ 
$\smat{
 0 & -2 & -4 & -4 \\
 -2 & -4 & -8 & -8 \\
 -2 & -4 & -6 & -6 \\
 -1 & -2 & -3 & -2 \\
}$
 & -2 } \\
 \stab{$H4_{5}$&
$\smat{
 2 & -2 & 0 & 0 \\
 -1 & 2 & -1 & 0 \\
 0 & -2 & 2 & -2 \\
 0 & 0 & -1 & 2 \\
}$
 & $\frac{1}{2}$ 
$\smat{
 0 & -2 & -2 & -2 \\
 -1 & -2 & -2 & -2 \\
 -2 & -4 & -2 & -2 \\
 -1 & -2 & -1 & 0 \\
}$
 & -4 } \\
 \stab{$H4_{6}$&
$\smat{
 2 & -2 & 0 & 0 \\
 -1 & 2 & -2 & 0 \\
 0 & -1 & 2 & -2 \\
 0 & 0 & -1 & 2 \\
}$
 & $\frac{1}{4}$ 
$\smat{
 0 & -4 & -8 & -8 \\
 -2 & -4 & -8 & -8 \\
 -2 & -4 & -4 & -4 \\
 -1 & -2 & -2 & 0 \\
}$
 & -4 } \\
 \stab{$H4_{7}$&
$\smat{
 2 & -3 & 0 & 0 \\
 -1 & 2 & -1 & 0 \\
 0 & -1 & 2 & -2 \\
 0 & 0 & -1 & 2 \\
}$
 & $\frac{1}{2}$ 
$\smat{
 -2 & -6 & -6 & -6 \\
 -2 & -4 & -4 & -4 \\
 -2 & -4 & -2 & -2 \\
 -1 & -2 & -1 & 0 \\
}$
 & -2 } \\
 \stab{$H4_{8}$&
$\smat{
 2 & -1 & 0 & 0 \\
 -3 & 2 & -1 & 0 \\
 0 & -1 & 2 & -3 \\
 0 & 0 & -1 & 2 \\
}$
 & $\frac{1}{3}$ 
$\smat{
 0 & -1 & -2 & -3 \\
 -3 & -2 & -4 & -6 \\
 -6 & -4 & -2 & -3 \\
 -3 & -2 & -1 & 0 \\
}$
 & -3 } \\
 \stab{$H4_{9}$&
$\smat{
 2 & -1 & 0 & 0 \\
 -1 & 2 & -1 & 0 \\
 0 & -1 & 2 & -3 \\
 0 & 0 & -1 & 2 \\
}$
 & 
$\smat{
 0 & -1 & -2 & -3 \\
 -1 & -2 & -4 & -6 \\
 -2 & -4 & -6 & -9 \\
 -1 & -2 & -3 & -4 \\
}$
 & -1 } \\
 \stab{$H4_{10}$&
$\smat{
 2 & -2 & 0 & 0 \\
 -1 & 2 & -1 & 0 \\
 0 & -1 & 2 & -3 \\
 0 & 0 & -1 & 2 \\
}$
 & $\frac{1}{2}$ 
$\smat{
 0 & -2 & -4 & -6 \\
 -1 & -2 & -4 & -6 \\
 -2 & -4 & -4 & -6 \\
 -1 & -2 & -2 & -2 \\
}$
 & -2 } \\
 \stab{$H4_{11}$&
$\smat{
 2 & -3 & 0 & 0 \\
 -1 & 2 & -1 & 0 \\
 0 & -1 & 2 & -3 \\
 0 & 0 & -1 & 2 \\
}$
 & $\frac{1}{3}$ 
$\smat{
 0 & -3 & -6 & -9 \\
 -1 & -2 & -4 & -6 \\
 -2 & -4 & -2 & -3 \\
 -1 & -2 & -1 & 0 \\
}$
 & -3 } \\
 \stab{$H4_{12}$&
$\smat{
 2 & -1 & 0 & 0 \\
 -1 & 2 & -3 & 0 \\
 0 & -1 & 2 & -1 \\
 0 & 0 & -1 & 2 \\
}$
 & $\frac{1}{3}$ 
$\smat{
 0 & -3 & -6 & -3 \\
 -3 & -6 & -12 & -6 \\
 -2 & -4 & -6 & -3 \\
 -1 & -2 & -3 & 0 \\
}$
 & -3 } \\
 \stab{$H4_{13}$&
$\smat{
 2 & -2 & 0 & 0 \\
 -1 & 2 & -2 & 0 \\
 0 & -1 & 2 & -1 \\
 0 & 0 & -1 & 2 \\
}$
 & $\frac{1}{2}$ 
$\smat{
 -2 & -6 & -8 & -4 \\
 -3 & -6 & -8 & -4 \\
 -2 & -4 & -4 & -2 \\
 -1 & -2 & -2 & 0 \\
}$
 & -2 } \\
 \stab{$H4_{14}$&
$\smat{
 2 & -2 & 0 & 0 \\
 -1 & 2 & -1 & 0 \\
 0 & -2 & 2 & -1 \\
 0 & 0 & -1 & 2 \\
}$
 & $\frac{1}{2}$ 
$\smat{
 -2 & -6 & -4 & -2 \\
 -3 & -6 & -4 & -2 \\
 -4 & -8 & -4 & -2 \\
 -2 & -4 & -2 & 0 \\
}$
 & -2 } \\
 \stab{$H4_{15}$&
$\smat{
 2 & -3 & 0 & 0 \\
 -1 & 2 & -1 & 0 \\
 0 & -1 & 2 & -1 \\
 0 & 0 & -1 & 2 \\
}$
 & 
$\smat{
 -4 & -9 & -6 & -3 \\
 -3 & -6 & -4 & -2 \\
 -2 & -4 & -2 & -1 \\
 -1 & -2 & -1 & 0 \\
}$
 & -1 } \\
 \stab{$H4_{16}$&
$\smat{
 2 & -2 & 0 & 0 \\
 -1 & 2 & -2 & 0 \\
 0 & -1 & 2 & -1 \\
 0 & 0 & -2 & 2 \\
}$
 & $\frac{1}{2}$ 
$\smat{
 0 & -2 & -4 & -2 \\
 -1 & -2 & -4 & -2 \\
 -1 & -2 & -2 & -1 \\
 -1 & -2 & -2 & 0 \\
}$
 & -4 } \\
 \stab{$H4_{17}$&
$\smat{
 2 & -3 & 0 & 0 \\
 -1 & 2 & -1 & 0 \\
 0 & -1 & 2 & -1 \\
 0 & 0 & -2 & 2 \\
}$
 & $\frac{1}{2}$ 
$\smat{
 -2 & -6 & -6 & -3 \\
 -2 & -4 & -4 & -2 \\
 -2 & -4 & -2 & -1 \\
 -2 & -4 & -2 & 0 \\
}$
 & -2 } \\
 \stab{$H4_{18}$&
$\smat{
 2 & -3 & 0 & 0 \\
 -1 & 2 & -1 & 0 \\
 0 & -1 & 2 & -1 \\
 0 & 0 & -3 & 2 \\
}$
 & $\frac{1}{3}$ 
$\smat{
 0 & -3 & -6 & -3 \\
 -1 & -2 & -4 & -2 \\
 -2 & -4 & -2 & -1 \\
 -3 & -6 & -3 & 0 \\
}$
 & -3 } \\
 \stab{$H4_{19}$&
$\smat{
 2 & -1 & 0 & -1 \\
 -1 & 2 & -1 & 0 \\
 0 & -2 & 2 & -2 \\
 -1 & 0 & -1 & 2 \\
}$
 & $\frac{1}{4}$ 
$\smat{
 0 & -2 & -2 & -2 \\
 -2 & -1 & -2 & -3 \\
 -4 & -4 & -2 & -4 \\
 -2 & -3 & -2 & -1 \\
}$
 & -8 } \\
 \stab{$NH4_{20}$&
$\smat{
 2 & -1 & 0 & -1 \\
 -2 & 2 & -1 & 0 \\
 0 & -2 & 2 & -2 \\
 -1 & 0 & -1 & 2 \\
}$
 & $\frac{1}{12}$ 
$\smat{
 0 & -4 & -4 & -4 \\
 -6 & -2 & -5 & -8 \\
 -12 & -8 & -2 & -8 \\
 -6 & -6 & -3 & 0 \\
}$
 & -12 } \\
 \stab{$NH4_{21}$&
$\smat{
 2 & -2 & 0 & -1 \\
 -1 & 2 & -1 & 0 \\
 0 & -2 & 2 & -2 \\
 -1 & 0 & -1 & 2 \\
}$
 & $\frac{1}{6}$ 
$\smat{
 0 & -3 & -3 & -3 \\
 -2 & -1 & -2 & -3 \\
 -4 & -5 & -1 & -3 \\
 -2 & -4 & -2 & 0 \\
}$
 & -12 } \\
 \stab{$H4_{22}$&
$\smat{
 2 & -1 & 0 & -1 \\
 -2 & 2 & -1 & 0 \\
 0 & -2 & 2 & -2 \\
 -2 & 0 & -1 & 2 \\
}$
 & $\frac{1}{4}$ 
$\smat{
 0 & -1 & -1 & -1 \\
 -2 & 0 & -1 & -2 \\
 -4 & -2 & 0 & -2 \\
 -2 & -2 & -1 & 0 \\
}$
 & -16 } \\
 \stab{$NH4_{23}$&
$\smat{
 2 & -2 & 0 & -1 \\
 -1 & 2 & -1 & 0 \\
 0 & -2 & 2 & -2 \\
 -2 & 0 & -1 & 2 \\
}$
 & $\frac{1}{6}$ 
$\smat{
 0 & -2 & -2 & -2 \\
 -2 & 0 & -1 & -2 \\
 -4 & -4 & 0 & -2 \\
 -2 & -4 & -2 & 0 \\
}$
 & -18 } \\
 \stab{$H4_{24}$&
$\smat{
 2 & -2 & 0 & -2 \\
 -1 & 2 & -1 & 0 \\
 0 & -2 & 2 & -2 \\
 -1 & 0 & -1 & 2 \\
}$
 & $\frac{1}{4}$ 
$\smat{
 0 & -2 & -2 & -2 \\
 -1 & 0 & -1 & -2 \\
 -2 & -2 & 0 & -2 \\
 -1 & -2 & -1 & 0 \\
}$
 & -16 } \\
 \stab{$H4_{25}$&
$\smat{
 2 & -1 & 0 & -1 \\
 -1 & 2 & -1 & 0 \\
 0 & -3 & 2 & -1 \\
 -3 & 0 & -1 & 2 \\
}$
 & $\frac{1}{6}$ 
$\smat{
 0 & -3 & -2 & -1 \\
 -3 & 0 & -1 & -2 \\
 -6 & -3 & 0 & -3 \\
 -3 & -6 & -3 & 0 \\
}$
 & -12 } \\
 \stab{$NH4_{26}$&
$\smat{
 2 & -1 & 0 & -1 \\
 -1 & 2 & -1 & 0 \\
 0 & -3 & 2 & -1 \\
 -1 & 0 & -1 & 2 \\
}$
 & $\frac{1}{4}$ 
$\smat{
 0 & -3 & -2 & -1 \\
 -2 & -2 & -2 & -2 \\
 -4 & -5 & -2 & -3 \\
 -2 & -4 & -2 & 0 \\
}$
 & -8 } \\
 \stab{$NH4_{27}$&
$\smat{
 2 & -1 & 0 & -3 \\
 -1 & 2 & -1 & 0 \\
 0 & -3 & 2 & -1 \\
 -1 & 0 & -1 & 2 \\
}$
 & $\frac{1}{4}$ 
$\smat{
 0 & -3 & -2 & -1 \\
 -1 & 0 & -1 & -2 \\
 -2 & -1 & 0 & -3 \\
 -1 & -2 & -1 & 0 \\
}$
 & -16 } \\
 \stab{$NH4_{28}$&
$\smat{
 2 & -1 & 0 & -2 \\
 -1 & 2 & -1 & 0 \\
 0 & -3 & 2 & -1 \\
 -1 & 0 & -1 & 2 \\
}$
 & $\frac{1}{12}$ 
$\smat{
 0 & -9 & -6 & -3 \\
 -4 & -2 & -4 & -6 \\
 -8 & -7 & -2 & -9 \\
 -4 & -8 & -4 & 0 \\
}$
 & -12 } \\
 \stab{$NH4_{29}$&
$\smat{
 2 & -1 & 0 & -1 \\
 -1 & 2 & -1 & 0 \\
 0 & -3 & 2 & -1 \\
 -2 & 0 & -1 & 2 \\
}$
 & $\frac{1}{10}$ 
$\smat{
 0 & -6 & -4 & -2 \\
 -5 & -2 & -3 & -4 \\
 -10 & -8 & -2 & -6 \\
 -5 & -10 & -5 & 0 \\
}$
 & -10 } \\
 \stab{\attn$NH4_{30}$&
$\smat{
 2 & -2 & 0 & -1 \\
 -1 & 2 & -1 & 0 \\
 0 & -1 & 2 & -1 \\
 -1 & 0 & -1 & 2 \\
}$
 & $\frac{1}{4}$ 
$\smat{
 -4 & -7 & -6 & -5 \\
 -4 & -4 & -4 & -4 \\
 -4 & -5 & -2 & -3 \\
 -4 & -6 & -4 & -2 \\
}$
 & -4 } \\
 \stab{$NH4_{31}$&
$\smat{
 2 & -1 & 0 & -2 \\
 -2 & 2 & -1 & 0 \\
 0 & -1 & 2 & -1 \\
 -1 & 0 & -1 & 2 \\
}$
 & $\frac{1}{9}$ 
$\smat{
 -4 & -5 & -6 & -7 \\
 -7 & -2 & -6 & -10 \\
 -6 & -3 & 0 & -6 \\
 -5 & -4 & -3 & -2 \\
}$
 & -9 } \\
 \stab{$H4_{32}$&
$\smat{
 2 & -1 & 0 & -1 \\
 -2 & 2 & -1 & 0 \\
 0 & -1 & 2 & -1 \\
 -2 & 0 & -1 & 2 \\
}$
 & $\frac{1}{4}$ 
$\smat{
 -2 & -2 & -2 & -2 \\
 -4 & -1 & -2 & -3 \\
 -4 & -2 & 0 & -2 \\
 -4 & -3 & -2 & -1 \\
}$
 & -8 } \\
 \stab{\attn$H4_{33}$&
$\smat{
 2 & -2 & 0 & -1 \\
 -1 & 2 & -1 & 0 \\
 0 & -1 & 2 & -1 \\
 -1 & 0 & -2 & 2 \\
}$
 & $\frac{1}{7}$ 
$\smat{
 -2 & -6 & -8 & -5 \\
 -3 & -2 & -5 & -4 \\
 -4 & -5 & -2 & -3 \\
 -5 & -8 & -6 & -2 \\
}$
 & -7 } \\
 \stab{\attn$NH4_{34}$&
$\smat{
 2 & -1 & 0 & -1 \\
 -2 & 2 & -1 & 0 \\
 0 & -1 & 2 & -1 \\
 -1 & 0 & -2 & 2 \\
}$
 & $\frac{1}{8}$ 
$\smat{
 -2 & -4 & -6 & -4 \\
 -5 & -2 & -7 & -6 \\
 -6 & -4 & -2 & -4 \\
 -7 & -6 & -5 & -2 \\
}$
 & -8 } \\
 \stab{$NH4_{35}$&
$\smat{
 2 & -1 & 0 & -2 \\
 -2 & 2 & -1 & 0 \\
 0 & -1 & 2 & -1 \\
 -1 & 0 & -2 & 2 \\
}$
 & $\frac{1}{15}$ 
$\smat{
 -2 & -6 & -10 & -7 \\
 -5 & 0 & -10 & -10 \\
 -6 & -3 & 0 & -6 \\
 -7 & -6 & -5 & -2 \\
}$
 & -15 } \\
 \stab{$NH4_{36}$&
$\smat{
 2 & -1 & 0 & -1 \\
 -2 & 2 & -1 & 0 \\
 0 & -1 & 2 & -2 \\
 -2 & 0 & -1 & 2 \\
}$
 & $\frac{1}{12}$ 
$\smat{
 -2 & -3 & -4 & -5 \\
 -8 & 0 & -4 & -8 \\
 -12 & -6 & 0 & -6 \\
 -8 & -6 & -4 & -2 \\
}$
 & -12 } \\
 \stab{$NH4_{37}$&
$\smat{
 2 & -1 & 0 & -2 \\
 -2 & 2 & -1 & 0 \\
 0 & -2 & 2 & -1 \\
 -1 & 0 & -2 & 2 \\
}$
 & $\frac{1}{5}$ 
$\smat{
 0 & -2 & -2 & -1 \\
 -1 & 0 & -2 & -2 \\
 -2 & -1 & 0 & -2 \\
 -2 & -2 & -1 & 0 \\
}$
 & -25 } \\
 \stab{$H4_{38}$&
$\smat{
 2 & -1 & 0 & 0 \\
 -2 & 2 & -2 & -2 \\
 0 & -1 & 2 & 0 \\
 0 & -1 & 0 & 2 \\
}$
 & $\frac{1}{2}$ 
$\smat{
 0 & -1 & -1 & -1 \\
 -2 & -2 & -2 & -2 \\
 -1 & -1 & 0 & -1 \\
 -1 & -1 & -1 & 0 \\
}$
 & -8 } \\
 \stab{$H4_{39}$&
$\smat{
 2 & -1 & 0 & 0 \\
 -1 & 2 & -2 & -2 \\
 0 & -1 & 2 & 0 \\
 0 & -1 & 0 & 2 \\
}$
 & $\frac{1}{2}$ 
$\smat{
 0 & -2 & -2 & -2 \\
 -2 & -4 & -4 & -4 \\
 -1 & -2 & -1 & -2 \\
 -1 & -2 & -2 & -1 \\
}$
 & -4 } \\
 \stab{$H4_{40}$&
$\smat{
 2 & -2 & 0 & 0 \\
 -1 & 2 & -2 & -2 \\
 0 & -1 & 2 & 0 \\
 0 & -1 & 0 & 2 \\
}$
 & $\frac{1}{4}$ 
$\smat{
 0 & -4 & -4 & -4 \\
 -2 & -4 & -4 & -4 \\
 -1 & -2 & 0 & -2 \\
 -1 & -2 & -2 & 0 \\
}$
 & -8 } \\
 \stab{$H4_{41}$&
$\smat{
 2 & -1 & 0 & 0 \\
 -1 & 2 & -3 & -1 \\
 0 & -1 & 2 & 0 \\
 0 & -1 & 0 & 2 \\
}$
 & $\frac{1}{2}$ 
$\smat{
 0 & -2 & -3 & -1 \\
 -2 & -4 & -6 & -2 \\
 -1 & -2 & -2 & -1 \\
 -1 & -2 & -3 & 0 \\
}$
 & -4 } \\
 \stab{$H4_{42}$&
$\smat{
 2 & -1 & 0 & 0 \\
 -1 & 2 & -1 & -1 \\
 0 & -1 & 2 & 0 \\
 0 & -3 & 0 & 2 \\
}$
 & $\frac{1}{2}$ 
$\smat{
 0 & -2 & -1 & -1 \\
 -2 & -4 & -2 & -2 \\
 -1 & -2 & 0 & -1 \\
 -3 & -6 & -3 & -2 \\
}$
 & -4 } \\
 \stab{$H4_{43}$&
$\smat{
 2 & -1 & 0 & 0 \\
 -2 & 2 & -1 & -1 \\
 0 & -1 & 2 & 0 \\
 0 & -2 & 0 & 2 \\
}$
 & $\frac{1}{2}$ 
$\smat{
 -1 & -2 & -1 & -1 \\
 -4 & -4 & -2 & -2 \\
 -2 & -2 & 0 & -1 \\
 -4 & -4 & -2 & -1 \\
}$
 & -4 } \\
 \stab{$H4_{44}$&
$\smat{
 2 & -2 & 0 & 0 \\
 -1 & 2 & -1 & -1 \\
 0 & -1 & 2 & 0 \\
 0 & -2 & 0 & 2 \\
}$
 & $\frac{1}{2}$ 
$\smat{
 -1 & -4 & -2 & -2 \\
 -2 & -4 & -2 & -2 \\
 -1 & -2 & 0 & -1 \\
 -2 & -4 & -2 & -1 \\
}$
 & -4 } \\
 \stab{$H4_{45}$&
$\smat{
 2 & -2 & 0 & 0 \\
 -1 & 2 & -2 & -1 \\
 0 & -1 & 2 & 0 \\
 0 & -2 & 0 & 2 \\
}$
 & $\frac{1}{4}$ 
$\smat{
 0 & -4 & -4 & -2 \\
 -2 & -4 & -4 & -2 \\
 -1 & -2 & 0 & -1 \\
 -2 & -4 & -4 & 0 \\
}$
 & -8 } \\
 \stab{$H4_{46}$&
$\smat{
 2 & -2 & 0 & 0 \\
 -1 & 2 & -1 & -1 \\
 0 & -2 & 2 & 0 \\
 0 & -2 & 0 & 2 \\
}$
 & $\frac{1}{2}$ 
$\smat{
 0 & -2 & -1 & -1 \\
 -1 & -2 & -1 & -1 \\
 -1 & -2 & 0 & -1 \\
 -1 & -2 & -1 & 0 \\
}$
 & -8 } \\
 \stab{$H4_{47}$&
$\smat{
 2 & -1 & -1 & -1 \\
 -1 & 2 & -1 & -1 \\
 -1 & -1 & 2 & -1 \\
 -1 & -1 & -1 & 2 \\
}$
 & $\frac{1}{3}$ 
$\smat{
 0 & -1 & -1 & -1 \\
 -1 & 0 & -1 & -1 \\
 -1 & -1 & 0 & -1 \\
 -1 & -1 & -1 & 0 \\
}$
 & -27 } \\
 \stab{$H4_{48}$&
$\smat{
 2 & -1 & 0 & -1 \\
 -1 & 2 & -1 & -1 \\
 0 & -1 & 2 & -1 \\
 -1 & -1 & -1 & 2 \\
}$
 & $\frac{1}{6}$ 
$\smat{
 0 & -3 & -3 & -3 \\
 -3 & -2 & -3 & -4 \\
 -3 & -3 & 0 & -3 \\
 -3 & -4 & -3 & -2 \\
}$
 & -12 } \\
 \stab{$H4_{49}$&
$\smat{
 2 & -1 & 0 & 0 \\
 -3 & 2 & -1 & -1 \\
 0 & -1 & 2 & -1 \\
 0 & -1 & -1 & 2 \\
}$
 & $\frac{1}{3}$ 
$\smat{
 0 & -1 & -1 & -1 \\
 -3 & -2 & -2 & -2 \\
 -3 & -2 & 0 & -1 \\
 -3 & -2 & -1 & 0 \\
}$
 & -9 } \\
 \stab{$H4_{50}$&
$\smat{
 2 & -1 & 0 & 0 \\
 -1 & 2 & -1 & -1 \\
 0 & -1 & 2 & -1 \\
 0 & -1 & -1 & 2 \\
}$
 & $\frac{1}{3}$ 
$\smat{
 0 & -3 & -3 & -3 \\
 -3 & -6 & -6 & -6 \\
 -3 & -6 & -4 & -5 \\
 -3 & -6 & -5 & -4 \\
}$
 & -3 } \\
 \stab{$H4_{51}$&
$\smat{
 2 & -3 & 0 & 0 \\
 -1 & 2 & -1 & -1 \\
 0 & -1 & 2 & -1 \\
 0 & -1 & -1 & 2 \\
}$
 & $\frac{1}{3}$ 
$\smat{
 0 & -3 & -3 & -3 \\
 -1 & -2 & -2 & -2 \\
 -1 & -2 & 0 & -1 \\
 -1 & -2 & -1 & 0 \\
}$
 & -9 } \\
 \stab{$H4_{52}$&
$\smat{
 2 & -2 & 0 & 0 \\
 -1 & 2 & -1 & -1 \\
 0 & -1 & 2 & -1 \\
 0 & -1 & -1 & 2 \\
}$
 & $\frac{1}{6}$ 
$\smat{
 0 & -6 & -6 & -6 \\
 -3 & -6 & -6 & -6 \\
 -3 & -6 & -2 & -4 \\
 -3 & -6 & -4 & -2 \\
}$
 & -6 } \\
 \stab{$H4_{53}$&
$\smat{
 2 & -1 & 0 & 0 \\
 -2 & 2 & -1 & -1 \\
 0 & -1 & 2 & -1 \\
 0 & -1 & -1 & 2 \\
}$
 & $\frac{1}{6}$ 
$\smat{
 0 & -3 & -3 & -3 \\
 -6 & -6 & -6 & -6 \\
 -6 & -6 & -2 & -4 \\
 -6 & -6 & -4 & -2 \\
}$
 & -6 } \\
 \stab{$H5_{1}$&
$\smat{
 2 & -1 & 0 & 0 & 0 \\
 -2 & 2 & -1 & 0 & 0 \\
 0 & -1 & 2 & -1 & 0 \\
 0 & 0 & -2 & 2 & -1 \\
 0 & 0 & 0 & -1 & 2 \\
}$
 & $\frac{1}{2}$ 
$\smat{
 -1 & -2 & -3 & -2 & -1 \\
 -4 & -4 & -6 & -4 & -2 \\
 -6 & -6 & -6 & -4 & -2 \\
 -8 & -8 & -8 & -4 & -2 \\
 -4 & -4 & -4 & -2 & 0 \\
}$
 & -2 } \\
 \stab{$H5_{2}$&
$\smat{
 2 & -1 & 0 & 0 & 0 \\
 -2 & 2 & -1 & 0 & 0 \\
 0 & -1 & 2 & -2 & 0 \\
 0 & 0 & -1 & 2 & -1 \\
 0 & 0 & 0 & -1 & 2 \\
}$
 & $\frac{1}{2}$ 
$\smat{
 -1 & -2 & -3 & -4 & -2 \\
 -4 & -4 & -6 & -8 & -4 \\
 -6 & -6 & -6 & -8 & -4 \\
 -4 & -4 & -4 & -4 & -2 \\
 -2 & -2 & -2 & -2 & 0 \\
}$
 & -2 } \\
 \stab{$H5_{3}$&
$\smat{
 2 & -1 & 0 & 0 & 0 \\
 -1 & 2 & -2 & 0 & 0 \\
 0 & -1 & 2 & -1 & 0 \\
 0 & 0 & -1 & 2 & -1 \\
 0 & 0 & 0 & -2 & 2 \\
}$
 & $\frac{1}{2}$ 
$\smat{
 0 & -2 & -4 & -4 & -2 \\
 -2 & -4 & -8 & -8 & -4 \\
 -2 & -4 & -6 & -6 & -3 \\
 -2 & -4 & -6 & -4 & -2 \\
 -2 & -4 & -6 & -4 & -1 \\
}$
 & -2 } \\
 \stab{$H5_{4}$&
$\smat{
 2 & -1 & 0 & 0 & 0 \\
 -1 & 2 & -1 & 0 & 0 \\
 0 & -2 & 2 & -1 & 0 \\
 0 & 0 & -1 & 2 & -1 \\
 0 & 0 & 0 & -2 & 2 \\
}$
 & $\frac{1}{2}$ 
$\smat{
 0 & -2 & -2 & -2 & -1 \\
 -2 & -4 & -4 & -4 & -2 \\
 -4 & -8 & -6 & -6 & -3 \\
 -4 & -8 & -6 & -4 & -2 \\
 -4 & -8 & -6 & -4 & -1 \\
}$
 & -2 } \\
 \stab{\attn$NH5_{5}$&
$\smat{
 2 & -1 & 0 & 0 & -1 \\
 -1 & 2 & -1 & 0 & 0 \\
 0 & -1 & 2 & -2 & 0 \\
 0 & 0 & -1 & 2 & -1 \\
 -1 & 0 & 0 & -1 & 2 \\
}$
 & $\frac{1}{5}$ 
$\smat{
 -1 & -3 & -5 & -7 & -4 \\
 -4 & -2 & -5 & -8 & -6 \\
 -7 & -6 & -5 & -9 & -8 \\
 -5 & -5 & -5 & -5 & -5 \\
 -3 & -4 & -5 & -6 & -2 \\
}$
 & -5 } \\
 \stab{$H5_{6}$&
$\smat{
 2 & -1 & 0 & 0 & -1 \\
 -2 & 2 & -1 & 0 & 0 \\
 0 & -1 & 2 & -2 & 0 \\
 0 & 0 & -1 & 2 & -1 \\
 -1 & 0 & 0 & -1 & 2 \\
}$
 & $\frac{1}{8}$ 
$\smat{
 -1 & -3 & -5 & -7 & -4 \\
 -6 & -2 & -6 & -10 & -8 \\
 -10 & -6 & -2 & -6 & -8 \\
 -7 & -5 & -3 & -1 & -4 \\
 -4 & -4 & -4 & -4 & 0 \\
}$
 & -8 } \\
 \stab{$NH5_{7}$&
$\smat{
 2 & -2 & 0 & 0 & -1 \\
 -1 & 2 & -1 & 0 & 0 \\
 0 & -1 & 2 & -2 & 0 \\
 0 & 0 & -1 & 2 & -1 \\
 -1 & 0 & 0 & -1 & 2 \\
}$
 & $\frac{1}{9}$ 
$\smat{
 -1 & -5 & -8 & -11 & -6 \\
 -4 & -2 & -5 & -8 & -6 \\
 -7 & -8 & -2 & -5 & -6 \\
 -5 & -7 & -4 & -1 & -3 \\
 -3 & -6 & -6 & -6 & 0 \\
}$
 & -9 } \\
 \stab{$H5_{8}$&
$\smat{
 2 & -1 & 0 & 0 & -2 \\
 -1 & 2 & -1 & 0 & 0 \\
 0 & -1 & 2 & -2 & 0 \\
 0 & 0 & -1 & 2 & -1 \\
 -1 & 0 & 0 & -1 & 2 \\
}$
 & $\frac{1}{8}$ 
$\smat{
 -1 & -4 & -7 & -10 & -6 \\
 -4 & 0 & -4 & -8 & -8 \\
 -7 & -4 & -1 & -6 & -10 \\
 -5 & -4 & -3 & -2 & -6 \\
 -3 & -4 & -5 & -6 & -2 \\
}$
 & -8 } \\
 \stab{$H5_{9}$&
$\smat{
 2 & -1 & 0 & 0 & 0 \\
 -1 & 2 & -1 & -1 & -2 \\
 0 & -1 & 2 & 0 & 0 \\
 0 & -1 & 0 & 2 & 0 \\
 0 & -1 & 0 & 0 & 2 \\
}$
 & $\frac{1}{2}$ 
$\smat{
 0 & -2 & -1 & -1 & -2 \\
 -2 & -4 & -2 & -2 & -4 \\
 -1 & -2 & 0 & -1 & -2 \\
 -1 & -2 & -1 & 0 & -2 \\
 -1 & -2 & -1 & -1 & -1 \\
}$
 & -8 } \\
 \stab{$H5_{10}$&
$\smat{
 2 & -2 & 0 & 0 & 0 \\
 -1 & 2 & -1 & -1 & -1 \\
 0 & -1 & 2 & 0 & 0 \\
 0 & -1 & 0 & 2 & 0 \\
 0 & -1 & 0 & 0 & 2 \\
}$
 & $\frac{1}{2}$ 
$\smat{
 -1 & -4 & -2 & -2 & -2 \\
 -2 & -4 & -2 & -2 & -2 \\
 -1 & -2 & 0 & -1 & -1 \\
 -1 & -2 & -1 & 0 & -1 \\
 -1 & -2 & -1 & -1 & 0 \\
}$
 & -8 } \\
 \stab{$H5_{11}$&
$\smat{
 2 & -1 & 0 & -1 & 0 \\
 -1 & 2 & -1 & 0 & -1 \\
 0 & -1 & 2 & -1 & 0 \\
 -1 & 0 & -1 & 2 & -1 \\
 0 & -1 & 0 & -1 & 2 \\
}$
 & $\frac{1}{4}$ 
$\smat{
 0 & -2 & -2 & -2 & -2 \\
 -2 & -1 & -2 & -3 & -2 \\
 -2 & -2 & 0 & -2 & -2 \\
 -2 & -3 & -2 & -1 & -2 \\
 -2 & -2 & -2 & -2 & 0 \\
}$
 & -16 } \\
 \stab{$H5_{12}$&
$\smat{
 2 & -1 & 0 & 0 & 0 \\
 -1 & 2 & -1 & 0 & -1 \\
 0 & -1 & 2 & -1 & 0 \\
 0 & 0 & -1 & 2 & -1 \\
 0 & -1 & 0 & -1 & 2 \\
}$
 & $\frac{1}{4}$ 
$\smat{
 0 & -4 & -4 & -4 & -4 \\
 -4 & -8 & -8 & -8 & -8 \\
 -4 & -8 & -5 & -6 & -7 \\
 -4 & -8 & -6 & -4 & -6 \\
 -4 & -8 & -7 & -6 & -5 \\
}$
 & -4 } \\
 \stab{$H5_{13}$&
$\smat{
 2 & -1 & 0 & 0 & 0 \\
 -2 & 2 & -1 & 0 & -1 \\
 0 & -1 & 2 & -1 & 0 \\
 0 & 0 & -1 & 2 & -1 \\
 0 & -1 & 0 & -1 & 2 \\
}$
 & $\frac{1}{4}$ 
$\smat{
 0 & -2 & -2 & -2 & -2 \\
 -4 & -4 & -4 & -4 & -4 \\
 -4 & -4 & -1 & -2 & -3 \\
 -4 & -4 & -2 & 0 & -2 \\
 -4 & -4 & -3 & -2 & -1 \\
}$
 & -8 } \\
 \stab{$H5_{14}$&
$\smat{
 2 & -2 & 0 & 0 & 0 \\
 -1 & 2 & -1 & 0 & -1 \\
 0 & -1 & 2 & -1 & 0 \\
 0 & 0 & -1 & 2 & -1 \\
 0 & -1 & 0 & -1 & 2 \\
}$
 & $\frac{1}{4}$ 
$\smat{
 0 & -4 & -4 & -4 & -4 \\
 -2 & -4 & -4 & -4 & -4 \\
 -2 & -4 & -1 & -2 & -3 \\
 -2 & -4 & -2 & 0 & -2 \\
 -2 & -4 & -3 & -2 & -1 \\
}$
 & -8 } \\
 \stab{$H5_{15}$&
$\smat{
 2 & -1 & 0 & 0 & 0 \\
 -1 & 2 & -1 & -2 & 0 \\
 0 & -1 & 2 & 0 & 0 \\
 0 & -1 & 0 & 2 & -1 \\
 0 & 0 & 0 & -1 & 2 \\
}$
 & $\frac{1}{4}$ 
$\smat{
 -1 & -6 & -3 & -8 & -4 \\
 -6 & -12 & -6 & -16 & -8 \\
 -3 & -6 & -1 & -8 & -4 \\
 -4 & -8 & -4 & -8 & -4 \\
 -2 & -4 & -2 & -4 & 0 \\
}$
 & -4 } \\
 \stab{$H5_{16}$&
$\smat{
 2 & -1 & 0 & 0 & 0 \\
 -1 & 2 & -1 & -1 & 0 \\
 0 & -1 & 2 & 0 & 0 \\
 0 & -2 & 0 & 2 & -1 \\
 0 & 0 & 0 & -1 & 2 \\
}$
 & $\frac{1}{4}$ 
$\smat{
 -1 & -6 & -3 & -4 & -2 \\
 -6 & -12 & -6 & -8 & -4 \\
 -3 & -6 & -1 & -4 & -2 \\
 -8 & -16 & -8 & -8 & -4 \\
 -4 & -8 & -4 & -4 & 0 \\
}$
 & -4 } \\
 \stab{$H5_{17}$&
$\smat{
 2 & -1 & 0 & 0 & 0 \\
 -1 & 2 & -2 & -1 & 0 \\
 0 & -1 & 2 & 0 & 0 \\
 0 & -1 & 0 & 2 & -2 \\
 0 & 0 & 0 & -1 & 2 \\
}$
 & $\frac{1}{2}$ 
$\smat{
 0 & -2 & -2 & -2 & -2 \\
 -2 & -4 & -4 & -4 & -4 \\
 -1 & -2 & -1 & -2 & -2 \\
 -2 & -4 & -4 & -2 & -2 \\
 -1 & -2 & -2 & -1 & 0 \\
}$
 & -4 } \\
 \stab{$H5_{18}$&
$\smat{
 2 & -1 & 0 & 0 & 0 \\
 -1 & 2 & -2 & -1 & 0 \\
 0 & -1 & 2 & 0 & 0 \\
 0 & -1 & 0 & 2 & -1 \\
 0 & 0 & 0 & -1 & 2 \\
}$
 & $\frac{1}{2}$ 
$\smat{
 -2 & -6 & -6 & -4 & -2 \\
 -6 & -12 & -12 & -8 & -4 \\
 -3 & -6 & -5 & -4 & -2 \\
 -4 & -8 & -8 & -4 & -2 \\
 -2 & -4 & -4 & -2 & 0 \\
}$
 & -2 } \\
 \stab{$H5_{19}$&
$\smat{
 2 & -1 & 0 & 0 & 0 \\
 -1 & 2 & -2 & -1 & 0 \\
 0 & -1 & 2 & 0 & 0 \\
 0 & -1 & 0 & 2 & -1 \\
 0 & 0 & 0 & -2 & 2 \\
}$
 & $\frac{1}{2}$ 
$\smat{
 0 & -2 & -2 & -2 & -1 \\
 -2 & -4 & -4 & -4 & -2 \\
 -1 & -2 & -1 & -2 & -1 \\
 -2 & -4 & -4 & -2 & -1 \\
 -2 & -4 & -4 & -2 & 0 \\
}$
 & -4 } \\
 \stab{$H5_{20}$&
$\smat{
 2 & -1 & 0 & 0 & 0 \\
 -1 & 2 & -1 & -1 & 0 \\
 0 & -2 & 2 & 0 & 0 \\
 0 & -1 & 0 & 2 & -2 \\
 0 & 0 & 0 & -1 & 2 \\
}$
 & $\frac{1}{2}$ 
$\smat{
 0 & -2 & -1 & -2 & -2 \\
 -2 & -4 & -2 & -4 & -4 \\
 -2 & -4 & -1 & -4 & -4 \\
 -2 & -4 & -2 & -2 & -2 \\
 -1 & -2 & -1 & -1 & 0 \\
}$
 & -4 } \\
 \stab{$H5_{21}$&
$\smat{
 2 & -1 & 0 & 0 & 0 \\
 -1 & 2 & -1 & -1 & 0 \\
 0 & -2 & 2 & 0 & 0 \\
 0 & -1 & 0 & 2 & -1 \\
 0 & 0 & 0 & -1 & 2 \\
}$
 & $\frac{1}{2}$ 
$\smat{
 -2 & -6 & -3 & -4 & -2 \\
 -6 & -12 & -6 & -8 & -4 \\
 -6 & -12 & -5 & -8 & -4 \\
 -4 & -8 & -4 & -4 & -2 \\
 -2 & -4 & -2 & -2 & 0 \\
}$
 & -2 } \\
 \stab{$H5_{22}$&
$\smat{
 2 & -1 & 0 & 0 & 0 \\
 -1 & 2 & -1 & -1 & 0 \\
 0 & -2 & 2 & 0 & 0 \\
 0 & -1 & 0 & 2 & -1 \\
 0 & 0 & 0 & -2 & 2 \\
}$
 & $\frac{1}{2}$ 
$\smat{
 0 & -2 & -1 & -2 & -1 \\
 -2 & -4 & -2 & -4 & -2 \\
 -2 & -4 & -1 & -4 & -2 \\
 -2 & -4 & -2 & -2 & -1 \\
 -2 & -4 & -2 & -2 & 0 \\
}$
 & -4 } \\
 \stab{$H6_{1}$&
$\smat{
 2 & -1 & 0 & 0 & 0 & 0 \\
 -1 & 2 & -1 & -1 & 0 & 0 \\
 0 & -1 & 2 & 0 & 0 & 0 \\
 0 & -1 & 0 & 2 & -1 & 0 \\
 0 & 0 & 0 & -2 & 2 & -1 \\
 0 & 0 & 0 & 0 & -1 & 2 \\
}$
 & $\frac{1}{2}$ 
$\smat{
 0 & -2 & -1 & -3 & -2 & -1 \\
 -2 & -4 & -2 & -6 & -4 & -2 \\
 -1 & -2 & 0 & -3 & -2 & -1 \\
 -3 & -6 & -3 & -6 & -4 & -2 \\
 -4 & -8 & -4 & -8 & -4 & -2 \\
 -2 & -4 & -2 & -4 & -2 & 0 \\
}$
 & -4 } \\
 \stab{$H6_{2}$&
$\smat{
 2 & -1 & 0 & 0 & 0 & 0 \\
 -1 & 2 & -1 & -1 & 0 & 0 \\
 0 & -1 & 2 & 0 & 0 & 0 \\
 0 & -1 & 0 & 2 & -2 & 0 \\
 0 & 0 & 0 & -1 & 2 & -1 \\
 0 & 0 & 0 & 0 & -1 & 2 \\
}$
 & $\frac{1}{2}$ 
$\smat{
 0 & -2 & -1 & -3 & -4 & -2 \\
 -2 & -4 & -2 & -6 & -8 & -4 \\
 -1 & -2 & 0 & -3 & -4 & -2 \\
 -3 & -6 & -3 & -6 & -8 & -4 \\
 -2 & -4 & -2 & -4 & -4 & -2 \\
 -1 & -2 & -1 & -2 & -2 & 0 \\
}$
 & -4 } \\
 \stab{$NH6_{3}$&
$\smat{
 2 & -1 & 0 & 0 & 0 & -1 \\
 -1 & 2 & -1 & 0 & 0 & 0 \\
 0 & -1 & 2 & -1 & 0 & 0 \\
 0 & 0 & -1 & 2 & -2 & 0 \\
 0 & 0 & 0 & -1 & 2 & -1 \\
 -1 & 0 & 0 & 0 & -1 & 2 \\
}$
 & $\frac{1}{6}$ 
$\smat{
 0 & -2 & -4 & -6 & -8 & -4 \\
 -3 & 0 & -3 & -6 & -9 & -6 \\
 -6 & -4 & -2 & -6 & -10 & -8 \\
 -9 & -8 & -7 & -6 & -11 & -10 \\
 -6 & -6 & -6 & -6 & -6 & -6 \\
 -3 & -4 & -5 & -6 & -7 & -2 \\
}$
 & -6 } \\
 \stab{$H6_{4}$&
$\smat{
 2 & -1 & 0 & 0 & 0 & -1 \\
 -2 & 2 & -1 & 0 & 0 & 0 \\
 0 & -1 & 2 & -1 & 0 & 0 \\
 0 & 0 & -1 & 2 & -2 & 0 \\
 0 & 0 & 0 & -1 & 2 & -1 \\
 -1 & 0 & 0 & 0 & -1 & 2 \\
}$
 & $\frac{1}{4}$ 
$\smat{
 0 & -1 & -2 & -3 & -4 & -2 \\
 -2 & 0 & -2 & -4 & -6 & -4 \\
 -4 & -2 & 0 & -2 & -4 & -4 \\
 -6 & -4 & -2 & 0 & -2 & -4 \\
 -4 & -3 & -2 & -1 & 0 & -2 \\
 -2 & -2 & -2 & -2 & -2 & 0 \\
}$
 & -8 } \\
 \stab{$NH6_{5}$&
$\smat{
 2 & -2 & 0 & 0 & 0 & -1 \\
 -1 & 2 & -1 & 0 & 0 & 0 \\
 0 & -1 & 2 & -1 & 0 & 0 \\
 0 & 0 & -1 & 2 & -2 & 0 \\
 0 & 0 & 0 & -1 & 2 & -1 \\
 -1 & 0 & 0 & 0 & -1 & 2 \\
}$
 & $\frac{1}{3}$ 
$\smat{
 0 & -1 & -2 & -3 & -4 & -2 \\
 -1 & 0 & -1 & -2 & -3 & -2 \\
 -2 & -2 & 0 & -1 & -2 & -2 \\
 -3 & -4 & -2 & 0 & -1 & -2 \\
 -2 & -3 & -2 & -1 & 0 & -1 \\
 -1 & -2 & -2 & -2 & -2 & 0 \\
}$
 & -9 } \\
 \stab{$H6_{6}$&
$\smat{
 2 & -1 & 0 & 0 & 0 & 0 \\
 -1 & 2 & -1 & -1 & -1 & -1 \\
 0 & -1 & 2 & 0 & 0 & 0 \\
 0 & -1 & 0 & 2 & 0 & 0 \\
 0 & -1 & 0 & 0 & 2 & 0 \\
 0 & -1 & 0 & 0 & 0 & 2 \\
}$
 & $\frac{1}{2}$ 
$\smat{
 0 & -2 & -1 & -1 & -1 & -1 \\
 -2 & -4 & -2 & -2 & -2 & -2 \\
 -1 & -2 & 0 & -1 & -1 & -1 \\
 -1 & -2 & -1 & 0 & -1 & -1 \\
 -1 & -2 & -1 & -1 & 0 & -1 \\
 -1 & -2 & -1 & -1 & -1 & 0 \\
}$
 & -16 } \\
 \stab{$H6_{7}$&
$\smat{
 2 & -2 & 0 & 0 & 0 & 0 \\
 -1 & 2 & -1 & 0 & 0 & 0 \\
 0 & -1 & 2 & -1 & -1 & 0 \\
 0 & 0 & -1 & 2 & 0 & 0 \\
 0 & 0 & -1 & 0 & 2 & -1 \\
 0 & 0 & 0 & 0 & -1 & 2 \\
}$
 & $\frac{1}{2}$ 
$\smat{
 -4 & -10 & -12 & -6 & -8 & -4 \\
 -5 & -10 & -12 & -6 & -8 & -4 \\
 -6 & -12 & -12 & -6 & -8 & -4 \\
 -3 & -6 & -6 & -2 & -4 & -2 \\
 -4 & -8 & -8 & -4 & -4 & -2 \\
 -2 & -4 & -4 & -2 & -2 & 0 \\
}$
 & -2 } \\
 \stab{$H6_{8}$&
$\smat{
 2 & -1 & 0 & 0 & 0 & 0 \\
 -2 & 2 & -1 & 0 & 0 & 0 \\
 0 & -1 & 2 & -1 & -1 & 0 \\
 0 & 0 & -1 & 2 & 0 & 0 \\
 0 & 0 & -1 & 0 & 2 & -1 \\
 0 & 0 & 0 & 0 & -1 & 2 \\
}$
 & $\frac{1}{2}$ 
$\smat{
 -4 & -5 & -6 & -3 & -4 & -2 \\
 -10 & -10 & -12 & -6 & -8 & -4 \\
 -12 & -12 & -12 & -6 & -8 & -4 \\
 -6 & -6 & -6 & -2 & -4 & -2 \\
 -8 & -8 & -8 & -4 & -4 & -2 \\
 -4 & -4 & -4 & -2 & -2 & 0 \\
}$
 & -2 } \\
 \stab{$H6_{9}$&
$\smat{
 2 & -2 & 0 & 0 & 0 & 0 \\
 -1 & 2 & -1 & 0 & 0 & 0 \\
 0 & -1 & 2 & -1 & -1 & 0 \\
 0 & 0 & -1 & 2 & 0 & 0 \\
 0 & 0 & -1 & 0 & 2 & -1 \\
 0 & 0 & 0 & 0 & -2 & 2 \\
}$
 & $\frac{1}{2}$ 
$\smat{
 0 & -2 & -4 & -2 & -4 & -2 \\
 -1 & -2 & -4 & -2 & -4 & -2 \\
 -2 & -4 & -4 & -2 & -4 & -2 \\
 -1 & -2 & -2 & 0 & -2 & -1 \\
 -2 & -4 & -4 & -2 & -2 & -1 \\
 -2 & -4 & -4 & -2 & -2 & 0 \\
}$
 & -4 } \\
 \stab{$H6_{10}$&
$\smat{
 2 & -1 & 0 & 0 & 0 & 0 \\
 -2 & 2 & -1 & 0 & 0 & 0 \\
 0 & -1 & 2 & -1 & -1 & 0 \\
 0 & 0 & -1 & 2 & 0 & 0 \\
 0 & 0 & -1 & 0 & 2 & -1 \\
 0 & 0 & 0 & 0 & -2 & 2 \\
}$
 & $\frac{1}{2}$ 
$\smat{
 0 & -1 & -2 & -1 & -2 & -1 \\
 -2 & -2 & -4 & -2 & -4 & -2 \\
 -4 & -4 & -4 & -2 & -4 & -2 \\
 -2 & -2 & -2 & 0 & -2 & -1 \\
 -4 & -4 & -4 & -2 & -2 & -1 \\
 -4 & -4 & -4 & -2 & -2 & 0 \\
}$
 & -4 } \\
 \stab{$H6_{11}$&
$\smat{
 2 & -1 & 0 & 0 & 0 & 0 \\
 -2 & 2 & -1 & 0 & 0 & 0 \\
 0 & -1 & 2 & -1 & -1 & 0 \\
 0 & 0 & -1 & 2 & 0 & 0 \\
 0 & 0 & -1 & 0 & 2 & -2 \\
 0 & 0 & 0 & 0 & -1 & 2 \\
}$
 & $\frac{1}{2}$ 
$\smat{
 0 & -1 & -2 & -1 & -2 & -2 \\
 -2 & -2 & -4 & -2 & -4 & -4 \\
 -4 & -4 & -4 & -2 & -4 & -4 \\
 -2 & -2 & -2 & 0 & -2 & -2 \\
 -4 & -4 & -4 & -2 & -2 & -2 \\
 -2 & -2 & -2 & -1 & -1 & 0 \\
}$
 & -4 } \\
 \stab{$H6_{12}$&
$\smat{
 2 & -1 & 0 & 0 & 0 & 0 \\
 -2 & 2 & -1 & 0 & 0 & 0 \\
 0 & -1 & 2 & -1 & 0 & 0 \\
 0 & 0 & -1 & 2 & -1 & 0 \\
 0 & 0 & 0 & -2 & 2 & -1 \\
 0 & 0 & 0 & 0 & -1 & 2 \\
}$
 & $\frac{1}{2}$ 
$\smat{
 0 & -1 & -2 & -3 & -2 & -1 \\
 -2 & -2 & -4 & -6 & -4 & -2 \\
 -4 & -4 & -4 & -6 & -4 & -2 \\
 -6 & -6 & -6 & -6 & -4 & -2 \\
 -8 & -8 & -8 & -8 & -4 & -2 \\
 -4 & -4 & -4 & -4 & -2 & 0 \\
}$
 & -2 } \\
 \stab{$H6_{13}$&
$\smat{
 2 & -1 & 0 & 0 & 0 & 0 \\
 -2 & 2 & -1 & 0 & 0 & 0 \\
 0 & -1 & 2 & -1 & 0 & 0 \\
 0 & 0 & -1 & 2 & -2 & 0 \\
 0 & 0 & 0 & -1 & 2 & -1 \\
 0 & 0 & 0 & 0 & -1 & 2 \\
}$
 & $\frac{1}{2}$ 
$\smat{
 0 & -1 & -2 & -3 & -4 & -2 \\
 -2 & -2 & -4 & -6 & -8 & -4 \\
 -4 & -4 & -4 & -6 & -8 & -4 \\
 -6 & -6 & -6 & -6 & -8 & -4 \\
 -4 & -4 & -4 & -4 & -4 & -2 \\
 -2 & -2 & -2 & -2 & -2 & 0 \\
}$
 & -2 } \\
 \stab{$H6_{14}$&
$\smat{
 2 & -1 & 0 & 0 & 0 & 0 \\
 -1 & 2 & -1 & 0 & 0 & 0 \\
 0 & -1 & 2 & -1 & 0 & 0 \\
 0 & 0 & -1 & 2 & -1 & 0 \\
 0 & 0 & 0 & -2 & 2 & -1 \\
 0 & 0 & 0 & 0 & -1 & 2 \\
}$
 & 
$\smat{
 0 & -1 & -2 & -3 & -2 & -1 \\
 -1 & -2 & -4 & -6 & -4 & -2 \\
 -2 & -4 & -6 & -9 & -6 & -3 \\
 -3 & -6 & -9 & -12 & -8 & -4 \\
 -4 & -8 & -12 & -16 & -10 & -5 \\
 -2 & -4 & -6 & -8 & -5 & -2 \\
}$
 & -1 } \\
 \stab{$H6_{15}$&
$\smat{
 2 & -1 & 0 & 0 & 0 & 0 \\
 -1 & 2 & -1 & 0 & 0 & 0 \\
 0 & -2 & 2 & -1 & 0 & 0 \\
 0 & 0 & -1 & 2 & -1 & 0 \\
 0 & 0 & 0 & -1 & 2 & -1 \\
 0 & 0 & 0 & 0 & -1 & 2 \\
}$
 & 
$\smat{
 -2 & -5 & -4 & -3 & -2 & -1 \\
 -5 & -10 & -8 & -6 & -4 & -2 \\
 -8 & -16 & -12 & -9 & -6 & -3 \\
 -6 & -12 & -9 & -6 & -4 & -2 \\
 -4 & -8 & -6 & -4 & -2 & -1 \\
 -2 & -4 & -3 & -2 & -1 & 0 \\
}$
 & -1 } \\
 \stab{$H6_{16}$&
$\smat{
 2 & -1 & 0 & 0 & 0 & 0 \\
 -1 & 2 & -1 & 0 & 0 & 0 \\
 0 & -1 & 2 & -2 & 0 & 0 \\
 0 & 0 & -1 & 2 & -1 & 0 \\
 0 & 0 & 0 & -1 & 2 & -1 \\
 0 & 0 & 0 & 0 & -1 & 2 \\
}$
 & $\frac{1}{2}$ 
$\smat{
 0 & -2 & -4 & -6 & -4 & -2 \\
 -2 & -4 & -8 & -12 & -8 & -4 \\
 -4 & -8 & -12 & -18 & -12 & -6 \\
 -3 & -6 & -9 & -12 & -8 & -4 \\
 -2 & -4 & -6 & -8 & -4 & -2 \\
 -1 & -2 & -3 & -4 & -2 & 0 \\
}$
 & -2 } \\
 \stab{$H6_{17}$&
$\smat{
 2 & -1 & 0 & 0 & 0 & 0 \\
 -1 & 2 & -2 & 0 & 0 & 0 \\
 0 & -1 & 2 & -1 & 0 & 0 \\
 0 & 0 & -1 & 2 & -1 & 0 \\
 0 & 0 & 0 & -1 & 2 & -1 \\
 0 & 0 & 0 & 0 & -2 & 2 \\
}$
 & $\frac{1}{2}$ 
$\smat{
 0 & -2 & -4 & -4 & -4 & -2 \\
 -2 & -4 & -8 & -8 & -8 & -4 \\
 -2 & -4 & -6 & -6 & -6 & -3 \\
 -2 & -4 & -6 & -4 & -4 & -2 \\
 -2 & -4 & -6 & -4 & -2 & -1 \\
 -2 & -4 & -6 & -4 & -2 & 0 \\
}$
 & -2 } \\
 \stab{$H6_{18}$&
$\smat{
 2 & -1 & 0 & 0 & 0 & 0 \\
 -1 & 2 & -1 & 0 & 0 & 0 \\
 0 & -2 & 2 & -1 & 0 & 0 \\
 0 & 0 & -1 & 2 & -1 & 0 \\
 0 & 0 & 0 & -1 & 2 & -1 \\
 0 & 0 & 0 & 0 & -2 & 2 \\
}$
 & $\frac{1}{2}$ 
$\smat{
 0 & -2 & -2 & -2 & -2 & -1 \\
 -2 & -4 & -4 & -4 & -4 & -2 \\
 -4 & -8 & -6 & -6 & -6 & -3 \\
 -4 & -8 & -6 & -4 & -4 & -2 \\
 -4 & -8 & -6 & -4 & -2 & -1 \\
 -4 & -8 & -6 & -4 & -2 & 0 \\
}$
 & -2 } \\
 \stab{$H6_{19}$&
$\smat{
 2 & -1 & 0 & 0 & 0 & 0 \\
 -1 & 2 & -1 & 0 & 0 & -1 \\
 0 & -1 & 2 & -1 & 0 & 0 \\
 0 & 0 & -1 & 2 & -1 & 0 \\
 0 & 0 & 0 & -1 & 2 & -1 \\
 0 & -1 & 0 & 0 & -1 & 2 \\
}$
 & $\frac{1}{5}$ 
$\smat{
 0 & -5 & -5 & -5 & -5 & -5 \\
 -5 & -10 & -10 & -10 & -10 & -10 \\
 -5 & -10 & -6 & -7 & -8 & -9 \\
 -5 & -10 & -7 & -4 & -6 & -8 \\
 -5 & -10 & -8 & -6 & -4 & -7 \\
 -5 & -10 & -9 & -8 & -7 & -6 \\
}$
 & -5 } \\
 \stab{$H6_{20}$&
$\smat{
 2 & -1 & 0 & 0 & 0 & 0 \\
 -1 & 2 & -1 & -1 & -1 & 0 \\
 0 & -1 & 2 & 0 & 0 & 0 \\
 0 & -1 & 0 & 2 & 0 & 0 \\
 0 & -1 & 0 & 0 & 2 & -1 \\
 0 & 0 & 0 & 0 & -1 & 2 \\
}$
 & $\frac{1}{2}$ 
$\smat{
 -2 & -6 & -3 & -3 & -4 & -2 \\
 -6 & -12 & -6 & -6 & -8 & -4 \\
 -3 & -6 & -2 & -3 & -4 & -2 \\
 -3 & -6 & -3 & -2 & -4 & -2 \\
 -4 & -8 & -4 & -4 & -4 & -2 \\
 -2 & -4 & -2 & -2 & -2 & 0 \\
}$
 & -4 } \\
 \stab{$H6_{21}$&
$\smat{
 2 & -1 & 0 & 0 & 0 & 0 \\
 -1 & 2 & -1 & -1 & -1 & 0 \\
 0 & -1 & 2 & 0 & 0 & 0 \\
 0 & -1 & 0 & 2 & 0 & 0 \\
 0 & -1 & 0 & 0 & 2 & -1 \\
 0 & 0 & 0 & 0 & -2 & 2 \\
}$
 & $\frac{1}{2}$ 
$\smat{
 0 & -2 & -1 & -1 & -2 & -1 \\
 -2 & -4 & -2 & -2 & -4 & -2 \\
 -1 & -2 & 0 & -1 & -2 & -1 \\
 -1 & -2 & -1 & 0 & -2 & -1 \\
 -2 & -4 & -2 & -2 & -2 & -1 \\
 -2 & -4 & -2 & -2 & -2 & 0 \\
}$
 & -8 } \\
 \stab{$H6_{22}$&
$\smat{
 2 & -1 & 0 & 0 & 0 & 0 \\
 -1 & 2 & -1 & -1 & -1 & 0 \\
 0 & -1 & 2 & 0 & 0 & 0 \\
 0 & -1 & 0 & 2 & 0 & 0 \\
 0 & -1 & 0 & 0 & 2 & -2 \\
 0 & 0 & 0 & 0 & -1 & 2 \\
}$
 & $\frac{1}{2}$ 
$\smat{
 0 & -2 & -1 & -1 & -2 & -2 \\
 -2 & -4 & -2 & -2 & -4 & -4 \\
 -1 & -2 & 0 & -1 & -2 & -2 \\
 -1 & -2 & -1 & 0 & -2 & -2 \\
 -2 & -4 & -2 & -2 & -2 & -2 \\
 -1 & -2 & -1 & -1 & -1 & 0 \\
}$
 & -8 } \\
 \stab{$H7_1$&
$\smat{
 2 & -1 & 0 & 0 & 0 & 0 & 0 \\
 -1 & 2 & -1 & 0 & 0 & 0 & 0 \\
 0 & -1 & 2 & -1 & -1 & 0 & 0 \\
 0 & 0 & -1 & 2 & 0 & 0 & 0 \\
 0 & 0 & -1 & 0 & 2 & -1 & 0 \\
 0 & 0 & 0 & 0 & -1 & 2 & -1 \\
 0 & 0 & 0 & 0 & 0 & -2 & 2 \\
}$
 & $\frac{1}{2}$ 
$\smat{
 0 & -2 & -4 & -2 & -4 & -4 & -2 \\
 -2 & -4 & -8 & -4 & -8 & -8 & -4 \\
 -4 & -8 & -12 & -6 & -12 & -12 & -6 \\
 -2 & -4 & -6 & -2 & -6 & -6 & -3 \\
 -4 & -8 & -12 & -6 & -10 & -10 & -5 \\
 -4 & -8 & -12 & -6 & -10 & -8 & -4 \\
 -4 & -8 & -12 & -6 & -10 & -8 & -3 \\
}$
 & -2 } \\
 \stab{$H7_2$&
$\smat{
 2 & -1 & 0 & 0 & 0 & 0 & 0 \\
 -1 & 2 & -1 & 0 & 0 & 0 & 0 \\
 0 & -1 & 2 & -1 & -1 & 0 & 0 \\
 0 & 0 & -1 & 2 & 0 & 0 & 0 \\
 0 & 0 & -1 & 0 & 2 & -1 & 0 \\
 0 & 0 & 0 & 0 & -1 & 2 & -2 \\
 0 & 0 & 0 & 0 & 0 & -1 & 2 \\
}$
 & $\frac{1}{2}$ 
$\smat{
 0 & -2 & -4 & -2 & -4 & -4 & -4 \\
 -2 & -4 & -8 & -4 & -8 & -8 & -8 \\
 -4 & -8 & -12 & -6 & -12 & -12 & -12 \\
 -2 & -4 & -6 & -2 & -6 & -6 & -6 \\
 -4 & -8 & -12 & -6 & -10 & -10 & -10 \\
 -4 & -8 & -12 & -6 & -10 & -8 & -8 \\
 -2 & -4 & -6 & -3 & -5 & -4 & -3 \\
}$
 & -2 } \\
 \stab{$H7_3$&
$\smat{
 2 & -1 & 0 & 0 & 0 & 0 & 0 \\
 -1 & 2 & -1 & 0 & 0 & 0 & -1 \\
 0 & -1 & 2 & -1 & 0 & 0 & 0 \\
 0 & 0 & -1 & 2 & -1 & 0 & 0 \\
 0 & 0 & 0 & -1 & 2 & -1 & 0 \\
 0 & 0 & 0 & 0 & -1 & 2 & -1 \\
 0 & -1 & 0 & 0 & 0 & -1 & 2 \\
}$
 & $\frac{1}{6}$ 
$\smat{
 0 & -6 & -6 & -6 & -6 & -6 & -6 \\
 -6 & -12 & -12 & -12 & -12 & -12 & -12 \\
 -6 & -12 & -7 & -8 & -9 & -10 & -11 \\
 -6 & -12 & -8 & -4 & -6 & -8 & -10 \\
 -6 & -12 & -9 & -6 & -3 & -6 & -9 \\
 -6 & -12 & -10 & -8 & -6 & -4 & -8 \\
 -6 & -12 & -11 & -10 & -9 & -8 & -7 \\
}$
 & -6 } \\
 \stab{$H7_4$&
$\smat{
 2 & -1 & 0 & 0 & 0 & 0 & 0 \\
 -1 & 2 & -1 & 0 & -1 & 0 & 0 \\
 0 & -1 & 2 & -1 & 0 & 0 & 0 \\
 0 & 0 & -1 & 2 & 0 & 0 & 0 \\
 0 & -1 & 0 & 0 & 2 & -1 & -1 \\
 0 & 0 & 0 & 0 & -1 & 2 & 0 \\
 0 & 0 & 0 & 0 & -1 & 0 & 2 \\
}$
 & $\frac{1}{4}$ 
$\smat{
 -4 & -12 & -8 & -4 & -12 & -6 & -6 \\
 -12 & -24 & -16 & -8 & -24 & -12 & -12 \\
 -8 & -16 & -8 & -4 & -16 & -8 & -8 \\
 -4 & -8 & -4 & 0 & -8 & -4 & -4 \\
 -12 & -24 & -16 & -8 & -20 & -10 & -10 \\
 -6 & -12 & -8 & -4 & -10 & -3 & -5 \\
 -6 & -12 & -8 & -4 & -10 & -5 & -3 \\
}$
 & -4 } \\
 \stab{$H8_1$&
$\smat{
 2 & -1 & 0 & 0 & 0 & 0 & 0 & 0 \\
 -1 & 2 & -1 & 0 & 0 & 0 & 0 & -1 \\
 0 & -1 & 2 & -1 & 0 & 0 & 0 & 0 \\
 0 & 0 & -1 & 2 & -1 & 0 & 0 & 0 \\
 0 & 0 & 0 & -1 & 2 & -1 & 0 & 0 \\
 0 & 0 & 0 & 0 & -1 & 2 & -1 & 0 \\
 0 & 0 & 0 & 0 & 0 & -1 & 2 & -1 \\
 0 & -1 & 0 & 0 & 0 & 0 & -1 & 2 \\
}$
 & $\frac{1}{7}$ 
$\smat{
 0 & -7 & -7 & -7 & -7 & -7 & -7 & -7 \\
 -7 & -14 & -14 & -14 & -14 & -14 & -14 & -14 \\
 -7 & -14 & -8 & -9 & -10 & -11 & -12 & -13 \\
 -7 & -14 & -9 & -4 & -6 & -8 & -10 & -12 \\
 -7 & -14 & -10 & -6 & -2 & -5 & -8 & -11 \\
 -7 & -14 & -11 & -8 & -5 & -2 & -6 & -10 \\
 -7 & -14 & -12 & -10 & -8 & -6 & -4 & -9 \\
 -7 & -14 & -13 & -12 & -11 & -10 & -9 & -8 \\
}$
 & -7 } \\
 \stab{$H8_2$&
$\smat{
 2 & -1 & 0 & 0 & 0 & 0 & 0 & 0 \\
 -1 & 2 & -1 & 0 & 0 & 0 & 0 & 0 \\
 0 & -1 & 2 & -1 & -1 & 0 & 0 & 0 \\
 0 & 0 & -1 & 2 & 0 & 0 & 0 & 0 \\
 0 & 0 & -1 & 0 & 2 & -1 & 0 & 0 \\
 0 & 0 & 0 & 0 & -1 & 2 & -1 & 0 \\
 0 & 0 & 0 & 0 & 0 & -1 & 2 & -1 \\
 0 & 0 & 0 & 0 & 0 & 0 & -2 & 2 \\
}$
 & $\frac{1}{2}$ 
$\smat{
 0 & -2 & -4 & -2 & -4 & -4 & -4 & -2 \\
 -2 & -4 & -8 & -4 & -8 & -8 & -8 & -4 \\
 -4 & -8 & -12 & -6 & -12 & -12 & -12 & -6 \\
 -2 & -4 & -6 & -2 & -6 & -6 & -6 & -3 \\
 -4 & -8 & -12 & -6 & -10 & -10 & -10 & -5 \\
 -4 & -8 & -12 & -6 & -10 & -8 & -8 & -4 \\
 -4 & -8 & -12 & -6 & -10 & -8 & -6 & -3 \\
 -4 & -8 & -12 & -6 & -10 & -8 & -6 & -2 \\
}$
 & -2 } \\
 \stab{$H8_3$&
$\smat{
 2 & -1 & 0 & 0 & 0 & 0 & 0 & 0 \\
 -1 & 2 & -1 & 0 & 0 & 0 & 0 & 0 \\
 0 & -1 & 2 & -1 & -1 & 0 & 0 & 0 \\
 0 & 0 & -1 & 2 & 0 & 0 & 0 & 0 \\
 0 & 0 & -1 & 0 & 2 & -1 & 0 & 0 \\
 0 & 0 & 0 & 0 & -1 & 2 & -1 & 0 \\
 0 & 0 & 0 & 0 & 0 & -1 & 2 & -2 \\
 0 & 0 & 0 & 0 & 0 & 0 & -1 & 2 \\
}$
 & $\frac{1}{2}$ 
$\smat{
 0 & -2 & -4 & -2 & -4 & -4 & -4 & -4 \\
 -2 & -4 & -8 & -4 & -8 & -8 & -8 & -8 \\
 -4 & -8 & -12 & -6 & -12 & -12 & -12 & -12 \\
 -2 & -4 & -6 & -2 & -6 & -6 & -6 & -6 \\
 -4 & -8 & -12 & -6 & -10 & -10 & -10 & -10 \\
 -4 & -8 & -12 & -6 & -10 & -8 & -8 & -8 \\
 -4 & -8 & -12 & -6 & -10 & -8 & -6 & -6 \\
 -2 & -4 & -6 & -3 & -5 & -4 & -3 & -2 \\
}$
 & -2 } \\
 \stab{$H8_4$&
$\smat{
 2 & -1 & 0 & 0 & 0 & 0 & 0 & 0 \\
 -1 & 2 & -1 & 0 & 0 & 0 & 0 & 0 \\
 0 & -1 & 2 & -1 & -1 & 0 & 0 & 0 \\
 0 & 0 & -1 & 2 & 0 & 0 & 0 & 0 \\
 0 & 0 & -1 & 0 & 2 & -1 & 0 & 0 \\
 0 & 0 & 0 & 0 & -1 & 2 & -1 & -1 \\
 0 & 0 & 0 & 0 & 0 & -1 & 2 & 0 \\
 0 & 0 & 0 & 0 & 0 & -1 & 0 & 2 \\
}$
 & $\frac{1}{2}$ 
$\smat{
 0 & -2 & -4 & -2 & -4 & -4 & -2 & -2 \\
 -2 & -4 & -8 & -4 & -8 & -8 & -4 & -4 \\
 -4 & -8 & -12 & -6 & -12 & -12 & -6 & -6 \\
 -2 & -4 & -6 & -2 & -6 & -6 & -3 & -3 \\
 -4 & -8 & -12 & -6 & -10 & -10 & -5 & -5 \\
 -4 & -8 & -12 & -6 & -10 & -8 & -4 & -4 \\
 -2 & -4 & -6 & -3 & -5 & -4 & -1 & -2 \\
 -2 & -4 & -6 & -3 & -5 & -4 & -2 & -1 \\
}$
 & -4 } \\
 \stab{$H8_5$&
$\smat{
 2 & -1 & 0 & 0 & 0 & 0 & 0 & 0 \\
 -1 & 2 & -1 & 0 & 0 & 0 & 0 & 0 \\
 0 & -1 & 2 & -1 & 0 & 0 & -1 & 0 \\
 0 & 0 & -1 & 2 & -1 & 0 & 0 & 0 \\
 0 & 0 & 0 & -1 & 2 & -1 & 0 & 0 \\
 0 & 0 & 0 & 0 & -1 & 2 & 0 & 0 \\
 0 & 0 & -1 & 0 & 0 & 0 & 2 & -1 \\
 0 & 0 & 0 & 0 & 0 & 0 & -1 & 2 \\
}$
 & $\frac{1}{3}$ 
$\smat{
 -2 & -7 & -12 & -9 & -6 & -3 & -8 & -4 \\
 -7 & -14 & -24 & -18 & -12 & -6 & -16 & -8 \\
 -12 & -24 & -36 & -27 & -18 & -9 & -24 & -12 \\
 -9 & -18 & -27 & -18 & -12 & -6 & -18 & -9 \\
 -6 & -12 & -18 & -12 & -6 & -3 & -12 & -6 \\
 -3 & -6 & -9 & -6 & -3 & 0 & -6 & -3 \\
 -8 & -16 & -24 & -18 & -12 & -6 & -14 & -7 \\
 -4 & -8 & -12 & -9 & -6 & -3 & -7 & -2 \\
}$
 & -3 } \\
 \stab{$H9_1$&
$\smat{
 2 & -1 & 0 & 0 & 0 & 0 & 0 & 0 & 0 \\
 -1 & 2 & -1 & 0 & 0 & 0 & 0 & 0 & 0 \\
 0 & -1 & 2 & -1 & 0 & 0 & 0 & 0 & 0 \\
 0 & 0 & -1 & 2 & -1 & -1 & 0 & 0 & 0 \\
 0 & 0 & 0 & -1 & 2 & 0 & 0 & 0 & 0 \\
 0 & 0 & 0 & -1 & 0 & 2 & -1 & 0 & 0 \\
 0 & 0 & 0 & 0 & 0 & -1 & 2 & -1 & 0 \\
 0 & 0 & 0 & 0 & 0 & 0 & -1 & 2 & -1 \\
 0 & 0 & 0 & 0 & 0 & 0 & 0 & -1 & 2 \\
}$
 & $\frac{1}{2}$ 
$\smat{
 -1 & -4 & -7 & -10 & -5 & -8 & -6 & -4 & -2 \\
 -4 & -8 & -14 & -20 & -10 & -16 & -12 & -8 & -4 \\
 -7 & -14 & -21 & -30 & -15 & -24 & -18 & -12 & -6 \\
 -10 & -20 & -30 & -40 & -20 & -32 & -24 & -16 & -8 \\
 -5 & -10 & -15 & -20 & -9 & -16 & -12 & -8 & -4 \\
 -8 & -16 & -24 & -32 & -16 & -24 & -18 & -12 & -6 \\
 -6 & -12 & -18 & -24 & -12 & -18 & -12 & -8 & -4 \\
 -4 & -8 & -12 & -16 & -8 & -12 & -8 & -4 & -2 \\
 -2 & -4 & -6 & -8 & -4 & -6 & -4 & -2 & 0 \\
}$
 & -2 } \\
 \stab{$H9_2$&
$\smat{
 2 & -1 & 0 & 0 & 0 & 0 & 0 & 0 & 0 \\
 -1 & 2 & -1 & 0 & 0 & 0 & 0 & 0 & 0 \\
 0 & -1 & 2 & -1 & -1 & 0 & 0 & 0 & 0 \\
 0 & 0 & -1 & 2 & 0 & 0 & 0 & 0 & 0 \\
 0 & 0 & -1 & 0 & 2 & -1 & 0 & 0 & 0 \\
 0 & 0 & 0 & 0 & -1 & 2 & -1 & 0 & 0 \\
 0 & 0 & 0 & 0 & 0 & -1 & 2 & -1 & 0 \\
 0 & 0 & 0 & 0 & 0 & 0 & -1 & 2 & -1 \\
 0 & 0 & 0 & 0 & 0 & 0 & 0 & -2 & 2 \\
}$
 & $\frac{1}{2}$ 
$\smat{
 0 & -2 & -4 & -2 & -4 & -4 & -4 & -4 & -2 \\
 -2 & -4 & -8 & -4 & -8 & -8 & -8 & -8 & -4 \\
 -4 & -8 & -12 & -6 & -12 & -12 & -12 & -12 & -6 \\
 -2 & -4 & -6 & -2 & -6 & -6 & -6 & -6 & -3 \\
 -4 & -8 & -12 & -6 & -10 & -10 & -10 & -10 & -5 \\
 -4 & -8 & -12 & -6 & -10 & -8 & -8 & -8 & -4 \\
 -4 & -8 & -12 & -6 & -10 & -8 & -6 & -6 & -3 \\
 -4 & -8 & -12 & -6 & -10 & -8 & -6 & -4 & -2 \\
 -4 & -8 & -12 & -6 & -10 & -8 & -6 & -4 & -1 \\
}$
 & -2 } \\
 \stab{$H9_3$&
$\smat{
 2 & -1 & 0 & 0 & 0 & 0 & 0 & 0 & 0 \\
 -1 & 2 & -1 & 0 & 0 & 0 & 0 & 0 & 0 \\
 0 & -1 & 2 & -1 & -1 & 0 & 0 & 0 & 0 \\
 0 & 0 & -1 & 2 & 0 & 0 & 0 & 0 & 0 \\
 0 & 0 & -1 & 0 & 2 & -1 & 0 & 0 & 0 \\
 0 & 0 & 0 & 0 & -1 & 2 & -1 & 0 & 0 \\
 0 & 0 & 0 & 0 & 0 & -1 & 2 & -1 & 0 \\
 0 & 0 & 0 & 0 & 0 & 0 & -1 & 2 & -2 \\
 0 & 0 & 0 & 0 & 0 & 0 & 0 & -1 & 2 \\
}$
 & $\frac{1}{2}$ 
$\smat{
 0 & -2 & -4 & -2 & -4 & -4 & -4 & -4 & -4 \\
 -2 & -4 & -8 & -4 & -8 & -8 & -8 & -8 & -8 \\
 -4 & -8 & -12 & -6 & -12 & -12 & -12 & -12 & -12 \\
 -2 & -4 & -6 & -2 & -6 & -6 & -6 & -6 & -6 \\
 -4 & -8 & -12 & -6 & -10 & -10 & -10 & -10 & -10 \\
 -4 & -8 & -12 & -6 & -10 & -8 & -8 & -8 & -8 \\
 -4 & -8 & -12 & -6 & -10 & -8 & -6 & -6 & -6 \\
 -4 & -8 & -12 & -6 & -10 & -8 & -6 & -4 & -4 \\
 -2 & -4 & -6 & -3 & -5 & -4 & -3 & -2 & -1 \\
}$
 & -2 } \\
 \stab{$H9_4$&
$\smat{
 2 & -1 & 0 & 0 & 0 & 0 & 0 & 0 & 0 \\
 -1 & 2 & -1 & 0 & 0 & 0 & 0 & 0 & 0 \\
 0 & -1 & 2 & -1 & -1 & 0 & 0 & 0 & 0 \\
 0 & 0 & -1 & 2 & 0 & 0 & 0 & 0 & 0 \\
 0 & 0 & -1 & 0 & 2 & -1 & 0 & 0 & 0 \\
 0 & 0 & 0 & 0 & -1 & 2 & -1 & 0 & 0 \\
 0 & 0 & 0 & 0 & 0 & -1 & 2 & -1 & -1 \\
 0 & 0 & 0 & 0 & 0 & 0 & -1 & 2 & 0 \\
 0 & 0 & 0 & 0 & 0 & 0 & -1 & 0 & 2 \\
}$
 & $\frac{1}{4}$ 
$\smat{
 0 & -4 & -8 & -4 & -8 & -8 & -8 & -4 & -4 \\
 -4 & -8 & -16 & -8 & -16 & -16 & -16 & -8 & -8 \\
 -8 & -16 & -24 & -12 & -24 & -24 & -24 & -12 & -12 \\
 -4 & -8 & -12 & -4 & -12 & -12 & -12 & -6 & -6 \\
 -8 & -16 & -24 & -12 & -20 & -20 & -20 & -10 & -10 \\
 -8 & -16 & -24 & -12 & -20 & -16 & -16 & -8 & -8 \\
 -8 & -16 & -24 & -12 & -20 & -16 & -12 & -6 & -6 \\
 -4 & -8 & -12 & -6 & -10 & -8 & -6 & -1 & -3 \\
 -4 & -8 & -12 & -6 & -10 & -8 & -6 & -3 & -1 \\
}$
 & -4 } \\
 \stab{$H9_5$&
$\smat{
 2 & -1 & 0 & 0 & 0 & 0 & 0 & 0 & 0 \\
 -1 & 2 & -1 & 0 & 0 & 0 & 0 & 0 & -1 \\
 0 & -1 & 2 & -1 & 0 & 0 & 0 & 0 & 0 \\
 0 & 0 & -1 & 2 & -1 & 0 & 0 & 0 & 0 \\
 0 & 0 & 0 & -1 & 2 & -1 & 0 & 0 & 0 \\
 0 & 0 & 0 & 0 & -1 & 2 & -1 & 0 & 0 \\
 0 & 0 & 0 & 0 & 0 & -1 & 2 & -1 & 0 \\
 0 & 0 & 0 & 0 & 0 & 0 & -1 & 2 & -1 \\
 0 & -1 & 0 & 0 & 0 & 0 & 0 & -1 & 2 \\
}$
 & $\frac{1}{8}$ 
$\smat{
 0 & -8 & -8 & -8 & -8 & -8 & -8 & -8 & -8 \\
 -8 & -16 & -16 & -16 & -16 & -16 & -16 & -16 & -16 \\
 -8 & -16 & -9 & -10 & -11 & -12 & -13 & -14 & -15 \\
 -8 & -16 & -10 & -4 & -6 & -8 & -10 & -12 & -14 \\
 -8 & -16 & -11 & -6 & -1 & -4 & -7 & -10 & -13 \\
 -8 & -16 & -12 & -8 & -4 & 0 & -4 & -8 & -12 \\
 -8 & -16 & -13 & -10 & -7 & -4 & -1 & -6 & -11 \\
 -8 & -16 & -14 & -12 & -10 & -8 & -6 & -4 & -10 \\
 -8 & -16 & -15 & -14 & -13 & -12 & -11 & -10 & -9 \\
}$
 & -8 } \\
 \stab{$H10_1$&
$\smat{
 2 & -1 & 0 & 0 & 0 & 0 & 0 & 0 & 0 & 0 \\
 -1 & 2 & -1 & 0 & 0 & 0 & 0 & 0 & 0 & 0 \\
 0 & -1 & 2 & -1 & -1 & 0 & 0 & 0 & 0 & 0 \\
 0 & 0 & -1 & 2 & 0 & 0 & 0 & 0 & 0 & 0 \\
 0 & 0 & -1 & 0 & 2 & -1 & 0 & 0 & 0 & 0 \\
 0 & 0 & 0 & 0 & -1 & 2 & -1 & 0 & 0 & 0 \\
 0 & 0 & 0 & 0 & 0 & -1 & 2 & -1 & 0 & 0 \\
 0 & 0 & 0 & 0 & 0 & 0 & -1 & 2 & -1 & 0 \\
 0 & 0 & 0 & 0 & 0 & 0 & 0 & -1 & 2 & -1 \\
 0 & 0 & 0 & 0 & 0 & 0 & 0 & 0 & -1 & 2 \\
}$
 & 
$\smat{
 -4 & -9 & -14 & -7 & -12 & -10 & -8 & -6 & -4 & -2 \\
 -9 & -18 & -28 & -14 & -24 & -20 & -16 & -12 & -8 & -4 \\
 -14 & -28 & -42 & -21 & -36 & -30 & -24 & -18 & -12 & -6 \\
 -7 & -14 & -21 & -10 & -18 & -15 & -12 & -9 & -6 & -3 \\
 -12 & -24 & -36 & -18 & -30 & -25 & -20 & -15 & -10 & -5 \\
 -10 & -20 & -30 & -15 & -25 & -20 & -16 & -12 & -8 & -4 \\
 -8 & -16 & -24 & -12 & -20 & -16 & -12 & -9 & -6 & -3 \\
 -6 & -12 & -18 & -9 & -15 & -12 & -9 & -6 & -4 & -2 \\
 -4 & -8 & -12 & -6 & -10 & -8 & -6 & -4 & -2 & -1 \\
 -2 & -4 & -6 & -3 & -5 & -4 & -3 & -2 & -1 & 0 \\
}$
 & -1 } \\
 \stab{$H10_2$&
$\smat{
 2 & -1 & 0 & 0 & 0 & 0 & 0 & 0 & 0 & 0 \\
 -1 & 2 & -1 & 0 & 0 & 0 & 0 & 0 & 0 & 0 \\
 0 & -1 & 2 & -1 & -1 & 0 & 0 & 0 & 0 & 0 \\
 0 & 0 & -1 & 2 & 0 & 0 & 0 & 0 & 0 & 0 \\
 0 & 0 & -1 & 0 & 2 & -1 & 0 & 0 & 0 & 0 \\
 0 & 0 & 0 & 0 & -1 & 2 & -1 & 0 & 0 & 0 \\
 0 & 0 & 0 & 0 & 0 & -1 & 2 & -1 & 0 & 0 \\
 0 & 0 & 0 & 0 & 0 & 0 & -1 & 2 & -1 & 0 \\
 0 & 0 & 0 & 0 & 0 & 0 & 0 & -1 & 2 & -1 \\
 0 & 0 & 0 & 0 & 0 & 0 & 0 & 0 & -2 & 2 \\
}$
 & $\frac{1}{2}$ 
$\smat{
 0 & -2 & -4 & -2 & -4 & -4 & -4 & -4 & -4 & -2 \\
 -2 & -4 & -8 & -4 & -8 & -8 & -8 & -8 & -8 & -4 \\
 -4 & -8 & -12 & -6 & -12 & -12 & -12 & -12 & -12 & -6 \\
 -2 & -4 & -6 & -2 & -6 & -6 & -6 & -6 & -6 & -3 \\
 -4 & -8 & -12 & -6 & -10 & -10 & -10 & -10 & -10 & -5 \\
 -4 & -8 & -12 & -6 & -10 & -8 & -8 & -8 & -8 & -4 \\
 -4 & -8 & -12 & -6 & -10 & -8 & -6 & -6 & -6 & -3 \\
 -4 & -8 & -12 & -6 & -10 & -8 & -6 & -4 & -4 & -2 \\
 -4 & -8 & -12 & -6 & -10 & -8 & -6 & -4 & -2 & -1 \\
 -4 & -8 & -12 & -6 & -10 & -8 & -6 & -4 & -2 & 0 \\
}$
 & -2 } \\
 \stab{$H10_3$&
$\smat{
 2 & -1 & 0 & 0 & 0 & 0 & 0 & 0 & 0 & 0 \\
 -1 & 2 & -1 & 0 & 0 & 0 & 0 & 0 & 0 & 0 \\
 0 & -1 & 2 & -1 & -1 & 0 & 0 & 0 & 0 & 0 \\
 0 & 0 & -1 & 2 & 0 & 0 & 0 & 0 & 0 & 0 \\
 0 & 0 & -1 & 0 & 2 & -1 & 0 & 0 & 0 & 0 \\
 0 & 0 & 0 & 0 & -1 & 2 & -1 & 0 & 0 & 0 \\
 0 & 0 & 0 & 0 & 0 & -1 & 2 & -1 & 0 & 0 \\
 0 & 0 & 0 & 0 & 0 & 0 & -1 & 2 & -1 & 0 \\
 0 & 0 & 0 & 0 & 0 & 0 & 0 & -1 & 2 & -2 \\
 0 & 0 & 0 & 0 & 0 & 0 & 0 & 0 & -1 & 2 \\
}$
 & $\frac{1}{2}$ 
$\smat{
 0 & -2 & -4 & -2 & -4 & -4 & -4 & -4 & -4 & -4 \\
 -2 & -4 & -8 & -4 & -8 & -8 & -8 & -8 & -8 & -8 \\
 -4 & -8 & -12 & -6 & -12 & -12 & -12 & -12 & -12 & -12 \\
 -2 & -4 & -6 & -2 & -6 & -6 & -6 & -6 & -6 & -6 \\
 -4 & -8 & -12 & -6 & -10 & -10 & -10 & -10 & -10 & -10 \\
 -4 & -8 & -12 & -6 & -10 & -8 & -8 & -8 & -8 & -8 \\
 -4 & -8 & -12 & -6 & -10 & -8 & -6 & -6 & -6 & -6 \\
 -4 & -8 & -12 & -6 & -10 & -8 & -6 & -4 & -4 & -4 \\
 -4 & -8 & -12 & -6 & -10 & -8 & -6 & -4 & -2 & -2 \\
 -2 & -4 & -6 & -3 & -5 & -4 & -3 & -2 & -1 & 0 \\
}$
 & -2 } \\
 \stab{$H10_4$&
$\smat{
 2 & -1 & 0 & 0 & 0 & 0 & 0 & 0 & 0 & 0 \\
 -1 & 2 & -1 & 0 & 0 & 0 & 0 & 0 & 0 & 0 \\
 0 & -1 & 2 & -1 & -1 & 0 & 0 & 0 & 0 & 0 \\
 0 & 0 & -1 & 2 & 0 & 0 & 0 & 0 & 0 & 0 \\
 0 & 0 & -1 & 0 & 2 & -1 & 0 & 0 & 0 & 0 \\
 0 & 0 & 0 & 0 & -1 & 2 & -1 & 0 & 0 & 0 \\
 0 & 0 & 0 & 0 & 0 & -1 & 2 & -1 & 0 & 0 \\
 0 & 0 & 0 & 0 & 0 & 0 & -1 & 2 & -1 & -1 \\
 0 & 0 & 0 & 0 & 0 & 0 & 0 & -1 & 2 & 0 \\
 0 & 0 & 0 & 0 & 0 & 0 & 0 & -1 & 0 & 2 \\
}$
 & $\frac{1}{2}$ 
$\smat{
 0 & -2 & -4 & -2 & -4 & -4 & -4 & -4 & -2 & -2 \\
 -2 & -4 & -8 & -4 & -8 & -8 & -8 & -8 & -4 & -4 \\
 -4 & -8 & -12 & -6 & -12 & -12 & -12 & -12 & -6 & -6 \\
 -2 & -4 & -6 & -2 & -6 & -6 & -6 & -6 & -3 & -3 \\
 -4 & -8 & -12 & -6 & -10 & -10 & -10 & -10 & -5 & -5 \\
 -4 & -8 & -12 & -6 & -10 & -8 & -8 & -8 & -4 & -4 \\
 -4 & -8 & -12 & -6 & -10 & -8 & -6 & -6 & -3 & -3 \\
 -4 & -8 & -12 & -6 & -10 & -8 & -6 & -4 & -2 & -2 \\
 -2 & -4 & -6 & -3 & -5 & -4 & -3 & -2 & 0 & -1 \\
 -2 & -4 & -6 & -3 & -5 & -4 & -3 & -2 & -1 & 0 \\
}$
 & -4 } \\

\normalsize

\subsection{$NS3_{82}-NS3_{86}$}\label{ss8.2} Here, the values of $\alpha$ can only be equal to
\[
\alpha=\begin{cases}\frac{\sqrt{21}-11}{10},
\ \frac{\sqrt{15}-8}{7},
\ \frac{\sqrt{13}-7}{6}&\text{in cases $NS3_{82}-NS3_{84}$},\\
\frac{2\sqrt{6}-7}{5}, \ \frac{\sqrt{5}-3}{2}&\text{in cases $NS3_{85}-NS3_{86}$}.
\end{cases}
\]
Each of these five values of $\alpha$ is irrational and satisfies $-1<\alpha<0$.
Each of these 5 Lie superalgebras $\fg$ has 4 systems of simple roots, pairwise
non-equivalent, and here is what they are. Let $A_1$ be the Cartan matrix of the system
of simple roots all three of which are odd isotropic, i.e., the Cartan matrix whose main
diagonal is $(0,0,0)$; set
\be\label{Afour}
A_{k+1}:=r_k(A_1)\ \text{ for $k=1,2,3$},
\ee
where $r_k$ is the reflection \eqref{oddR} in the $k$th node. The Cartan matrices of the
four systems of simple roots of $\fg$ are precisely $A_1,\ldots, A_4$; the only isotropic
node of $A_{k+1}$ is the $k$th one, whereas all three nodes of $A_1$ are isotropic. These
four matrices are pairwise non-equivalent for $NS3_{82}-NS3_{85}$; for $NS3_{86}$, the
matrices $A_2$, $A_3$, $A_4$ are pairwise equivalent, see Statement~\ref{n86}.

The columns of the tables below are: the algebra and which of its four Cartan matrices
$A_i$ the row concerns, then $A_i$ itself, then $A_i^{-1}$, and the last column contains
$\det A_i$. Since $-1<\alpha<0$, we see that $\det A_1<0$, that $\det A_2=-\det A_1>0$
and that $\det A_3>0$; for the sign of $\det A_4$, which is negative, and for what the
four determinants have in common, see Subsection~\ref{Galois}.

We also observe that, unlike the inverses of Cartan matrices of hyperbolic Lie algebras,
each of the inverse matrices below has both positive and negative elements. As above, the
inverse matrix without zero elements is marked by \fbox{!!!}.
%\AO{DL: the discussion deleted here is the sentence that observed that the determinants of Subsections~\ref{ss8.1} and \ref{ssHyp} are all negative whereas some of the ones below are positive, and promised to indicate this in the table. Item E$)$ of Subsection~\ref{NP} states its first half, Corollary~\ref{GalCor} explains its second, and the annotations it announced are the ones you asked to remove.}
\[\arraycolsep=1.5pt
\begin{array}{@{}l@{\ }c@{\ \ }c@{\ \ }l@{}}
\attn NS3_{82}\text{--}NS3_{84},\ A_1&
\smat{
 0 & 1 & \alpha  \\
 -1 & 0 & 2+\alpha \\
 -1 & 2+\frac{1}{\alpha } & 0 \\
}&
-\frac{1}{3(\alpha+1)}\smat{
-2 \alpha -5-\frac{2}{\alpha }& 2 \alpha +1& \alpha +2 \\
-\alpha -2& \alpha & -\alpha\\
-2-\frac{1}{\alpha }& -1 & 1\\
}& -3(\alpha+1)\\[1mm]
\attn NS3_{82}\text{--}NS3_{84},\ A_2&
\smat{
 0 & 1 & \alpha  \\
 -1 & 2 & -1 \\
 -1 & -1 & 2 \\
}&
-\frac{1}{3 (\alpha +1)} \smat{
 -3 & \alpha +2 & 2 \alpha +1 \\
 -3 & -\alpha  & \alpha  \\
 -3 & 1 & -1 \\
}& 3(\alpha+1)\\[1mm]
\attn NS3_{82},\ A_3&
\smat{
 2 & -1 & -1 \\
 -1 & 0 & 2+\alpha \\
-2 & -1 & 2 \\
}&
\frac{1}{4 \alpha +5}\smat{
 \alpha +2 & 3 & -\alpha -2 \\
 -2 \alpha -2 & 2 & -2 \alpha -3 \\
 1 & 4 & -1 \\
}& 4\alpha+5\\[1mm]
\attn NS3_{82},\ A_4&
\smat{
 2 & -1 & -1 \\
 -2 & 2 & -1 \\
-1 &2+\frac{1}{\alpha} & 0 \\
}&
\frac{1}{5 \alpha +4} \smat{
 2 \alpha +1 & -2 \alpha -1 & 3 \alpha  \\
 \alpha  & -\alpha  & 4 \alpha  \\
 -2 (\alpha +1) & -3 \alpha -2 & 2 \alpha  \\
}& 5+\frac{4}{\alpha}\\[1mm]
\attn NS3_{83},\ A_3&
\smat{
 2 & -1 & -1 \\
 -1 & 0 & 2+\alpha \\
-3 & -1 & 2 \\
}&
\frac{1}{5 \alpha+7 } \smat{
 \alpha +2 & 3 & -\alpha -2 \\
 -3 \alpha -4 & 1 & -2 \alpha -3 \\
 1 & 5 & -1 \\
}& 5\alpha+7\\[1mm]
\attn NS3_{83},\ A_4&
\smat{
 2 & -1 & -1 \\
 -3 & 2 & -1 \\
 -1 &2+\frac{1}{\alpha} & 0 \\
}&
\frac{1}{7 \alpha +5} \smat{
 2 \alpha +1 & -2 \alpha -1 & 3 \alpha  \\
 \alpha  & -\alpha  & 5 \alpha  \\
 -4\alpha -3 & -3 \alpha -2 & \alpha  \\
}& 7+\frac{5}{\alpha}\\[1mm]
NS3_{84},\ A_3&
\smat{
 2 & -1 & -1 \\
 -1 & 0 & 2+\alpha \\
-4 & -1 & 2 \\
}&
\frac{1}{3 (2 \alpha +3)} \smat{
 \alpha +2 & 3 & -\alpha -2 \\
 -2 (2 \alpha +3) & 0 & -2 \alpha -3 \\
 1 & 6 & -1 \\
}& 3(2\alpha+3)\\[1mm]
NS3_{84},\ A_4&
\smat{
 2 & -1 & -1 \\
 -4 & 2 & -1 \\
-1 &2+\frac{1}{\alpha} & 0 \\
}&
\frac{1}{3(3 \alpha +2)} \smat{
 2 \alpha +1 & -2 \alpha -1 & 3 \alpha  \\
 \alpha  & -\alpha  & 6 \alpha  \\
 -2(3\alpha+2) & -3 \alpha -2 & 0  \\
}& 3\left(3+\frac{2}{\alpha}\right)\\
\end{array}
\]
\[\arraycolsep=1.5pt
\begin{array}{@{}l@{\ }c@{\ \ }c@{\ \ }l@{}}
\attn NS3_{85}\text{--}NS3_{86},\ A_1&
\smat{
 0 & 1 & \alpha  \\
 -1 & 0 & 3+\alpha \\
 -1 & 3+\frac{1}{\alpha } & 0 \\
}&
-\frac{1}{4(\alpha+1)}
\smat{
 -3 \alpha -10-\frac{3}{\alpha } & 3 \alpha +1& \alpha +3 \\
 -\alpha -3 & \alpha & -\alpha\\
 -3-\frac{1}{\alpha } & -1&1 \\
}& -4(\alpha+1)\\[1mm]
NS3_{85}\text{--}NS3_{86},\ A_2&
\smat{
 0 & 1 & \alpha  \\
 -1 & 2 & -2 \\
 -1 & -2 & 2 \\
}&
-\frac{1}{4 (\alpha +1)} \smat{
 0 & 2 (\alpha +1) & 2 (\alpha +1) \\
 -4 & -\alpha  & \alpha  \\
 -4 & 1 & -1 \\
}& 4(\alpha+1)\\[1mm]
\attn NS3_{85},\ A_3&
\smat{
 2 & -1 & -2 \\
 -1 & 0 & 3+\alpha \\
-1 & -1 & 2 \\
}&
\frac{1}{3 \alpha +5} \smat{
 \alpha +3 & 4 & -\alpha -3 \\
 -\alpha -1 & 2 & -2 \alpha -4 \\
 1 & 3 & -1 \\
}& 3\alpha+5\\[1mm]
\attn NS3_{85},\ A_4&
\smat{
 2 & -2 & -1 \\
 -1 & 2 & -1 \\
-1 &3+\frac{1}{\alpha} & 0 \\
}&
\frac{1}{5 \alpha +3} \smat{
 3 \alpha +1 & -3 \alpha -1 & 4 \alpha  \\
 \alpha  & -\alpha  & 3 \alpha  \\
 -\alpha -1 & -2 (2 \alpha +1) & 2 \alpha  \\
}& 5+\frac{3}{\alpha }\\[1mm]
NS3_{86},\ A_3&
\smat{
 2 & -1 & -2 \\
 -1 & 0 & 3+\alpha \\
-2 & -1 & 2 \\
}&
\frac{1}{4 (\alpha +2)} \smat{
 \alpha +3 & 4 & -\alpha -3 \\
 -2 (\alpha +2) & 0 & -2 (\alpha +2) \\
 1 & 4 & -1 \\
}& 4(\alpha+2)\\[1mm]
NS3_{86},\ A_4&
\smat{
 2 & -2 & -1 \\
 -2 & 2 & -1 \\
-1 &3+\frac{1}{\alpha} & 0 \\
}&
\frac{1}{4 (2 \alpha +1)} \smat{
 3 \alpha +1 & -3 \alpha -1 & 4 \alpha  \\
 \alpha  & -\alpha  & 4 \alpha  \\
 -2 (2 \alpha +1) & -2 (2 \alpha +1) & 0 \\
}& 8+\frac{4}{\alpha}\\
\end{array}
\]

\sssbegin{Remark}[The lower left entries of the matrices $A_4$]\label{items34}
The entries of the matrices $A_4$, and of the matrix $A_3$ of $NS3_{84}$, are the ones
that \eqref{Afour} yields, and this settles two of them that might be read otherwise. In
particular, the lower left entry of the matrix $A_4$ of $NS3_{83}$ is $-1$, and hence
$\det A_4=7+\nfrac5\alpha$ is the Galois conjugate of $\det A_3=5\alpha+7$, as
Subsection~\ref{Galois} requires; a matrix with $-\alpha$ in that corner is not a Cartan
matrix of $NS3_{83}$. Likewise, the matrix $A_4$ of $NS3_{84}$ has $-4$ in the lower left
corner: the matrix with $-2$ there is a Cartan matrix of $NS3_{82}$ alone, that entry
arising only for $\alpha=\nfrac{\sqrt{21}-11}{10}$, whereas $NS3_{84}$ has
$\alpha=\nfrac{\sqrt{13}-7}{6}$.
\end{Remark}

\subsection{The Galois involution, and the determinant when $A$ is not
symmetrizable}\label{Galois} None of the Cartan matrices of Subsection~\ref{ss8.2} is
symmetrizable. Hence, there is no invariant symmetric bilinear form on the span of the
roots, no signature, and formula \eqref{signdet} has nothing to act upon; the entries are
not rational either, so the index \eqref{PQ} is not available as it stands. What replaces
both is the arithmetic of the field the entries live in.

Each of the five values of $\alpha$ satisfies a quadratic equation whose constant term is
equal to its leading coefficient:
\be\label{alphaeq}
\begin{array}{lll}
5\alpha^2+11\alpha+5=0,& 7\alpha^2+16\alpha+7=0,& 3\alpha^2+7\alpha+3=0,\\[1mm]
5\alpha^2+14\alpha+5=0,& \alpha^2+3\alpha+1=0&
\end{array}
\ee
for $NS3_{82}-NS3_{86}$, respectively. Hence,
\be\label{unit}
\alpha\cdot\frac1\alpha=1\text{~~and~~}\alpha+\frac1\alpha\in\Qee,
\ee
i.e., \textbf{$\alpha$ is a unit of norm $1$ of the real quadratic field
$\Kee:=\Qee(\alpha)$, and $\nfrac1\alpha$ is precisely the Galois conjugate
$\bar\alpha$ of $\alpha$}. The entry $t+\nfrac1\alpha$ of the matrices $A_1$, where $t=2$
for $NS3_{82}-NS3_{84}$ and $t=3$ for $NS3_{85}, NS3_{86}$, is therefore the conjugate of
the entry $t+\alpha$ standing above it, and the Galois involution $\sigma$ of $\Kee$ over
$\Qee$ acts on the Cartan matrices themselves.

\sssbegin{Lemma}[The Galois involution permutes the Cartan matrices]\label{Gal} For each
of $NS3_{82}-NS3_{86}$, the involution $\sigma$ maps every Cartan matrix of $\fg$ to a
Cartan matrix of the same $\fg$. On $A_1$, it acts by the transposition of the $2$nd and
the $3$rd nodes followed by the rescaling of the first row by $\nfrac1\alpha$; it
interchanges $A_3$ and $A_4$, and there
\be\label{conj34}
\det A_4=\sigma(\det A_3).
\ee
\end{Lemma}

\begin{proof} Write $\sigma(A)$ for the matrix obtained from $A$ by conjugating every
entry, and let $s$ be the transposition of the $2$nd and the $3$rd indices. For
$A=A_1$, we have
\[
\sigma(A)=\smat{0&1&\nfrac1\alpha\\ -1&0&t+\nfrac1\alpha\\ -1&t+\alpha&0},\quad\text{so}\quad
(\sigma(A))_{s(i)s(j)}=\smat{0&\nfrac1\alpha&1\\ -1&0&t+\alpha\\ -1&t+\nfrac1\alpha&0},
\]
which is $A$ with its first row multiplied by $\nfrac1\alpha$. For $A=A_3$, the same
transposition turns $\sigma(A)$ into $A_4$ entry by entry; since a simultaneous
permutation of the rows and the columns does not change the determinant, \eqref{conj34}
follows. The matrix $A_2$ is symmetric in its $2$nd and $3$rd nodes, so $\sigma$ maps it
to itself up to the rescaling of the first row. Finally, every Cartan matrix of $\fg$ is
one of $A_1,\ldots, A_4$.\end{proof}

\sssbegin{Corollary}[The sign is not an invariant, the norm is]\label{GalCor} For each of
$NS3_{82}-NS3_{86}$, two of the four determinants are positive and two are negative, and
$\det A_3$ and $\det A_4$ multiply to a \textbf{negative rational} number,
namely to the norm
\be\label{norms}
N_{\Kee/\Qee}(\det A)=-3,\ -6,\ -9,\ -8,\ -16\text{~~for $NS3_{82}-NS3_{86}$,
respectively.}
\ee
\end{Corollary}

\begin{proof} Immediate from Lemma~\ref{Gal} and \eqref{unit}: for
$\det A_3=u\alpha+v$ with $u,v\in\Qee$, we have
\[
N_{\Kee/\Qee}(\det A_3)=(u\alpha+v)\left(u\nfrac1\alpha+v\right)=u^2+v^2+uv\left(\alpha+
\nfrac1\alpha\right),
\]
and the five values of $\alpha+\nfrac1\alpha$ are
$-\nfrac{11}{5}$, $-\nfrac{16}{7}$, $-\nfrac73$, $-\nfrac{14}{5}$, $-3$. The same
computation applied to $A_1$ and $A_2$, whose determinants differ by a sign, gives the
remaining pair.\end{proof}

So, for these five Lie superalgebras, \textbf{the sign of $\det A$ is not an invariant of
the algebra}: it is reversed by the Galois involution, which merely moves one Cartan
matrix of $\fg$ to another one of the same $\fg$. What is an invariant is the norm
$N_{\Kee/\Qee}(\det A)$, and \eqref{norms} says that this norm is negative in all five
cases. That negativity is what remains, in the non-symmetrizable case, of the Lorentzian
sign that \eqref{signdet} gives in the symmetrizable one: the determinant of a Cartan
matrix of an almost affine Lie superalgebra is negative --- if not itself, then in norm.

\sssbegin{Statement}[For $NS3_{86}$, two of the four matrices are equivalent]\label{n86}
Computed over $\Kee$ for each of the five algebras: for $NS3_{82}-NS3_{85}$, the four
Cartan matrices printed in Subsection~$\ref{ss8.2}$ are pairwise non-equivalent in the
sense of Subsection~$\ref{normA}$. For $NS3_{86}$, they are not$:$ the matrices $A_2$,
$A_3$, $A_4$ are pairwise equivalent, so only two equivalence classes remain, that of
$A_1$ and that of $A_2$.\end{Statement}

Indeed, for $NS3_{86}$ we have $3+\alpha=-\nfrac1\alpha$ and $3+\nfrac1\alpha=-\alpha$,
by the last of the equations \eqref{alphaeq}; this coincidence, which occurs for none of
the other four values of $\alpha$, makes $A_4$ into $A_2$ after the renumbering of the
nodes by $1\mapsto 3$, $2\mapsto 1$, $3\mapsto 2$ and the change of the sign of one row,
and $A_3$ into $A_2$ after the renumbering $1\mapsto2$, $2\mapsto3$, $3\mapsto1$ and the
rescaling of the row with the $0$ on the main diagonal by $-\nfrac1\alpha$, which is
allowed by the Warning~\ref{War}.

Statement~\ref{n86} concerns the Cartan \textit{matrices}, not the systems of simple
roots: the three reflections do produce four pairwise different systems of simple roots
for $NS3_{86}$ as well, but two of them happen to have equivalent Cartan matrices. Thus,
\textbf{the Cartan matrix does not always distinguish the systems of simple roots}.

%%%%%%%%%%%%%%%%%%%%%%%%%%%%%

\end{document}